%% file: main.tex
\documentclass{article}

\usepackage{mystyle}

\input{graphics-prelude.tex}

\title{Uniformly \texorpdfstring{polynomial\=/time}{polynomial-time} classification of surface homeomorphisms}
\author{Filippo Baroni}
\date{}

\begin{document}

\newconstant{g}{g}{0}{\text{genus}}
\newconstant{p}{p}{1}{\text{number of punctures}}
\newconstant{chi}{\chi}{2}{\text{Euler characteristic}}

\maketitle
\thispagestyle{fancy}

\begin{abstract}
We describe an algorithm which, given two essential curves on a surface $S$, computes their distance in the curve graph of $S$, up to multiplicative and additive errors.
As an application, we present an algorithm to decide the Nielsen\=/Thurston type (periodic, reducible, or pseudo\=/Anosov) of a mapping class of $S$.
The novelty of our algorithms lies in the fact that their running time is polynomial in the size of the input \emph{and} in the complexity of $S$ -- say, its Euler characteristic.
This is in contrast with previously known algorithms, which run in polynomial time in the size of the input \emph{for any fixed surface $S$}.
\end{abstract}

\include*{introduction}

\include*{curves-and-train-tracks}

\include*{distance-algorithm}

\include*{nielsen-thurston-classification}

\section*{Constants}
\label{sec:constants}
\addcontentsline{toc}{section}{Constants}

For ease of navigation, we have collected all the constants introduced in this article in the following table.

\constantstable{}

\printbibliography[heading=bibintoc]

\end{document}

%% file: graphics-prelude.tex
\NewExpandableDocumentCommand{\colourmakebackground}{m}{palette #1!20}
\NewExpandableDocumentCommand{\colourmakedark}{m}{palette #1!60!black}

\foreach \c[count=\i] in {8900FF,FF0000,00A3FF,FFA600,00D230} {
    \xglobal\definecolor{palette \i}{HTML}{\c}
}
\colorlet{grey}{gray}
\colorlet{grey background}{grey!10}

\usepackage[outline]{contour}
\contourlength{1pt}
\NewDocumentCommand{\tikzcontour}{m}{\contour{white}{#1}}

\tikzset{
line/.style={line width=.5pt,line cap=round},
line over/.style={line,preaction={draw=#1,line width=2pt,line cap=butt}},
pin line/.style={line cap=round,line width=.3,grey,postaction={decorate,decoration={markings,mark=at position 1 with {\fill[grey!80!black] circle (.5pt);}}}},
dot/.pic={\fill circle(1.5pt);},
label node/.style={inner sep=2pt,name=label}
}

\input{train-tracks.tex}
\input{one-vertex-triangulations.tex}
\input{annular-subsurface-projection.tex}


%% file: introduction.tex
\section{Introduction}

It was proved by \textcite{masur-minsky-1} that the curve graph $\C(S)$ of a surface $S$ is hyperbolic, and that a mapping class of $S$ is pseudo\=/Anosov if and only if it acts loxodromically on $\C(S)$.
In this article, we aim to exploit this characterisation to provide an efficient algorithm for recognising pseudo\=/Anosov surface homeomorphisms.
Naturally, the crucial building block of this classification algorithm will be a procedure to (coarsely) compute distances in the curve graph.

There is a wide collection of algorithms in the literature addressing the computation of distances in $\C(S)$ (for instance, \cite{watanabe-curve-complex,shackleton-tightness,birman-efficient-geodesics}) and the Nielsen\=/Thurston classification of surface homeomorphisms (for instance, \cite{bestvina-train-tracks,bell-mapping-classes,koberda-classification,takarajima-train-tracks,hamidi-tehrani-diffeomorphisms}).
These lists include \cite{bell-webb-algorithms}, where \citeauthor{bell-webb-algorithms} describe, for each \emph{fixed} surface, polynomial\=/time algorithms to solve the two mentioned problems.
In contrast, the main contribution of the present article is the focus on providing algorithms that run in polynomial time in the size of the input, \emph{uniformly with respect to the complexity of the surface} -- by which we mean, for example, its Euler characteristic.

In order to estimate the distance between two curves in $\C(S)$, we rely on a result of \textcite{masur-minsky-quasiconvexity}, stating that a splitting sequence of train tracks corresponds to an unparametrised quasi\=/geodesic in the curve graph.
By applying a combinatorial criterion of \textcite{masur-minsky-1} to a carefully constructed splitting sequence, we are able to find the correct quasi\=/geodesic parametrisation, thus proving the following.

\newconstant{carried curves l+}{L_+^0}{L+0}{\constant{one switch d}+2\constant{carried curves d}+8}
\newconstant{carried curves l*}{L_\times}{L*}{2\constant{carried curves d}+12}
\newconstant{carried curves d}{D}{D}{\constant{deep nesting d}+\constant{proximity d}\cdot\constant{exponential decay n}}
\newconstant{curves l+}{L_+}{L+}{\constant{carried curves l+}+36}
\begin{theorem*}
Let $a$ and $b$ be essential curves on $S$. There is an algorithm to compute an estimate $d$ of the distance $\dist[\C(S)](a,b)$ between $a$ and $b$ in $\C(S)$, such that
\[
d-\constant{curves l+}\le\dist[\C(S)](a,b)\le\constant{carried curves l*}\cdot d+\constant{curves l+}.
\]
The constants $\constant{carried curves l*}$ and $\constant{curves l+}$ are polynomial in the complexity of $S$, and the running time of the algorithm is polynomial in the complexities of $S$, $a$, and $b$.
\end{theorem*}

We will describe in detail what we mean by complexity of curves in \cref{sec:surfaces curves train tracks:surfaces}; for now, let us say that the surface $S$ is given as input in the form of a triangulation, and the complexity of a curve is logarithmic in the number of times it intersects the edges of the triangulation.
Note that the algorithm in the theorem above does not compute the exact distance between two curves, but rather a reasonably accurate estimate thereof.
However, this approximation is enough to detect the Nielsen\=/Thurston type of a surface homeomorphism $f$.
In fact, for an arbitrary curve $a$ on $S$, the quantity $\dist[\C(S)](a,f^k(a))$ grows linearly in $k$ if $f$ is pseudo\=/Anosov, and remains bounded if $f$ is periodic or reducible.
This statement can be made effective, thanks to the work of \textcite{gadre-min-translation-length} and \textcite{webb-effective-geometry}, and leads to a proof of the following theorem.

\begin{theorem*}
Let $f$ be a homeomorphism of $S$. There is an algorithm to decide if $f$ is periodic, reducible (and non\=/periodic), or pseudo\=/Anosov. The running time of the algorithm is polynomial in the complexities of $S$ and $f$.
\end{theorem*}

Again, we will explain our notion of complexity of homeomorphisms in \cref{sec:surfaces curves train tracks:surfaces}; for now, we can take it to mean the word length in the mapping class group of $S$ with respect to a generating set consisting of ``reasonably short'' Dehn twists.

Let us emphasise, once more, that the novelty of this algorithm resides in the good dependence of the running time on the complexity of the surface.
To better illustrate this point, let us consider the braid group $B_n$ on $n$ strands with the standard generating set $\{\sigma_1,\ldots,\sigma_{n-1}\}$.
The following result, which is almost a direct generalisation of \citeauthor{calvarez-braids}'s \cite[Theorem 1]{calvarez-braids}, is an easy consequence of the classification theorem above.

\begin{corollary*}
There is an algorithm to decide the Nielsen\=/Thurston type of a braid in polynomial time in its word length and in the number $n$ of strands.
\end{corollary*}

This article is organised as follows.
In \cref{sec:surfaces curves train tracks} we set up the definitions and notation we will use throughout the article.
Most of this section is dedicated to combinatorial properties of train tracks, which will be crucial for our coarse distance algorithm.
We also discuss in detail how to represent surfaces, curves, homeomorphisms, and train tracks in a discrete fashion which is suitable to algorithmic manipulations, as well as the associated notions of complexity.
\Cref{sec:distance algorithm} follows, with a description of the algorithm to estimate distances in the curve graph.
The bulk of this section is devoted to the construction of specific splitting sequences of train tracks and their quasi\=/geodesic subsequences.
Finally, in \cref{sec:classification}, we apply the distance algorithm to the Nielsen\=/Thurston classification problem.

We have made an effort to keep the constants in our algorithms as explicit as possible.
For the sake of convenience, we have collected all of them in a table, which can be found \hyperref[sec:constants]{at the end of the article}.

%% file: curves-and-train-tracks.tex
\section{Surfaces, curves, and train tracks}
\label{sec:surfaces curves train tracks}

In this section, we introduce definitions and properties of the objects we will be using in our algorithms, including surfaces, curves, homeomorphisms, and train tracks.
These objects need to be fed as input to, manipulated by, and returned as output by algorithmic routines.
For this reason, alongside the usual, ``theoretical'' definitions, we will also provide ways to represent these objects combinatorially, and define complexities of these representations.

As anticipated before, in this article we focus on good dependence of the running times on the complexity of the underlying surface.
Therefore, all the constants we hide in the notation will be \emph{universal}.
For instance, for real\=/valued functions $F$ and $G$, we will say that ``$F$ is $\bigO(G)$'' if $F\le \alpha\cdot G+\beta$ for some universal constants $\alpha,\beta\ge 0$.
Similarly, we will write ``$\poly(X)$'' to signify a universal, albeit implicit, polynomial in the variable $X$; by ``$Y$ is polynomial in $X_1,\ldots,X_n$'', we will mean that
\[
    Y=\bigO(\poly(X_1)\cdots\poly(X_n)).
\]

\subsection{Triangulated surfaces}
\label{sec:surfaces curves train tracks:surfaces}

\step{Basic definitions.} In the context of this article, by ``surface'' we always mean a connected, orientable surface $S$ of finite genus $g\ge 0$, with $p\ge 0$ punctures and empty boundary. We will say that $S$ is \emph{closed} if $p=0$, and \emph{punctured} otherwise.

By ``curve on $S$'' we always mean a simple closed curve embedded in $S$. A curve $a\subs S$ is \emph{essential} if it does not bound a disc or a once\=/punctured disc in $S$. Given two curves $a$ and $b$ on $S$, we denote their (geometric) intersection number by $i(a,b)$. A \emph{multicurve} on $S$ is a disjoint union of finitely many curves on $S$, and it is \emph{essential} if all its components are.

By ``homeomorphism of $S$'' we always mean an orientation\=/preserving homeomorphism $\umap{S}{S}$. The \emph{mapping class group} of $S$, denoted by $\MCG(S)$, is the group of homeomorphisms of $S$ modulo isotopies.

\let\oldxi\xi
\newconstant{xi}{\oldxi}{4}{3g-3+p}
\gdef\xi{\constant{xi}}
\step{Complexity of surfaces.}
In this article, we aim to provide algorithms whose running times depend polynomially on the complexity of a surface. There are two essentially equivalent measures of complexity of a surface $S$: its Euler characteristic $\chi=2-2g-p$, and the integer
\[
    \declareconstant{xi}=\constantvalue{xi},
\]
which counts, for instance, the maximum number of disjoint non\=/parallel essential curves on $S$ (assuming that $\xi\ge 1$). For the sake of consistency, we will employ $\xi$ as our notion of complexity of $S$, although we remark that, for sufficiently complex surfaces, the two integers $\xi$ and $-\chi$ are equivalent up to a multiplicative constant.

For surfaces with $\xi\le 1$ -- namely, the sphere with at most four punctures and the torus with at most one puncture -- the problems we are trying to solve are already well\=/understood.
In fact, for these surfaces, the mapping class group and the curve graph are either essentially trivial (for the sphere with at most three punctures), or easily describable in terms of $\operatorname{SL(2,\ZZ)}$ and the Farey graph (see \cite[\S2.2]{farb-margalit} and \cite{beardon-continued-fractions}).
Therefore, unless otherwise stated, we will work with the underlying assumption that our ambient surface (always denoted by $S$) has complexity $\xi\ge 2$.

\step{Triangulations.}
To provide a surface as input to an algorithm, it is necessary to represent it combinatorially, and triangulations frequently offer the most suitable approach for this purpose.
We begin by describing the pertinent concepts for punctured surfaces, while deferring the discussion of closed surfaces to a subsequent paragraph.

Let us think of the punctures of a punctured surface $S$ as marked points on a closed surface. An \emph{ideal triangulation} (or simply a \emph{triangulation}) of $S$ is an ordered set $\TTT=\{e_1,\ldots,e_n\}$ where each $e_i$ is an arc embedded in $S$, such that
\begin{enumroman}
\item the endpoints of each arc are punctures,
\item the interiors of any two arcs are disjoint, and
\item each component $Z$ of $S\setminus(e_1\cup\ldots\cup e_n)$ is a \emph{triangle}, in the sense that there is a map $\umap{\Delta}{S}$ from a Euclidean $2$\=/simplex $\Delta$ to $S$ which sends each vertex of $\Delta$ to a puncture of $S$, the interior of each edge of $\Delta$ homeomorphically to the interior of an arc of $\TTT$, and the interior of $\Delta$ homeomorphically to $Z$, and moreover $Z$ does not contain any punctures.
\end{enumroman}

The arcs $e_1,\ldots,e_n$ are called \emph{edges} of $\TTT$; an easy counting argument shows that there are exactly $n=-3\chi$ edges in any ideal triangulation of $S$.
We will denote by $\card{\TTT}$ the number of edges of $\TTT$.

Note that a triangulation can be described combinatorially by providing a set of punctures, an ordered list of edges -- along with the information of which punctures are the endpoints of each edge -- and a set of triangles -- along with the information of the sequence of edges forming the boundary of each triangle.
We will always assume that surfaces are given as inputs to our algorithms in this form.

When $S$ is a closed surface, we can introduce a fictitious marked point and consider triangulations having this point as their only vertex. We call them \emph{one\=/vertex triangulations} (or simply \emph{triangulations}) of $S$; this is the combinatorial representation of closed surfaces that our algorithms will accept as input.

\step{Normal curves.}
Continuing our endeavour to establish combinatorial representations of the topological objects we wish to manipulate, we now shift our focus to curves on the surface $S$.

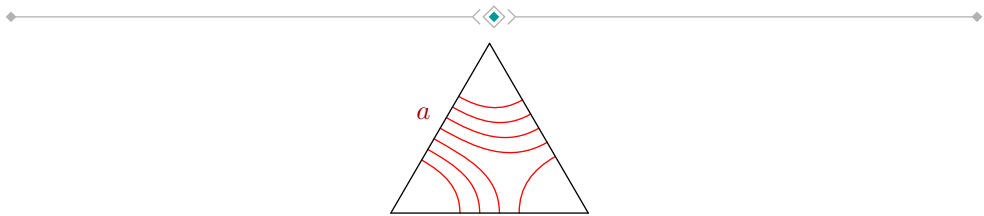
\begin{figure}
\centering
    \tikzsetnextfilename{normal-curve-example}%
    \input{figures-source/normal-curve-example.tex}%

\caption{The typical intersection between a normal curve $a$ and a triangle of the triangulation $\TTT$.}
\label{fig:normal curve example}
\end{figure}

Let us fix a triangulation $\TTT$ of $S$, one\=/vertex or ideal depending on whether $S$ is closed or not.
If $a$ is an essential curve (or, in fact, an essential multicurve) on $S$, then it can be isotoped so that its intersection with each triangle $\Delta$ of $\TTT$ is a collection of arcs, each of which connects distinct edges of $\Delta$; we call such a (multi)curve \emph{normal}.
The typical intersection of a normal curve with a triangle of $\TTT$ is illustrated in \cref{fig:normal curve example}.

A normal curve $a$ can be efficiently encoded by recording how many times it intersects each edge $e$ of $\TTT$; we denote this number by $\TTT(a)_e$.
Recalling that edges of $\TTT$ are ordered, we can collect the non\=/negative integers $\{\TTT(a)_e:e\in\TTT\}$ into a vector $\TTT(a)\in\ZZ_{\ge 0}^n$, where $n=\card{\TTT}$ is the number of edges of $\TTT$.

This representation is unambiguous, in the sense that we can recover the curve $a$ up to isotopy from the vector $\TTT(a)$.
It is also efficient, in the sense that it can be used to represent exponentially complicated curves with a linear amount of information.
More precisely, we define the \emph{complexity} of a normal curve $a$ with respect to the triangulation $\TTT$ to be the number
\[
    \norm[\TTT]{a}=\sum_{e\in\TTT}\log(\TTT(a)_e+1),
\]
where here -- and everywhere else -- we take the base of the logarithm to be $2$. We will omit the subscript $\TTT$ if the triangulation can be inferred from the context.

Note that the complexity of $a$ is roughly proportional to the number of digits that are required to store the vector $\TTT(a)$, but is clearly logarithmic in the number of intersections between $a$ and the edges of $\TTT$.
When we say that an algorithm takes as input an essential curve $a$ on a surface $S$, we will actually mean that it takes the vector $\TTT(a)$, where $\TTT$ is the triangulation representing $S$.
Likewise, all actions performed by our algorithms on essential curves will only operate on this vector representation.

Sometimes, it is convenient to use an alternative notion of complexity of normal curves, defined by
\[
\norm[\TTT,1]{a}=\log\Bigl(\sum_{e\in\TTT}\TTT(a)_e\Bigr)
\]
for a normal curve $a$; like above, we will omit the subscript $\TTT$ unless needed. Note that the two complexities are related by the inequalities
\[
\norm{a}_1\le\norm{a}\le\card{\TTT}\cdot\norm{a}_1.
\]

We conclude by remarking that, when $S$ is punctured, the vector $\TTT(a)$ is invariant under isotopy of $a$, in the sense that if $a$ and $b$ are normal curves on $S$ and are isotopic, then $\TTT(a)=\TTT(b)$.
This property immediately gives a linear\=/time algorithm in $\norm{a}+\norm{b}$ to decide whether two curves $a$ and $b$ are isotopic.
The situation is not as straightforward when $S$ is closed -- essentially because we allow isotopies across the vertex of $\TTT$ -- but a result of \textcite[Theorem 1.2]{lackenby-minimal-position} gives us the following proposition.

\begin{proposition}\label{thm:deciding isotopy of curves}
Let $a$ and $b$ be essential curves on a surface $S$. There is an algorithm to decide if $a$ and $b$ are isotopic. The running time of the algorithm is polynomial in $\xi$, $\norm{a}$, and $\norm{b}$.
\end{proposition}

\step{Flip and twist graph.}
A convenient framework to describe the relations between different triangulations of a surface is given by the flip and twist graph, introduced for algorithmic purposes by \textcite{bell-simplifying-triangulations}.
Suppose first that the surface $S$ is punctured.
The \emph{flip and twist graph} of $S$ is a graph $\GGG(S)$ whose vertices are  isotopy classes of triangulations of $S$, and edges connect triangulations which are related by certain elementary moves, namely reorderings, flips, and powers of Dehn twists.
Roughly speaking, the length of an edge is proportional to the computational complexity of performing the corresponding move.
More precisely, we can describe the three types of edges in $\GGG(S)$ as follows.
\begin{itemize}
\item If two triangulations $\TTT$ and $\TTT'$ have the same edges up to reordering, they are joined by an edge of length $1$ in $\GGG(S)$.
\item If $e$ is an edge of a triangulation $\TTT$ which is adjacent to two different triangles $\Delta_1$ and $\Delta_2$, a \emph{flip} about $e$ yields a new triangulation $\TTT'$ obtained by replacing $e$ with the opposite diagonal of the square $\Delta_1\cup\Delta_2$ (see \cite[\S2.2]{bell-mapping-classes}).
In the graph $\GGG(S)$, the triangulations $\TTT$ and $\TTT'$ are joined by an edge of length $1$.
\item Given an essential curve $a$ on $S$, denote by $\twist_a$ the Dehn twist of $S$ about $a$.
If $\TTT$ is a triangulation of $S$ and $\TTT'$ is the image of $\TTT$ under $\twist_a^k$ for some integer $k\neq 0$, then there is an edge of length $\norm[\TTT,1]{a}+\log|k|$ in $\GGG(S)$ between $\TTT$ and $\TTT'$. 
\end{itemize}
When $S$ is closed, its flip and twist graph is the flip and twist graph of the surface obtained by adding a fictitious marked point to $S$.

The graph $\GGG(S)$ is connected, as proved by \textcite{hatcher-triangulations}, even if we ignore the edges corresponding to powers of Dehn twists.
If we have a path of length $\ell$ in $\GGG(S)$ between two triangulations $\TTT$ and $\TTT'$, then we can perform a change of coordinates of normal curves between the two triangulations in polynomial time in $\ell$.
More precisely, if $a$ is a normal curve on $S$ with respect to the triangulation $\TTT$, then we can compute the vector $\TTT'(a)$ in polynomial time in $\xi$, $\ell$, and $\norm[\TTT]{a}$.
This follows from \cite[Proposition 2.2.1]{bell-mapping-classes} for flips and \cite[Theorem 4.1]{schaefer-computing-dehn-twists} for powers of Dehn twists.

\step{Representing homeomorphisms.}
There are a few reasonable combinatorial representations of homeomorphisms $\map{f}{S}{S}$ that are amenable to algorithmic manipulations. The algorithms we describe in \cref{sec:classification} are largely agnostic when it comes to this choice; our only requirement is that the selected representation comes equipped with a notion of complexity $\norm{-}$ satisfying the following properties.
\begin{enumroman}
\item\label[condition]{itm:homeomorphism representation:1} Given a homeomorphism $\map{f}{S}{S}$ of complexity $\norm{f}$ and an essential curve $a$ on $S$, there is an algorithm to compute the curve $f(a)$ in polynomial time in $\xi$, $\norm{f}$, and $\norm{a}$.
\item\label[condition]{itm:homeomorphism representation:2} For a homeomorphism $\map{f}{S}{S}$ and an essential curve $a$ on $S$, the complexities of $f$, $a$, and $f(a)$ are related by the inequality
\[
    \norm{f(a)}_1\le\poly(\xi)\cdot\poly(\norm{f})+\norm{a}_1.
\]
\end{enumroman}
Note that, in \cref{itm:homeomorphism representation:2}, we are using the complexity $\norm{-}_1$ of curves.

In order to provide a concrete example, we now describe how homeomorphisms can be represented as paths in the flip and twist graph of $S$.
Let $\TTT$ be a fixed triangulation of the surface $S$. If $\map{f}{S}{S}$ is a homeomorphism, then $f(\TTT)$ is also a triangulation of $S$; if $S$ is closed, up to isotopy we may assume that $f$ fixes the vertex of $\TTT$.
Since the flip and twist graph $\GGG(S)$ is connected, there is a path in $\GGG(S)$ from $\TTT$ to $f(\TTT)$, and we can use this path as a combinatorial representation of the homeomorphism $f$.
We then define the complexity $\norm{f}$ to be the length of this path.

We have already discussed how, given an essential curve $a$ on $S$, we can compute the curve $f(a)$ in polynomial time in $\xi$, $\norm{f}$, and $\norm{a}$.
As far as \cref{itm:homeomorphism representation:2} is concerned, it is not hard to check that the inequality
\[
    \norm{f(a)}_1\le3\norm{f}+\norm{a}_1
\]
is satisfied for every essential curve $a$ on $S$.

\subsection{Curve graph}

One of the main goals of this article is to give an efficient algorithm to coarsely compute distances in the curve graph of a surface, which we now define following \textcite{harvey-curve-complex}. Let $S$ be a surface with $\xi(S)\ge 2$.
The \emph{curve graph} of $S$ is a graph whose vertex set is the set $\C_0(S)$ of isotopy classes of essential curves in $S$.
Two classes of curves are connected by an edge if they admit disjoint representatives.
By stipulating that every edge has length $1$, the curve graph becomes a metric space, and we can talk about the distance $\dist[\C(S)](a,b)$ between two (isotopy classes of) curves $a$ and $b$.

The curve graph is connected, and in fact we can give a more precise upper bound for the distance between two curves in terms of their intersection number. 
The following proposition combines several known upper bounds of this form.


\begin{proposition}\label{thm:distance bound from intersections}
Let $a,b\in\C_0(S)$ be essential curves with $i(a,b)>0$. Then
\[
    \dist[\C(S)](a,b)\le\min\left\{i(a,b)+1,2\log i(a,b)+2,\frac{6}{\log\xi}\cdot \log i(a,b)+2\right\}.
\]
If $i(a,b)\le 2$ then we have the better bound
\[
    \dist[\C(S)](a,b)\le 2.
\]
\end{proposition}
\begin{proof}
Firstly, note that if $a$ and $b$ intersect at most twice, then they cannot fill a surface of complexity $\xi\ge 2$, and hence $\dist[\C(S)](a,b)\le 2$.
In the general case, the inequality $\dist[\C(S)](a,b)\le i(a,b)+1$ is proved in \cite[Lemma 1.1]{bowditch-intersection-numbers}.
The logarithmic bound $\dist[\C(S)](a,b)\le 2\log i(a,b)+2$ is known as \citeauthor{hempel-curve-complex}'s lemma \cite[Lemma 2.1]{hempel-curve-complex}, and is proved in \cite[Lemma 1.21]{schleimer-curve-complex} in the punctured case.

Finally, the last bound is a consequence of Lemma 2.1 of \cite{bowditch-uniform-hyperbolicity}. In fact, it follows from \citeauthor{bowditch-uniform-hyperbolicity}'s result that
\[
    \dist[\C(S)](a,b)\le\frac{2}{\log(2g+p-3)-1}\cdot\log i(a,b)+2;
\]
we get the desired inequality by noting that $\log\xi\le 3\max\{1,\log(2g+p-3)-1\}$.
\end{proof}

When $S$ is the sphere with $4$ punctures or the torus with $0$ or $1$ punctures, the definition above would yield a totally disconnected graph. Therefore, in these cases, we change the definition by stipulating that two classes of curves are connected by an edge if their intersection number is $1$ -- in the case of the torus with at most one puncture -- or $2$ -- in the case of the sphere with $4$ punctures.

When $S$ is an annulus (that is, a twice\=/punctured sphere), we adopt yet another definition. In the scope of this article, we will only be concerned with curve graphs of annuli when dealing with the subsurface projection map; therefore, we defer the definition of $\C(S)$ in the annular case to \cref{sec:classification:subsurface projection}.

Finally, let us introduce the following notation: if $A$ and $B$ are non\=/empty subsets of $\C_0(S)$, we define
\begin{align*}
    \dist[\C(S)](A,B)&=\min\{\dist[\C(S)](a,b):a\in A,b\in B\},\\
    \diam[\C(S)](A)&=\sup\{\dist[\C(S)](a,a'):a,a'\in A\},\\
    \diam[\C(S)](A,B)&=\sup\{\dist[\C(S)](a,b):a\in A,b\in B\}.
\end{align*}

\subsection{Combinatorics of train tracks}

Train tracks are one of the main tools to compute lower bounds for distances in the curve graph. We give a concise review of the pertinent definitions and properties, and then introduce new notations tailored to our algorithmic purposes. We will try to clearly differentiate between concepts derived from established literature and novel definitions. Our primary reference for train tracks, and in particular for their role in the computation of distances in the curve graph, is \cite[\S4]{masur-minsky-1}.

\step{Definition.} A \emph{train track} on a surface $S$ is an embedded $1$\=/complex $\tau\subs S$ such that every edge (called a \emph{branch}) is smooth and has well\=/defined tangent vectors at its endpoints; at each vertex (called a \emph{switch}) the incident branches are mutually tangent, and there is at least one branch for each of the two possible tangent directions. The valence of each switch is at least $3$, except for exactly one bivalent switch on each simple closed curve component of $\tau$. Finally, the connected components of $S\setminus\tau$ (which we call \emph{complementary regions}) are surfaces whose boundaries consist of smooth arcs interleaved with cusps; we require that no complementary region is a nullgon (that is, a smooth disc), a monogon (that is, a disc with one cusp), a bigon (that is, a disc with two cusps), a once\=/punctured nullgon, or a smooth annulus.

\step{Directions and ends.} We denote by $\switches(\tau)$ and $\branches(\tau)$ the sets of switches and branches of $\tau$ respectively. A simple counting argument (see \cite[Corollary 1.1.3]{penner-train-tracks}) shows that
\[
    \card{\switches(\tau)}\le 4\xi\qquad\text{and}\qquad\card{\branches(\tau)}\le6\xi.
\]
A \emph{direction} of a switch $s\in\switches(\tau)$ is a non\=/zero tangent vector to $\tau$ at $s$, defined up to multiplication by a positive scalar; in particular, each switch has two directions. If $\eta$ is a direction of $s$, we denote the opposite direction by $\bar{\eta}$.

We say that a branch of $\tau$ is \emph{returning} if it is adjacent to a single switch $s$, and its two endpoints have opposite directions at $s$.

A \emph{half\=/branch} of a branch $b\in\branches(\tau)$ is a component of $b\setminus q$ for a point $q$ in the interior of $b$.
An \emph{end} of $b$ is an equivalence class of half\=/branches of $b$, where two half\=/branches are equivalent if their intersection is still a half\=/branch of $b$.
It is clear that every branch $b$ has two ends; we will sometimes call one of them $b_\bullet$, and implicitly refer to the other one as $b^\bullet$ -- or vice versa.

The switch \emph{incident} to an end $e$, denoted by $\switch(e)$, is the unique switch of $\tau$ lying in the closure of a half\=/branch representing $e$.
The \emph{direction} of $e$, denoted by $\direction(e)$, is the direction of $\switch(e)$ represented by the tangent vector at the endpoint of a half\=/branch representing $e$.

Denote the set of ends of $\tau$ -- that is, the set of ends of branches of $\tau$ -- by $\ends(\tau)$.
For every switch $s$ of $\tau$ and every direction $\eta$ of $s$, we define
\[
    \listofends{\tau}{\eta}=\{e\in\ends(\tau):\direction(e)=\eta\}.
\]

\step{Subtracks and carrying.} For our purposes, a \emph{carrying map} between two train tracks (which may be multicurves) $\tau'$ and $\tau$ is a map $\map{\carrymap}{S}{S}$ homotopic to the identity such that:
\begin{enumroman}
\item $\carrymap(\tau')\subs\tau$;
\item $\carrymap(s')\in\switches(\tau)$ for every $s'\in\switches(\tau')$;\footnote{This condition is not part of the usual definition of carrying map in the literature.}
\item the differential of the restriction $\map{\carrymap|_{\tau'}}{\tau'}{\tau}$ is nowhere vanishing.
\end{enumroman}

If there exists a carrying map between $\tau'$ and $\tau$, we say that $\tau'$ is \emph{carried} by $\tau$, and we write $\tau'\carried\tau$. A very special case of carrying occurs when $\tau'$ is a subset of $\tau$; in this case, we say that $\tau'$ is a \emph{subtrack} of $\tau$, and we write $\tau'<\tau$. Whenever we write $\tau'\carried\tau$, we always implicitly assume that the carrying is induced by a specific carrying map, which we denote by $\carrymap[\tau'][\tau]$.

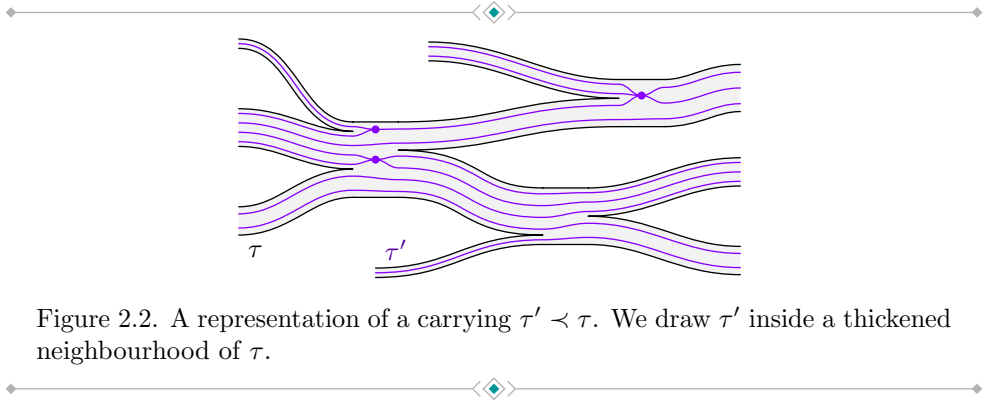
\begin{figure}
\centering
    \tikzsetnextfilename{carrying-example}%
    \input{figures-source/carrying-example.tex}%

\caption{A representation of a carrying $\tau'\carried\tau$. We draw $\tau'$ inside a thickened neighbourhood of $\tau$.}
\label{fig:carrying example}
\end{figure}

Our preferred way of depicting a carrying $\tau'\carried\tau$ is shown in \cref{fig:carrying example}: we draw a thickened neighbourhood of $\tau$ instead of $\tau$ itself, and we represent $\tau'$ so that it is contained in this neighbourhood.

\step{Pre\=/measures.} A \emph{pre\=/measure} on $\tau$ is a function
\[
    \map{\mu}{\branches(\tau)}{\ZZ_{\ge 0}},
\]
where $\ZZ_{\ge 0}$ denotes the set of non\=/negative integers.
We say that $\mu$ \emph{fills} $\tau$ if $\mu(b)>0$ for every $b\in\branches(\tau)$. 
With slight abuse of notation, if $e$ is an end of a branch $b$, we will write $\mu(e)$ for $\mu(b)$.

For every branch $b$ of $\tau$, define $\one{b}$ to be the pre\=/measure assigning weight $1$ to $b$ and weight $0$ to every other branch of $\tau$.
 Moreover, if $\mu$ and $\mu'$ are pre\=/measures on $\tau$, we will write $\mu+\mu'$ for the pre\=/measure obtained by adding $\mu$ and $\mu'$ pointwise, and $\mu\ge\mu'$ as a shorthand for ``$\mu(b)\ge\mu'(b)$ for every $b\in\branches(\tau)$''.

Two remarks are in order. Firstly, the notion of pre\=/measure is not found in the literature. In fact, the only reason why we give this definition instead of directly jumping to that of measure (see below) is that it is sometimes convenient in the description of our algorithms to refer to the pre\=/measures $\one{b}$. Secondly, it is common in the literature to define measures as taking non\=/negative real values. However, since we will focus solely on curves (and not on laminations), we have chosen a more combinatorial approach by restricting our definition to integral (pre\=/)measures.

\step{Measures.} A pre\=/measure $\mu$ on $\tau$ is a \emph{measure} if for every switch $s\in\switches(\tau)$ and every direction $\eta$ of $\tau$ the equality
\[
    \sum_{e\in\listofends{\tau}{\eta}}\mu(e)=\sum_{e\in\listofends{\tau}{\bar{\eta}}}\mu(e)
\]
holds.
For a measure $\mu$ and a switch $s$, we denote the quantity defined by either side of the above equation by $\mu(s)$.

Every non\=/zero measure $\mu$ on $\tau$ represents a non\=/empty essential multicurve carried by $\tau$. This multicurve is constructed by taking $\mu(b)$ parallel copies of $b$ for each branch $b\in\branches(\tau)$, and gluing corresponding arcs at switches. Every essential multicurve carried by $\tau$ is represented by some measure on $\tau$, and different measures always represent non\=/isotopic multicurves (see, for instance, Proposition 2.7.4 in \textcite{penner-train-tracks}'s book). For this reason, we will often blur the distinction between measures on a train track and essential multicurves carried by it.

We call a pair $(\tau,\mu)$ a \emph{measured train track} if $\tau$ is a train track and $\mu$ is a measure on $\tau$; we say that $(\tau,\mu)$ is \emph{\fullymeasured{}} if $\mu$ fills $\tau$ and represents a connected curve. A customary way to graphically represent a measured train track $(\tau,\mu)$ is to draw a thickened neighbourhood of $\tau$ instead of $\tau$ itself, where the thickness of a branch is proportional to the value assigned to it by $\mu$.

For a train track $\tau$, we define $\P(\tau)$ to be the set of measures on $\tau$. Moreover, we let
\[
    \intP(\tau)=\{\mu\in\P(\tau):\text{$\mu$ fills $\tau$}\}.
\]
A train track $\tau$ is \emph{recurrent} if $\intP(\tau)$ is non\=/empty.

Given two train tracks $\tau'\carried\tau$, we say that $\tau'$ \emph{fills} $\tau$ if $\intP(\tau')\subs\intP(\tau)$; note that talking about containment of these two sets makes sense if we interpret them as sets of isotopy classes of multicurves on $S$.

\step{Fundamental curves.} A \emph{fundamental measure}\footnote{Note that this definition is less restrictive than that of \emph{vertex measure} usually found in the literature; for this reason, we have decided to employ a different name.} on a train track $\tau$ is a non\=/zero measure $\nu\in\P(\tau)$ that cannot be written as the sum of two non\=/zero measures on $\tau$. We denote the set of fundamental measures on $\tau$ by $\fundamental(\tau)$. A fundamental measure will obviously represent a connected curve on $S$; therefore, we will sometimes refer to fundamental measures as ``fundamental curves''.

It is easy to see that, if $\nu$ is a fundamental measure on $\tau$, then $\nu(s)\le 2$ for every switch $s\in\switches(\tau)$. In fact, assume for contradiction that $\nu(s)\ge 3$ for some switch $s$ of $\tau$, and fix an orientation for the curve $a$ represented by $\nu$. The curve $a$ will run through $s$ at least twice with the same orientation, as depicted in \cref{fig:fundamental measures are small:a}. With a local cut\=/and\=/paste operation near $s$, as shown in \cref{fig:fundamental measures are small:b}, we transform $a$ into two immersed curves $a_1$ and $a_2$; these curves may have self\=/intersections, but they give rise to two non\=/zero measures $\nu_1$ and $\nu_2$ by counting how many times they run along branches of $\tau$. It is clear that $\nu=\nu_1+\nu_2$, hence $\nu$ is not fundamental.

\begin{figure}
\centering
\begin{subcaptionblock}{.45\linewidth}
\centering
    \tikzsetnextfilename{fundamental-measures-are-small-a}%
    \input{figures-source/fundamental-measures-are-small-a.tex}%

\caption{}
\label{fig:fundamental measures are small:a}
\end{subcaptionblock}
\begin{subcaptionblock}{.45\linewidth}
\centering
    \tikzsetnextfilename{fundamental-measures-are-small-b}%
    \input{figures-source/fundamental-measures-are-small-b.tex}%

\caption{}
\label{fig:fundamental measures are small:b}
\end{subcaptionblock}
\caption{\subref{fig:fundamental measures are small:a} If $\nu(s)\ge 3$, then $a$ runs through $s$ at least twice with the same orientation (towards the right in this picture). \subref{fig:fundamental measures are small:b} A cut\=/and\=/paste operation near $s$ turns $a$ into two (possibly self\=/intersecting) curves $a_1$ and $a_2$ representing two measures $\nu_1$ and $\nu_2$ such that $\nu=\nu_1+\nu_2$.}
\label{fig:fundamental measures are small}
\end{figure}

\step{Largeness and diagonal extensions.} A train track $\tau$ is \emph{large} if its complementary regions are topological discs with at most one puncture; in other words, if every essential curve in $S$ intersects $\tau$. Note that if $\tau$ is large and carried by a train track $\tau'$, then $\tau'$ is large as well.

We say that a train track $\rho$ is a \emph{diagonal extension} of a train track $\tau$ if $\tau<\rho$ and $\rho\setminus\tau$ consists of branches whose ends terminate in cusps of complementary regions of $\tau$. If $\tau'$ is a large recurrent train track which is carried by $\tau$ and fills $\tau$, then every diagonal extension of $\tau'$ is carried by a diagonal extension of $\tau$ (see the proof of Lemma 4.2 in \cite{masur-minsky-1}).

For a large train track $\tau$, we define the following sets of measures/multicurves:
\bgroup\allowdisplaybreaks
\begin{align*}
\PE(\tau)&=\bigcup\{\P(\rho):\text{$\rho$ is a recurrent diagonal extension of $\tau$}\},\\
\intPE(\tau)&=\{\mu\in\PE(\tau):\text{$\mu(b)>0$ for every $b\in\branches(\tau)$}\},\\
\PN(\tau)&=\bigcup\{\PE(\sigma):\text{$\sigma$ is a large recurrent subtrack of $\tau$}\},\\
\intPN(\tau)&=\bigcup\{\intPE(\sigma):\text{$\sigma$ is a large recurrent subtrack of $\tau$}\}.
\end{align*}
\egroup

\step{Distances in $\C(S)$.} Let us begin with a trivial observation: any two curves carried by a non\=/large train track lie at distance at most two in $\C(S)$. This follows from the fact that a train track $\tau$ is not large precisely when there is an essential curve on $S$ disjoint from $\tau$.

\begin{figure}
\centering
    \tikzsetnextfilename{measures-intersection-number}%
    \input{figures-source/measures-intersection-number.tex}%

\caption{The curves $a$ and $a'$ can be realised so that they run parallel to each other along branches of $\tau$, and intersect at most $\mu(s)\cdot\mu'(s)$ near a switch $s$.}
\label{fig:measures intersection number}
\end{figure}
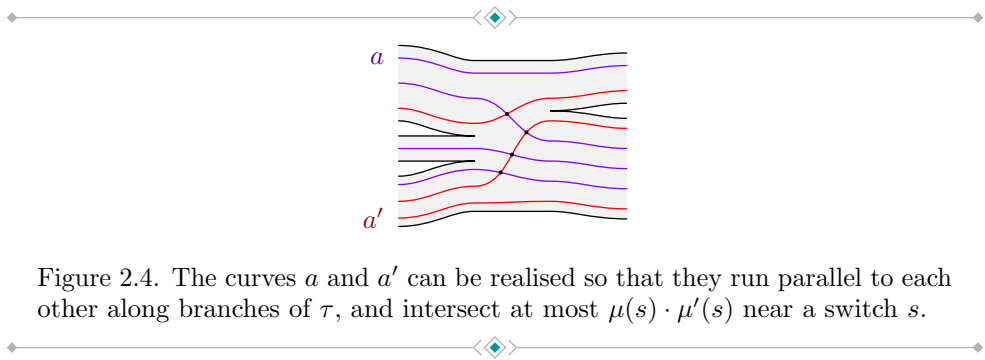

Next, let us give an upper bound for the intersection number of curves carried by a train track $\tau$. If $\mu$ and $\mu'$ are measures on $\tau$ representing curves $a$ and $a'$ respectively, then
\[
    i(a,a')\le\sum_{s\in\switches(\tau)}\mu(s)\cdot\mu'(s).
\]
In fact, $a$ and $a'$ can be isotoped so that they run parallel to the branches of $\tau$, and only intersect near the switches; in particular, near a switch $s\in\switches(\tau)$, they intersect at most $\mu(s)\cdot\mu'(s)$ times, as depicted in \cref{fig:measures intersection number}.

We now state what is probably the most crucial result on which our algorithm for computing coarse distances in $\C(S)$ relies.

\begin{proposition}\label{thm:nesting lemma}
Let $\tau$ be a large recurrent train track, and let $a\in\intPN(\tau)$ be an essential curve. Let $b\in\C_0(S)$ be any curve such that $i(a,b)=0$. Then $b\in\PN(\tau)$.
\end{proposition}

This proposition first appeared as Lemma 4.4 in \cite{masur-minsky-1}, with the additional constraint that $\tau$ should be \emph{birecurrent}. However, the assumption of \emph{transverse recurrence} can be dropped, as shown by \textcite[Lemma 3.2]{gadre-min-translation-length}; they attribute the proof to Leininger. Note that \citeauthor{gadre-min-translation-length} state the lemma for $\intPE$ and $\PE$ instead of $\intPN$ and $\PN$. However, our \cref{thm:nesting lemma} is an immediate consequence of their result. In fact, $a\in\intPN(\tau)$ means that $a\in\intPE(\sigma)$ for a large recurrent subtrack $\sigma<\tau$. By Lemma 3.2 of \cite{gadre-min-translation-length}, this implies that $b\in\PE(\sigma)$, and by definition $\PE(\sigma)\subs\PN(\tau)$.

An easy consequence of \cref{thm:nesting lemma} allows us to provide lower bounds for distances in the curve graph in terms of specific nested sequences of train tracks.

\begin{corollary}\label{thm:nested sequence distance}
Let $\tau_0\carries\tau_1\carries\ldots\carries\tau_n$ be a sequence of large recurrent train tracks such that
\[
    \PN(\tau_i)\subs\intPN(\tau_{i-1})\qquad\text{for every $1\le i\le n$}.
\]
Then there is a curve $a\in\fundamental(\tau_0)$ such that
\[
    \dist[\C(S)](a,\fundamental(\tau_n))\ge n.
\]
\end{corollary}
\begin{proof}
Let $b$ be a fundamental curve of $\tau_n$; in particular, note that $b\in\PN(\tau_n)$. Making use of \cref{thm:nesting lemma}, we can prove by induction on $0\le k<n$ that every essential curve whose distance from $b$ is at most $k$ lies in $\intPN(\tau_{n-k-1})$.
Let us now fix a curve $a$ representing a \emph{vertex measure} of $\tau_0$ according to the definition in \cite[\S4.1]{masur-minsky-1}. Note that $a$ is a fundamental curve of $\tau_0$. By Lemma 3.5 in \cite{masur-minsky-quasiconvexity}, the curve $a$ does not lie in $\intPN(\tau_0)$; by the above discussion, this implies that $\dist[\C(S)](a,b)\ge n$. Since $a$ can be chosen independently of $b$, this concludes the proof.
\end{proof}

\step{Co\=/orientations and height.} Let us remark that all the terminology we will introduce from now on is not standard in the literature; as anticipated, some of these definitions are specifically tailored to our algorithms, and would probably not be useful in most other contexts.

Loosely speaking, a co\=/orientation of a switch $s\in\switches(\tau)$ is a unit tangent vector orthogonal to $\tau$ at $s$, as depicted in \cref{fig:co-orientation definition}. More formally, a \emph{co\=/orientation} of $s$ is a tangent vector to $S$ at the point $s$ which is not tangent to $\tau$, defined up to the following equivalence relation: two co\=/orientations $\omega$ and $\omega'$ are the same if, for every direction $\eta$ of $s$, the bases $\{\eta,\omega\}$ and $\{\eta,\omega'\}$ of the tangent space to $S$ at $s$ are related by an orientation\=/preserving linear isomorphism. Clearly, every switch has two co\=/orientations; if $\omega$ is a co\=/orientation of $s$, we denote the opposite co\=/orientation by $\bar{\omega}$.

Fix a switch $s$ and a direction $\eta$ of $s$. A co\=/orientation $\omega$ of $s$ defines a linear order on the set of ends $\listofends{\tau}{\eta}$, given by enumerating the ends following the vector $\omega$; see \cref{fig:co-orientation order} for an example. For two ends $e_1,e_2\in\ends(\tau)$, we will say that $e_1$ is \emph{$\omega$\=/smaller} than $e_2$ (and write $e_1\endl{\omega}e_2$) if $\switch(e_1)=\switch(e_2)=s$, $\direction(e_1)=\direction(e_2)$, and $e_1$ is smaller than $e_2$ according to the order defined above.

\begin{figure}
\centering
\begin{subcaptionblock}{.45\linewidth}
\centering
    \tikzsetnextfilename{co-orientation-definition}%
    \input{figures-source/co-orientation-definition.tex}%

\caption{}
\label{fig:co-orientation definition}
\end{subcaptionblock}
\begin{subcaptionblock}{.45\linewidth}
\centering
    \tikzsetnextfilename{co-orientation-order}%
    \input{figures-source/co-orientation-order.tex}%

\caption{}
\label{fig:co-orientation order}
\end{subcaptionblock}
\caption{\subref{fig:co-orientation definition} A co\=/orientation $\omega$ of a switch $s$ can be informally thought of as a unit tangent vector orthogonal to $\tau$ at $s$. \subref{fig:co-orientation order} The linear order induced on $\listofends{\tau}{\eta}$ by the co\=/orientation $\omega$.}
\end{figure}

We denote by $\listofends{\tau}{\eta}[\omega]$ the set $\listofends{\tau}{\eta}$ endowed with the linear order $\endl{\omega}$. For $1\le i\le\card{\listofends{\tau}{\eta}}$, we denote by $\listofends{\tau}{\eta}[\omega]_i$ the $i$\=/th element of $\listofends{\tau}{\eta}[\omega]$. We will say that two ends of $\tau$ are \emph{consecutive} if they are of the form $\listofends{\tau}{\eta}[\omega]_i$ and $\listofends{\tau}{\eta}[\omega]_{i+1}$ for some switch $s\in\switches(\tau)$, some direction $\eta$ of $s$, some co\=/orientation $\omega$ of $s$, and some integer $i$; note that this condition does not depend on the choice of $\omega$.

If $\mu$ is a pre\=/measure on $\tau$, $s$ a switch of $\tau$, $\omega$ a co\=/orientation of $s$, and $e$ an end with $\switch(e)=s$, then we define the \emph{height} of $e$ as
\[
    \height{\omega}(e)=\sum_{\substack{e'\in\ends(\tau)\\e'\endg{\omega}e}}\mu(e').
\]
The pre\=/measure $\mu$ is omitted from the notation, but will always be clear from the context. The name ``height'' suggests the convention that, in our pictures, co\=/orientations will always point downwards.

\step{Carrying revisited.} Recall that, according to our definition, a carrying map $\carrymap$ between two train tracks $\tau'$ and $\tau$ sends switches to switches. As a consequence, there are induced maps -- which we keep calling $\carrymap$ -- sending directions, co\=/orientations, ends, and pre\=/measures on $\tau'$ to directions, co\=/orientations, ends, and pre\=/measures on $\tau$. Let us now state two properties of these maps.

\begin{itemize}
\item If $s'$ is a switch of $\tau'$, $\eta'$ a direction of $s'$, $\omega'$ a co\=/orientation of $s'$, $\eta=\carrymap(\eta')$, and $\omega=\carrymap(\omega')$, then the map
\[
    \map{\carrymap}{\listofends{\tau'}{\eta'}[\omega']}{\listofends{\tau}{\eta}[\omega]}
\]
is order\=/preserving. This is essentially obvious, as the reader can see from \cref{fig:carrying preserves ends order}.
\item The map $\carrymap$ sending pre\=/measures on $\tau'$ to pre\=/measures on $\tau$ is $\ZZ$\=/linear, and is uniquely determined by the pre\=/measures $\{\carrymap(\one{b'}):b'\in\branches(\tau')\}$. As a consequence, it sends measures to measures, and there is an induced $\ZZ$\=/linear map
\[
    \map{\carrymap}{\P(\tau')}{\P(\tau)}.
\]
\end{itemize}

\begin{figure}
\centering
    \tikzsetnextfilename{carrying-preserves-ends-order}%
    \input{figures-source/carrying-preserves-ends-order.tex}%

\caption{The function $\umap{\listofends{\tau'}{\eta'}[\omega']}{\listofends{\tau}{\eta}[\omega]}$ induced by a carrying map $\carrymap$ is order\=/preserving.}
\label{fig:carrying preserves ends order}
\end{figure}

In the context of a carrying $\tau'\carried\tau$, we will say that a branch $b\in\branches(\tau)$ is \emph{untouched} in $\tau'$ if $\carrymap[\tau'][\tau](\one{b'})=\one{b}$ for some branch $b'\in\branches(\tau')$. We will say that an end $e\in\ends(\tau)$ is \emph{persistent} in $\tau'$ if $\carrymap[\tau'][\tau](e')=e$ for some end $e'\in\ends(\tau')$.

If $(\tau,\mu)$ and $(\tau',\mu')$ are measured train tracks, we write $(\tau',\mu')\carried(\tau,\mu)$ to signify that $\tau'\carried\tau$ and $\carrymap[\tau'][\tau](\mu')=\mu$.

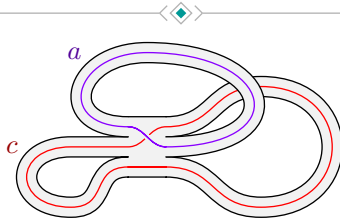
\begin{figure}
\centering
    \tikzsetnextfilename{one-switch-fundamental-curves}%
    \input{figures-source/one-switch-fundamental-curves.tex}%

\caption{The two types of fundamental curves of a one\=/switch train track: the curve $a$ is represented by $\one{b}$ for a returning branch $b$, while $c$ is represented by $\one{b_1}+\one{b_2}$ for non\=/returning branches $b_1$ and $b_2$.}
\label{fig:one-switch fundamental curves}
\end{figure}

\step{One\=/switch train tracks.} A \emph{one\=/switch train track} is a train track with only one switch. One\=/switch train tracks will be particularly relevant for our algorithms, since they are more convenient to work with than general train tracks. The main reason for this is that they only have $\bigO(\xi^2)$ fundamental curves, which can be described explicitly. In fact, if $\tau$ is a one\=/switch train track, then
\begin{align*}
\fundamental(\tau)&=\{\one{b}:\text{$b\in\branches(\tau)$ is returning}\}\\
&\phantom{{}={}}\cup\{\one{b_1}+\one{b_2}:\text{$b_1,b_2\in\branches(\tau)$ are not returning and $\direction(b_1^\bullet)\neq\direction(b_2^\bullet)$}\}.
\end{align*}
\Cref{fig:one-switch fundamental curves} shows the two types of fundamental curves of a one\=/switch train track.

We observe that two fundamental curves of a one\=/switch train track $\tau$ intersect at most $4$ times. The bound in \cref{thm:distance bound from intersections} implies that
\begin{equation}\label{eqn:one-switch train track fundamental diameter}
    \diam[\C(S)](\fundamental(\tau))\le 5.
\end{equation}

\step{Deep nesting.} Let $\tau'\carried\tau$ be a pair of one\=/switch train tracks. We say that $\tau'$ is \emph{deeply nested} in $\tau$ if exactly two ends of $\tau$ are persistent in $\tau'$. \Cref{fig:deeply nested example} shows what two deeply nested train tracks typically look like near their unique switches. In particular, the train track $\tau$ has two branches $b_1$ and $b_2$ (not necessarily distinct) with ends $b_1^\bullet$ and $b_2^\bullet$ satisfying $\direction(b_1^\bullet)\neq\direction(b_2^\bullet)$ and $\carrymap[\tau'][\tau](e')\in\{b_1^\bullet,b_2^\bullet\}$ for every end $e'$ of $\tau'$.

We now state two properties of deeply nested train tracks which will be useful in the proof of \cref{thm:distance algorithm for carried curves}.
We use the notation we introduced in the previous paragraph.

\begin{lemma}\label{thm:deep nesting lemma}
Let $\tau'\carried\tau$ be deeply nested one\=/switch train tracks.
\begin{enumarabic}
\item There is a fundamental measure $\nu$ on $\tau$, supported on the branches $b_1$ and $b_2$, such that $\carrymap[\tau'][\tau](\mu')\ge\nu$ for every non\=/zero measure $\mu'$ on $\tau'$.
\item Suppose that $\tau'$ is large, recurrent, and fills $\tau$. Then every recurrent diagonal extension $\rho'$ of $\tau'$ is carried by and deeply nested in a recurrent diagonal extension $\rho$ of $\tau$; the two ends of $\rho$ which are persistent in $\rho'$ are $b_1^\bullet$ and $b_2^\bullet$.
\end{enumarabic}
\end{lemma}
\begin{proof}
We prove the two statements separately.
\begin{enumarabic}
\item Let $\mu'$ be a non\=/zero measure on $\tau'$. There must be two ends $e_1'$ and $e_2'$ of $\tau'$ with $\direction(e_1')\neq\direction(e_2')$ such that $\mu'(e_1')>0$ and $\mu'(e_2')>0$. Without loss of generality, assume that $\carrymap[\tau'][\tau](e_1')=b_1^\bullet$ and $\carrymap[\tau'][\tau](e_2')=b_2^\bullet$. If we let $\mu=\carrymap[\tau'][\tau](\mu')$, it immediately follows that $\mu(b_1)>0$ and $\mu(b_2)>0$. If $b_1$ (respectively $b_2$) is returning, we can take $\nu=\one{b_1}$ (respectively $\one{b_2}$). Otherwise, we take $\nu=\one{b_1}+\one{b_2}$. Either way, we have that $\nu$ is a fundamental measure on $\tau$ and $\mu\ge\nu$.
\item We already know, from the proof of \cite[Lemma 4.2]{masur-minsky-1}, that $\rho'$ is carried by a recurrent diagonal extension $\rho$ of $\tau$. Now, let $b'$ be a branch in $\branches(\rho')\setminus\branches(\tau')$, and let $b'_\bullet$ be one of its ends. Since $\rho'$ is a diagonal extension of $\tau'$, there are two consecutive ends $e_1'$ and $e_2'$ of $\tau'$ such that $e_1'\endl{\omega'}b'_\bullet\endl{\omega'}e_2'$ for some orientation $\omega'$ of the unique switch of $\tau'$; in other words, $b'_\bullet$ is contained in the cusp of $S\setminus \tau'$ defined by the two ends $e_1'$ and $e_2'$, as shown in \cref{fig:deeply nested diagonal extension}. But $\carrymap[\rho'][\rho](e_1')=\carrymap[\rho'][\rho](e_2')=b_i^\bullet$ for some $i\in\{1,2\}$, and $\carrymap[\rho'][\rho]$ preserves the order of ends, so $\carrymap[\rho'][\rho](b'_\bullet)=b_i^\bullet$ as well. It follows that $\rho'$ is deeply nested in $\rho$, and the only two persistent ends of $\rho$ in $\rho'$ are $b_1^\bullet$ and $b_2^\bullet$.\qedhere
\end{enumarabic}
\end{proof}

\begin{figure}
\centering
\begin{subcaptionblock}[c]{.45\linewidth}
\centering
    \tikzsetnextfilename{deeply-nested-example}%
    \input{figures-source/deeply-nested-example.tex}%

\caption{}
\label{fig:deeply nested example}
\end{subcaptionblock}
\begin{subcaptionblock}[c]{.45\linewidth}
\centering
    \tikzsetnextfilename{deeply-nested-diagonal-extension}%
    \input{figures-source/deeply-nested-diagonal-extension.tex}%

\caption{}
\label{fig:deeply nested diagonal extension}
\end{subcaptionblock}
\caption{\subref{fig:deeply nested example} Two deeply nested train tracks $\tau'\carried\tau$; we have that $\carrymap[\tau'][\tau](e')\in\{b_1^\bullet,b_2^\bullet\}$ for every end $e'$ of $\tau'$. \subref{fig:deeply nested diagonal extension} If the end $b'_\bullet$ is contained in the cusp of $\tau'$ defined by $e_1'$ and $e_2'$, and $\carrymap[\tau'][\tau](e_1')=\carrymap[\tau'][\tau](e_2')=b_i^\bullet$, then $\carrymap[\rho'][\rho](b'_\bullet)=b_i^\bullet$ as well.}
\end{figure}

While the significance of the ``deeply nested'' condition will become apparent in the proof of \cref{thm:distance algorithm for carried curves}, it is easy to get lost in the technical details of the argument, and the reader might be left wondering why we have introduced this concept in the first place.
At a very high level, the reason is exemplified by the following situation.

Let $\tau''\carried\tau$ be large recurrent one\=/switch train tracks.
In light of \cref{thm:nested sequence distance}, we might be interested in deciding whether the condition $\PN(\tau'')\subs\intPN(\tau)$ holds.
This containment is challenging to check in polynomial time in $\xi$, mainly because the number of diagonal extensions of large subtracts of $\tau''$ could a priori be doubly exponential in $\xi$.
Suppose, however, that there is a train track $\tau'$ with $\tau''\carried\tau'\carried\tau$ such that $\tau''$ is deeply nested in $\tau'$; let $b_1'$ and $b_2'$ be the two (not necessarily distinct) branches of $\tau'$ having an end which is persistent in $\tau''$.
The content of \cref{thm:deep nesting lemma} is essentially the following: every curve carried by a diagonal extension of a large subtrack of $\tau''$ is also carried by a diagonal extension of a large subtrack of $\tau'$, in such a way that it runs along the branches $b_1'$ and $b_2'$.
Therefore, in order to verify that $\PN(\tau'')\subs\intPN(\tau)$, it suffices\footnote{The reader will have surely realised that this is not an equivalent condition. To be more precise, the argument in the proof of \cref{thm:distance algorithm for carried curves} proceeds as follows. The intermediate train track $\tau'$ is chosen so that it is uniformly close to $\tau''$ in $\C(S)$. Now, if $\carrymap[\tau'][\tau](\nu)$ fills a large subtrack of $\tau$ -- where $\nu$ is the fundamental measure on $\tau'$ given by \cref{thm:deep nesting lemma} -- then $\PN(\tau'')\subs\intPN(\tau)$ as described. Otherwise, it is easy to see that $\tau'$ is close to $\tau$ in $\C(S)$, and hence so is $\tau''$.} to check whether $\carrymap[\tau'][\tau](\one{b_1'})$ and $\carrymap[\tau'][\tau](\one{b_2'})$ fill a large subtrack of $\tau$, which is straightforward to do in polynomial time in $\xi$.
In other words, the existence of an intermediate train track $\tau'$ in which $\tau''$ is deeply nested allows us to bypass the issue of having to iterate over (potentially too many) diagonal extensions of large subtracks of $\tau''$.

\step{Implementation and complexity.}
Finally, we address the matter of how to represent and manipulate train tracks on a surface from an algorithmic point of view.
We think of train tracks as abstract graphs, with vertices given by the switches and edges given by the branches; for each switch, we keep track of which branches are incident to each direction, as well as the order in which they appear along a co\=/orientation.
Additionally, we store some information about the complementary regions, such as their topological type and the sequence of branches that appears along their boundary.

In short, we represent train tracks on $S$ in such a way that we can reconstruct them up to homeomorphism of $S$, although not necessarily up to isotopy.
While it would be theoretically possible to encode the specific way in which a train track is embedded in $S$ -- using a technique similar to that of normal curves -- this information would be completely irrelevant for our purposes.
In fact, all we need is a representation that is detailed enough to encode largeness of a train track, and supports the \splitmove{} and \twistmove{} moves described in \cref{sec:distance algorithm:aht moves}.

We remark that, for any sensible representation of train tracks, all the operations we will need to perform -- and in particular the \splitmove{} and \twistmove{} moves -- can be implemented in polynomial time in the complexity $\xi$ of the surface (recall that the number of switches and branches of a train track on $S$ is $\bigO(\xi)$).
However, the exact \emph{degree} of this polynomial will depend on the specific data structure used to represent train tracks.
Therefore, in order to avoid having to deal with implementation details, we will be content with stating a generic polynomial dependence on $\xi$ of the running time of our algorithms.
Instead, we will try and be more precise about the degree of the dependence on the complexities of curves and surface homeomorphisms.

A pre\=/measure $\mu$ on a train track $\tau$ is represented by a function $\umap{\branches(\tau)}{\ZZ_{\ge 0}}$; we define the \emph{complexity} of $\mu$ to be 
\[
    \norm{\mu}=\sum_{b\in\branches(\tau)}\log(\mu(b)+1).
\]

To describe a carrying $\tau'\carried\tau$ combinatorially, we simply keep track of the action of the carrying map $\carrymap[\tau'][\tau]$ on switches, orientations, and co\=/orientations of $\tau'$, together with the pre\=/measures $\{\carrymap[\tau'][\tau](\one{b'}):b'\in\branches(\tau')\}$ on $\tau$.
This information might not be enough to recover $\tau'$ from $\tau$ up to isotopy, but it is sufficient for computing the image $\carrymap[\tau'][\tau](\mu')\in\P(\tau)$ of any measure $\mu'\in\P(\tau')$.
The \emph{complexity} of the carrying map $\carrymap[\tau'][\tau]$ is defined as
\[
    \norm{\carrymap[\tau'][\tau]}=\sum_{b'\in\branches(\tau')}\norm{\carrymap[\tau'][\tau](\one{b'})}.
\]

%% file: figures-source/normal-curve-example.tex
\begin{tikzpicture}
\foreach \i in {1,2,3} {\coordinate (\i) at (-30+120*\i:1.5);}

\foreach \i/\j/\k in {1/2/7,2/3/4,3/1/5} {
    \foreach \l in {1,...,\k} {
        \tikzmath{\t=\l/(\k+1)*.5+.25;}
        \coordinate (\i-\l) at ($(\i)!\t!(\j)$);
    }
}
\tikzset{every path/.style={line,draw=palette 2}}
\foreach \i in {1,...,4} {
    \tikzmath{\j=6-\i;}
    \draw (1-\i) to[out=-30,in=-150] (3-\j);
}
\foreach \i in {1,...,3} {
    \tikzmath{\j=8-\i;}
    \draw (1-\j) to[out=-30,in=90] (2-\i);
}
\draw (2-4) to[out=90,in=-150] (3-1);

\draw[black] (1) -- (2) -- (3) -- (1);
\node[above left,\colourmakedark{2}] at ($(1)!.5!(2)$) {$a$};

\end{tikzpicture}

%% file: figures-source/carrying-example.tex
\begin{tikzpicture}
\tikzset{
tt/new=a,
tt/a/.cd,
set sep=.8cm,
set width=.6cm,
add switch={0,0},
add switch={3.5,.75},
add switch={2.5,-.75},
add branch standard={1}{2},
add branch standard={1}{3},
add branch={++(a-3) to[out=180,in=0] +(-2.5,-.75)},
add bunch={1}{left}{3,spread=1.5cm,yshift=.3cm},
add bunch={2}{above 1}{1,spread=2.5cm},
add bunch={2}{right}{1,spread=1cm,yshift=.2cm},
add bunch={3}{right}{2,spread=2cm},
set switch thickness={1}{1cm},
set branch weights={4=1,5=4,6=3,1=3,2=5,7=2,8=5,3=1,9=3,10=3},
finalise,
}
\fill[grey background] [tt/a/use branch background all,tt/a/use switch background all];
\draw[line,black] [tt/a/use branch contour all,tt/a/use switch contour all];
\tikzset{
tt/a/carried={
name=b,
switch={1}{.4cm}{1,2}{1},
switch={1}{0cm}{4,5}{3,4},
switch={2}{.1cm}{1,2,3}{1,2},
4=1,5=4,6=2,1=2,2=4,7=2,8=3,3=1,9=3,10=2,
},
tt/b/finalise
}
\draw[line,palette 1] [tt/b/use branch all];
\foreach \i in {1,2,3} {
    \pic[palette 1] at (b-\i) {dot};
}
\node[black,anchor=west,yshift=-.4cm] at (spath cs:{\ttbranch{a}{6}} 1) {$\tau$};
\node[\colourmakedark{1},anchor=west,yshift=.3cm] at (spath cs:{\ttbranch{a}{3}} 1) {$\tau'$};
\end{tikzpicture}

%% file: figures-source/fundamental-measures-are-small-a.tex
\begin{tikzpicture}
\tikzset{
tt/new=a,
tt/a/.cd,
set sep=.2cm,
set width=1cm,
set spread=1cm,
add switch={0,0},
add bunch={1}{left}{3},
add bunch={1}{right}{2},
set switch thickness={1}{2cm},
set branch weights={1=2,2=2,3=2,4=3,5=3},
finalise,
carried={
    name=b,
    1=1,2=1,3=1,4=2,5=1,
}
}
\fill[grey background] [tt/a/use branch background all,tt/a/use switch background all];
\draw[black,line] [tt/a/use branch contour all,tt/a/use switch contour all];
\draw[line,palette 1,postaction={decorate,decoration={markings,mark=at position -.2 with {\arrowreversed{Stealth[]}}}}] [tt/b/use branch={1}];
\draw[line,palette 1,postaction={decorate,decoration={markings,mark=at position .2 with {\arrow{Stealth[]}}}}] [tt/b/use branch={2}];
\draw[line,palette 1,postaction={decorate,decoration={markings,mark=at position -.2 with {\arrowreversed{Stealth[]}}}}] [tt/b/use branch={3}];
\node[\colourmakedark{1},anchor=east,outer xsep=2pt] at (spath cs:{\ttbranch{b}{1}} 1) {$a$};
\end{tikzpicture}

%% file: figures-source/fundamental-measures-are-small-b.tex
\begin{tikzpicture}
\tikzset{
tt/new=a,
tt/a/.cd,
set sep=.2cm,
set width=1cm,
set spread=1cm,
add switch={0,0},
add bunch={1}{left}{3},
add bunch={1}{right}{2},
set switch thickness={1}{2cm},
set branch weights={1=2,2=2,3=2,4=3,5=3},
finalise,
carried={
    name=b,
    1=1,2=1,3=1,4=2,5=1,
    matches={1}{1/3,2/2,3/1},
}
}
\fill[grey background] [tt/a/use branch background all,tt/a/use switch background all];
\draw[black,line] [tt/a/use branch contour all,tt/a/use switch contour all];
\draw[line,palette 2,line over=grey background,postaction={decorate,decoration={markings,mark=at position -.2 with {\arrowreversed{Stealth[]}}}}] [tt/b/use branch={1}];
\draw[line,palette 2,line over=grey background,postaction={decorate,decoration={markings,mark=at position .2 with {\arrow{Stealth[]}}}}] [tt/b/use branch={2}];
\draw[line,palette 3,line over=grey background,postaction={decorate,decoration={markings,mark=at position -.2 with {\arrowreversed{Stealth[]}}}}] [tt/b/use branch={3}];
\node[\colourmakedark{2},anchor=east,outer xsep=2pt] at (spath cs:{\ttbranch{b}{1}} 1) {$a_1$};
\node[\colourmakedark{3},anchor=east,outer xsep=2pt] at (spath cs:{\ttbranch{b}{3}} 1) {$a_2$};
\end{tikzpicture}

%% file: figures-source/measures-intersection-number.tex
\begin{tikzpicture}
\tikzset{
tt/new=a,
tt/a/.cd,
set sep=.2cm,
set width=1cm,
set spread=1cm,
add switch={0,0},
add bunch={1}{left}{3},
add bunch={1}{right}{2},
set switch thickness={1}{2cm},
set branch weights={1=3,2=1,3=2,4=2,5=4},
finalise,
carried={
    name=b,
    1=3,2=1,3=3,4=2,5=5,
    matches={1}{2/4,4/5,5/6},
}
}
\fill[grey background] [tt/a/use branch background all,tt/a/use switch background all];
\draw[black,line] [tt/a/use branch contour all,tt/a/use switch contour all];
\draw[palette 1,line,name path=a] [tt/b/use branch={1,2,4,5}];
\draw[palette 2,line,name path=b] [tt/b/use branch={3,6,7}];
\fill[black,name intersections={of=a and b,total=\t}] \foreach \i in {1,...,\t} {(intersection-\i) circle (.75pt)};
\node[\colourmakedark{1},anchor=east,outer xsep=2pt] at (spath cs:{\ttbranch{b}{1}} 1) {$a$};
\node[\colourmakedark{2},anchor=east,outer xsep=2pt] at (spath cs:{\ttbranch{b}{7}} 1) {$a'$};
\end{tikzpicture}

%% file: figures-source/co-orientation-definition.tex
\begin{tikzpicture}
\tikzset{
tt/new=a,
tt/a/.cd,
set sep=1cm,
set spread=1cm,
add switch={0,0},
add bunch={1}{left}{3,spread=1.5cm,sep=1cm},
add bunch={1}{right}{4,spread=2cm,sep=.8cm},
finalise,
}
\draw[black,line] [tt/a/use branch all];
\draw[palette 4,line,-{Triangle[]}] ++(a-1) -- +(0,-.75) node[\colourmakedark{4},below] {$\omega$};
\pic[palette 1] at (a-1) {dot} node[\colourmakedark{1},above,yshift=2pt] {$s$};
\end{tikzpicture}

%% file: figures-source/co-orientation-order.tex
\begin{tikzpicture}
\tikzset{
tt/new=a,
tt/a/.cd,
set sep=1cm,
set spread=1cm,
add switch={0,0},
add bunch={1}{left}{3,spread=1.5cm,sep=1cm},
add bunch={1}{right}{4,spread=2cm,sep=.8cm},
finalise,
}
\draw[black,line] [tt/a/use branch all];
\draw[palette 4,line,-{Triangle[]}] ++(a-1) -- +(0,-.75) node[\colourmakedark{4},below] {$\omega$};
\draw[palette 4,line,-{Triangle[]}] ++(a-1) -- +(.75,0) node[pos=.5,\colourmakedark{4},above,yshift=2pt] {$\eta$};
\pic[palette 1] at (a-1) {dot};
\foreach \b[count=\i] in {4,5,6,7} {
    \tikzset{tt/a/point along branch={\b}{begin,1.5cm,name=\i}}
    \node[\colourmakedark{3},circle,fill=white,inner sep=0] at (\i) {$\i$};
}
\foreach \i[evaluate=\i as \j using int(\i+1)] in {1,2,3} {
    \node[\colourmakedark{3},rotate=-90] at ($(\i)!.5!(\j)$) {$<$};
}
\end{tikzpicture}

%% file: figures-source/carrying-preserves-ends-order.tex
\begin{tikzpicture}
\tikzset{
tt/new=a,
tt/a/.cd,
set sep=.2cm,
set width=1cm,
set spread=1cm,
add switch={0,0},
add bunch={1}{left}{2},
add bunch={1}{right}{3},
set switch thickness={1}{2cm},
set branch weights={1=3,2=3,3=2,4=2,5=2},
finalise,
carried={
    name=b,
    switch={1}{0pt}{1,2,3}{1,...,5},
    1=2,2=1,3=1,4=2,5=2,
}
}
\fill[grey background] [tt/a/use branch background all,tt/a/use switch background all];
\draw[black,line] [tt/a/use branch contour all,tt/a/use switch contour all];
\draw[palette 1,line] [tt/b/use branch all];
\draw[palette 4,line,-{Triangle[]}] ++(b-1) -- +(.75,0) node[\colourmakedark{4},right] {\contour{grey background}{$\eta'$}};
\draw[line,palette 4,-{Triangle[]}] ++(b-1) pic[palette 1] {dot} -- +(0,-.75) node[\colourmakedark{4},above left] {\contour{grey background}{$\omega'$}};
\coordinate (c) at ($(a-1)+(0,-1.2cm)$);
\draw[palette 4,line,-{Triangle[]}] ++(c) -- +(.75,0) node[\colourmakedark{4},below right] {$\eta$};
\draw[palette 4,line,-{Triangle[]}] ++(c) -- +(0,-.75) node[\colourmakedark{4},above left] {$\omega$};
\fill[black] (c) circle (.3pt);
\node[black,yshift=-.4cm,left] at (spath cs:\ttbranch{a}{2} 1) {$\tau$};
\node[\colourmakedark{3},left] at (spath cs:\ttbranch{b}{1} 1) {$\tau'$};
\end{tikzpicture}

%% file: figures-source/one-switch-fundamental-curves.tex
\begin{tikzpicture}
\tikzset{
tt/new=a,
tt/a/.cd,
set sep=.8cm,
set spread=1cm,
add switch={0,0},
add bunch={1}{left}{3},
add bunch={1}{right}{3,spread=1.5cm},
finalise,
}
\path[spath/save=branch 1] (a-1) to[out=180,in=180,looseness=2] (0,1.25) to[out=0,in=0,looseness=3.5] (a-1);
\path[spath/save=branch 2] [tt/a/use branch=2] (spath cs:{\ttbranch{a}{2}} 1) to[out=180,in=180,looseness=2.5]  (spath cs:{\ttbranch{a}{3}} 1) [spath/use={\ttbranch{a}{3},weld,reverse}];
\path[spath/save=branch 3] [tt/a/use branch=4] (spath cs:{\ttbranch{a}{4}} 1) to[out=0,in=0,looseness=2] (spath cs:{\ttbranch{a}{6}} 1) [spath/use={\ttbranch{a}{6},weld,reverse}];
\tikzset{
tt/new=b,
tt/b/.cd,
set width=.5cm,
add switch={0,0},
add branch={[spath/use=branch 1]},
add branch={[spath/use=branch 2]},
add branch={[spath/use=branch 3]},
set switch thickness={1}{.8cm},
set branch weights={1=1,2=1,3=1},
finalise,
carried={
    name=c,
    1=1,
},
carried={
    name=d,
    2=1,3=1,
}
}
\fill[grey background] [tt/b/use branch background all but=1,tt/b/use switch background all];
\draw[line,black] [tt/b/use branch contour all but=1,tt/b/use switch contour all];
\draw[line,palette 2] [tt/d/use branch all];
\fill[grey background] [tt/b/use branch background=1];
\draw[line,black] [tt/b/use branch contour=1];
\draw[line over={grey background},palette 1] [tt/c/use branch all];
\node[above left,\colourmakedark{2},yshift=4pt] at (spath cs:{\ttbranch{d}{1}} .7) {$c$};
\node[\colourmakedark{1},above left,xshift=-2pt] at (spath cs:{\ttbranch{c}{1}} .4) {$a$};
\end{tikzpicture}

%% file: figures-source/deeply-nested-example.tex
\begin{tikzpicture}
\tikzset{
tt/new=a,
tt/a/.cd,
set sep=.2cm,
set width=1cm,
set spread=1cm,
add switch={0,0},
add bunch={1}{left}{3},
add bunch={1}{right}{2},
set switch thickness={1}{2cm},
set branch weights={1=2,2=2,3=2,4=3,5=3},
finalise,
carried={
    name=b,
    switch={1}{.075cm}{3,4}{3,4,5},
    1=1,2=4,3=2,4=5,5=3,
}
}
\fill[grey background] [tt/a/use branch background all but={2,4},tt/a/use switch background all];
\fill[\colourmakebackground{2}] [tt/a/use branch background={2}];
\fill[\colourmakebackground{3}] [tt/a/use branch background={4}];
\draw[black,line] [tt/a/use branch contour all but={2,4},tt/a/use switch contour all];
\draw[\colourmakedark{2},line] [tt/a/use branch contour={2}];
\draw[\colourmakedark{3},line] [tt/a/use branch contour={4}];
\draw[line,palette 1] [tt/b/use branch all];
\pic[palette 1] at (b-1) {dot};
\node[\colourmakedark{1},left,xshift=-2pt] at (spath cs:{\ttbranch{b}{1}} 1) {$\tau'$};
\node[black,left,xshift=-2pt,yshift=-4pt] at (spath cs:{\ttbranch{a}{3}} 1) {$\tau$};
\node[\colourmakedark{2},left,xshift=-2pt] at (spath cs:{\ttbranch{a}{2}} 1) {$b_1^\bullet$};
\node[\colourmakedark{3},right,xshift=2pt] at (spath cs:{\ttbranch{a}{4}} 1) {$b_2^\bullet$};
\end{tikzpicture}

%% file: figures-source/deeply-nested-diagonal-extension.tex
\begin{tikzpicture}
\tikzset{
tt/new=a,
tt/a/.cd,
set sep=.2cm,
set width=1cm,
set spread=1cm,
add switch={0,0},
add bunch={1}{left}{3},
add bunch={1}{right}{2},
set switch thickness={1}{2cm},
set branch weights={1=2,2=2,3=2,4=3,5=3},
finalise,
carried={
    name=b,
    switch={1}{.3cm}{}{1,2,3},
    4=3
}
}
\fill[grey background] [tt/a/use branch background all but=4,tt/a/use switch background all];
\fill[\colourmakebackground{2}] [tt/a/use branch background={4}];
\draw[black,line] [tt/a/use branch contour all but=4,tt/a/use switch contour all];
\draw[\colourmakedark{2},line] [tt/a/use branch contour={4}];
\draw[line,palette 1] [tt/b/use branch={1,3}];
\draw[line,palette 3] [tt/b/use branch=2];
\pic[palette 1] at (b-1) {dot};
\node[left,xshift=-2pt,\colourmakedark{1}] at (spath cs:\ttbranch{b}{3} .25) {$e'_1$};
\node[below,\colourmakedark{1}] at (spath cs:\ttbranch{b}{1} .2) {$e'_2$};
\node[right,xshift=2pt,\colourmakedark{3}] at (spath cs:\ttbranch{b}{2} 1) {$b'_{\bullet}$};
\node[above,yshift=.5cm,\colourmakedark{2}] at (spath cs:\ttbranch{a}{4} .5) {$b^\bullet_i$};
\end{tikzpicture}

%% file: distance-algorithm.tex
\section{Coarse distance algorithm}
\label{sec:distance algorithm}

For a train track $\tau$ on a surface $S$, the set $\fundamental(\tau)$ of fundamental curves of $\tau$ has universally bounded diameter in $\C(S)$, as it can be easily seen from \cref{thm:distance bound from intersections}.
Therefore, we can think of train tracks on $S$ as ``coarse points'' in $\C(S)$, and interpret sequences of train tracks as ``coarse paths'' in $\C(S)$.

\citeauthor{masur-minsky-quasiconvexity} show in \cite[Theorem 1.3]{masur-minsky-quasiconvexity} that a splitting sequence $\tau_0\carries\ldots\carries\tau_n$ of train tracks, if considered as a coarse path in $\C(S)$, is an unparametrised quasi\=/geodesic, and in fact lies uniformly close to an actual geodesic.
\textcite{bell-webb-algorithms} use this fact to construct a finite neighbourhood of this quasi\=/geodesic that is guaranteed to contain a (tight) geodesic between a fundamental curve of $\tau_0$ and a fundamental curve of $\tau_n$.
Since the size of this neighbourhood is polynomial in $n$, they employ this construction to give a polynomial\=/time algorithm to compute distances in the curve graph.

The algorithm of \citeauthor{bell-webb-algorithms}, however, does not run in polynomial time in the complexity of $S$.
In this section, we present an algorithm that offers a trade\=/off between accuracy and good dependence of the running time on $\xi$.
More precisely, given two essential curves $a$ and $b$ on $S$, our algorithm -- as stated in \cref{thm:distance algorithm} -- computes an estimate $d$ of the distance $\dist[\C(S)](a,b)$ such that
\[
d-\constant{curves l+}\le\dist[\C(S)](a,b)\le\constant{carried curves l*}\cdot d+\constant{curves l+}.
\]
The running time of the algorithm is polynomial in the parameters $\norm{a}$, $\norm{b}$, and $\xi$, and the constants $\constant{carried curves l*}$ and $\constant{curves l+}$ are polynomial in $\xi$.

The core idea of our algorithm is to use the combinatorial criterion of \cref{thm:nesting lemma} to detect when the unparametrised quasi\=/geodesic given by $\tau_0\carries\ldots\carries\tau_n$ has made definite progress in $\C(S)$, and extract a subsequence corresponding to an actual quasi\=/geodesic parametrisation.
The issue with this strategy is that, as already noted, the condition $\PN(\tau')\subs\intPN(\tau)$ is hard to verify in polynomial time in $\xi$ for an arbitrary pair $\tau'\carried\tau$ of train tracks.
To overcome this problem, we devise a procedure to construct a sequence of nested train tracks for which the condition above can be efficiently checked; as anticipated, the key property is that consecutive elements of the sequence should be deeply nested.

\subsection{\texorpdfstring{Agol\==Hass\==Thurston}{Agol-Hass-Thurston} moves}
\label{sec:distance algorithm:aht moves}

\textcite{agol-hass-thurston} introduce an algorithm to compute orbits of interval isometries in polynomial time; the algorithm works in complete generality, but the main application in their article is to efficiently compute properties of normal surfaces in $3$\=/manifolds.
We describe a specialised version of the Agol\=/Hass\=/Thurston algorithm for measured train tracks.
This algorithm is simpler than the general one, in the sense that it involves only two types of moves, which can easily be interpreted in a purely topological fashion.
The increased simplicity provides better running time bounds, and allows us to reason about combinatorial properties of the resulting train tracks.


We start by describing the two aforementioned moves.
The input for both of them is a \fullymeasured{} train track $(\tau,\mu)$ -- where we assume that $\tau$ is not a simple closed curve -- and a switch $s$ of $\tau$, together with a co\=/orientation $\omega$ of $s$.

\step{The \splitmove{} move.} Choose a direction $\eta$ of $s$ such that
\[
    \mu(\listofends{\tau}{\eta}[\omega]_1)\le\mu(\listofends{\tau}{\bar{\eta}}[\omega]_1).
\]
For ease of notation, let $e_i=\listofends{\tau}{\eta}[\omega]_i$. Moreover, denote by $b_i$ the branch of which $e_i$ is an end, and by $\bar{b}$ the branch of which $\listofends{\tau}{\bar{\eta}}[\omega]_1$ is an end; we call this end $\bar{b}^\bullet$.
Define
\[
    h=\max\{1\le i\le\card{\listofends{\tau}{\eta}}:\mu(b_1)+\ldots+\mu(b_i)\le\mu(\bar{b})\}.
\]
We call the branches $b_1,\ldots,b_h$ \emph{splitting branches}, and the ends $e_1,\ldots,e_h$ \emph{splitting ends}.

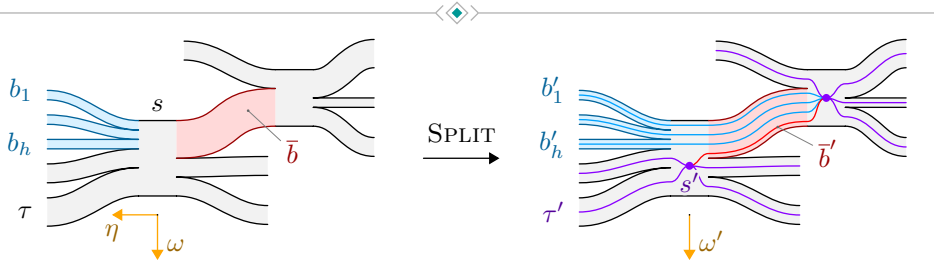
\begin{figure}
\centering
    \tikzsetnextfilename{split-move}%
    \input{figures-source/split-move.tex}%

\caption{The effect of applying a \splitmove{} move to a measured train track $(\tau,\mu)$ with switch $s\in\switches(\tau)$ and co\=/orientation $\omega$ of $s$. The measured train track $(\tau,\mu)$ is represented by a thickened neighbourhood of $\tau$, where the thickness of the branches is proportional to the weight assigned by $\mu$. On the right, we only draw the train track $\tau'\carried\tau$, ignoring the measure $\mu'$.}
\label{fig:split move}
\end{figure}

The output of the \splitmove{} move is a measured train track $(\tau',\mu')\carried(\tau,\mu)$, obtained by pushing the ends $e_1,\ldots,e_h$ along $\bar{b}$ all the way to the other end $\bar{b}_\bullet$.
\Cref{fig:split move} illustrates the construction of $\tau'$, and shows how
\[
    \tau'=\tau\setminus(b_1\cup\ldots\cup b_h)\cup(b_1'\cup\ldots\cup b_h'),
\]
where $b_i'$ is a branch of $\tau'$ that follows $b_i\cup\bar{b}$ very closely. 
The carrying $\tau'\carried\tau$ satisfies the following properties: every non\=/splitting branch of $\tau$ is untouched in $\tau'$, and for every splitting branch $b_i$ of $\tau$ there is a unique branch $b_i'$ of $\tau'$ such that $\carrymap[\tau'][\tau](\one{b_i'})\ge\one{b_i}$.
In fact, it is easy to see that $\carrymap[\tau'][\tau](\one{b_i'})$ is either $\one{b_i}+\one{\bar{b}}$ or $\one{b_i}+2\cdot\one{\bar{b}}$, depending on whether $b_i$ has one or two splitting ends.

The measure $\mu'$ is chosen so that $\carrymap[\tau'][\tau](\mu')=\mu$, and is obtained from $\mu$ as follows:
\begin{itemize}
\item if $b$ is a non\=/splitting branch of $\tau$ different from $\bar{b}$, then $\mu'(b')=\mu(b)$, where $b'$ is the branch of $\tau'$ such that $\carrymap[\tau'][\tau](\one{b'})=\one{b}$;
\item $\mu'(\bar{b}')=\mu(\bar{b})-\mu(b_1)-\ldots-\mu(b_h)$, where $\bar{b}'$ is the branch of $\tau'$ such that $\carrymap[\tau'][\tau](\one{\bar{b}'})=\one{\bar{b}}$;
\item $\mu'(b_i')=\mu(b_i)$ for every $1\le i\le h$.
\end{itemize}

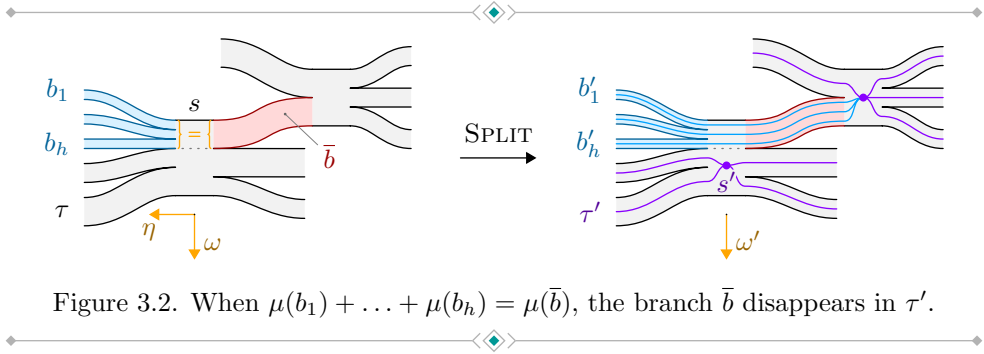
\begin{figure}
\centering
    \tikzsetnextfilename{split-move-branch-disappear}%
    \input{figures-source/split-move-branch-disappear.tex}%

\caption{When $\mu(b_1)+\ldots+\mu(b_h)=\mu(\bar{b})$, the branch $\bar{b}$ disappears in $\tau'$.}
\label{fig:split move branch disappear}
\end{figure}

There is one exceptional case in which the construction is slightly different.
If $\mu(\bar{b})=\mu(b_1)+\ldots+\mu(b_h)$, then the measure $\mu'$ would assign weight $0$ to the branch $\bar{b}'$ of $\tau'$.
To ensure that the measured train track $(\tau',\mu')$ is \fullymeasured{}, we remove the branch $\bar{b}'$ from $\tau'$, with the understanding that in this case the branch $\bar{b}$ of $\tau$ is not untouched in $\tau'$.
In this case, we say that the branch $\bar{b}$ \emph{disappears} in $\tau'$; the situation is depicted in \cref{fig:split move branch disappear}. 
Moreover, when $\mu(\bar{b})=\mu(s)$ -- that is, when $\bar{b}$ is the only branch incident to $s$ in the direction $\bar{\eta}$ -- not only does $\bar{b}$ disappear in $\tau'$, but there is no switch of $\tau'$ mapped to $s$ by $\carrymap[\tau'][\tau]$.
If this happens, we say that the switch $s$ \emph{disappears} in $\tau'$.

Without delving into implementation details, we remark that, from the description we have provided, it is easy to algorithmically construct the \fullymeasured{} train track $(\mu',\tau')$ and the carrying map $\carrymap[\tau'][\tau]$.
We observe that the complexity of the carrying map is $\norm{\carrymap[\tau'][\tau]}=\bigO(\xi)$, and the measure $\mu'$ satisfies $\norm{\mu'}\le\norm{\mu}$.
As long as a reasonable data structure is used to implement train tracks, the construction can be performed with running time $\bigO(\poly(\xi)\cdot\norm{\mu})$.

\step{The \twistmove{} move.} This move is an optimisation to perform multiple \splitmove{} moves in a row.
In the setting above, suppose that $h<\card{\listofends{\tau}{\eta}}$ and $\listofends{\tau}{\eta}[\omega]_{h+1}=\bar{b}_\bullet$.
We still call the branches $b_1,\ldots,b_h$ \emph{splitting branches}, and the ends $e_1,\ldots,e_h$ \emph{splitting ends}.
The output of the \twistmove{} move is the measured train track $(\tau',\mu')\carried(\tau,\mu)$ obtained from $k$ successive applications of the \splitmove{} move, where
\[
    k=\left\lfloor\frac{\mu(\bar{b})}{\mu(b_1)+\ldots+\mu(b_h)}\right\rfloor.
\]

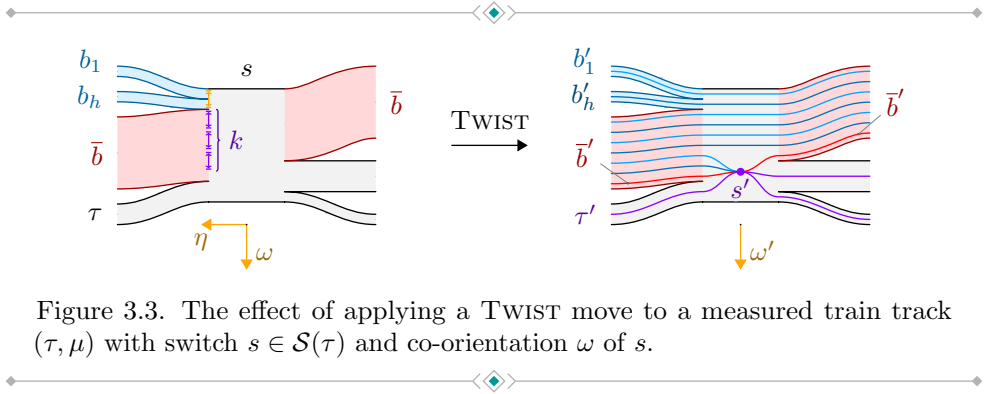
\begin{figure}
\centering
    \tikzsetnextfilename{twist-move}%
    \input{figures-source/twist-move.tex}%

\caption{The effect of applying a \twistmove{} move to a measured train track $(\tau,\mu)$ with switch $s\in\switches(\tau)$ and co\=/orientation $\omega$ of $s$.}
\label{fig:twist move}
\end{figure}

To perform this move efficiently, we first note that the train track $\tau'$ is the result of applying $\twist_{\bar{b}}^{\pm k}$ to $\tau$ -- which justifies the name ``\twistmove{}''; the construction is depicted in \cref{fig:twist move}.
Just like in the \splitmove{} case, non\=/splitting branches of $\tau$ are untouched in $\tau'$, and for every splitting branch $b_i$ of $\tau$ there is a unique branch $b_i'$ of $\tau'$ such that $\carrymap[\tau'][\tau](\one{b_i'})\ge\one{b_i}$.
More precisely, the pre\=/measure $\carrymap[\tau'][\tau](\one{b_i'})$ on $\tau$ is either $\one{b_i}+k\cdot\one{\bar{b}}$ or $\one{b_i}+2k\cdot\one{\bar{b}}$, depending on whether $b_i$ has one or two splitting ends.

The measure $\mu'$ on $\tau'$ is defined exactly like in the \splitmove{} case, the only difference being the formula
\begin{equation}\label{eqn:twist move formula}
    \mu'(\bar{b}')=\mu(\bar{b})-k\cdot(\mu(b_1)+\ldots+\mu(b_h)).
\end{equation}
Just like in the \splitmove{} case, if the measure $\mu'$ would assign weight $0$ to the branch $\bar{b}'$, we remove said branch from $\tau'$, and say that $\bar{b}$ \emph{disappears} in $\tau'$.
In this case, in particular, the branch $\bar{b}$ will not be untouched in $\tau'$.
Note that the switch $s$ never disappears in $\tau'$; in other words, there is a unique switch $s'$ of $\tau'$ such that $\carrymap[\tau'][\tau](s')=s$.

The discussion above gives a complete recipe to construct the measured train track $(\tau',\mu')$ and the carrying map $\carrymap[\tau'][\tau]$.\footnote{In fact, since $\tau$ and $\tau'$ differ by a homeomorphism of $S$, they are represented by the same abstract graph.}
We observe that the complexity of the carrying map is $\norm{\carrymap[\tau'][\tau]}=\bigO(\xi\cdot\log k)=\bigO(\xi\cdot\norm{\mu})$, and the measure $\mu'$ satisfies $\norm{\mu'}\le\norm{\mu}$.
Moreover, the \twistmove{} move can be performed with running time
\[
    \bigO(\poly(\xi)\cdot\norm{\mu}\cdot\log\norm{\mu}),
\]
using the algorithm of \textcite{harvey-multiplication} to compute the integers $k$ and $\mu'(\bar{b}')$ in time $\bigO(\norm{\mu}\cdot\log\norm{\mu})$.

\step{Proximity.} The first elementary property of the \splitmove{} and \twistmove{} moves we observe is that they do not move train tracks too far in $\C(S)$. More precisely, we have the following bound.

\newconstant{proximity d}{D_{\mathrm{prox}}}{Dprox}{17}
\begin{proposition}[Proximity]\label{thm:proximity}
If $(\tau',\mu')$ is obtained from $(\tau,\mu)$ by a \splitmove{} move or a \twistmove{} move, then
\[
    \diamtt{\tau}{\tau'}\le\constant{proximity d},
\]
where
\[
    \declareconstant{proximity d}=\constantvalue{proximity d}.
\]
\end{proposition}
\begin{proof}
Let us employ the same notation we have used above when defining the two moves.
\begin{substeps}
\item Suppose first that $(\tau',\mu')$ is obtained from $(\tau,\mu)$ by a \splitmove{} move.
Let $\nu$ and $\nu'$ be fundamental measures on $\tau$ and $\tau'$ respectively.
It is easy to see that $\carrymap[\tau'][\tau](\nu')(z)\le 2$ for every switch $z\in\switches(\tau)\setminus\{s\}$, while
\[
    \carrymap[\tau'][\tau](\nu')(s)=\nu'(s')+\nu'(b_1')+\ldots+\nu'(b_h')\le4,
\]
where $s'$ is the switch of $\tau'$ such that $\carrymap[\tau'][\tau](s')=s$ -- or, say, $\nu'(s')=0$ if $s$ disappears in $\tau'$. As a consequence, we have the inequalities
\[
    i(\nu,\nu')\le\sum_{z\in\switches(\tau)}\nu(z)\cdot\carrymap[\tau'][\tau](\nu')(z)\le 4\card{\switches(\tau)}+4\le 16\xi+4,
\]
which by \cref{thm:distance bound from intersections} (and some easy algebra) imply that $\dist[\C(S)](\nu,\nu')<17$.
\item Suppose now that $(\tau',\mu')$ is obtained from $(\tau,\mu)$ by a \twistmove{} move.
Note that the branches $\bar{b}$ and $\bar{b}'$ are isotopic simple closed curves.
Moreover, every fundamental curve of $\tau$ (respectively $\tau'$) intersects $\bar{b}$ (respectively $\bar{b}'$) at most twice.
This immediately gives the upper bound
\[
    \diamtt{\tau}{\tau'}\le\diam[\C(S)](\fundamental(\tau),\bar{b})+\diam[\C(S)](\bar{b}',\fundamental(\tau'))\le 2+2=4.\qedhere
\]
\end{substeps}
\end{proof}

\step{AHT sequences.} By iteratively applying \splitmove{} and \twistmove{} moves, we obtain a sequence of nested measured train tracks.
More precisely, let $(\tau,\mu)$ be a \fullymeasured{} train track, $s$ a switch of $\tau$, and $\omega$ a co\=/orientation of $s$.
The \emph{Agol\==Hass\==Thurston sequence} (or \emph{AHT sequence} for short) generated by this data is a sequence of measured train tracks obtained from the following procedure.

We start by setting  $(\tau_0,\mu_0)=(\tau,\mu)$, $s_0=s$, and $\omega_0=\omega$.
On the $i$\=/th step, for $i\ge 0$, consider the measured train track $(\tau_i,\mu_i)$.
If $\tau_i$ is a simple closed curve, then the procedure terminates.
If we can apply a \twistmove{} move to the switch $s_i$ with co\=/orientation $\omega_i$, then we do so, denoting by $(\tau_{i+1},\mu_{i+1})$ the resulting measured train track.
Otherwise, we apply a \splitmove{} move to the same co\=/oriented switch, and we call the resulting measured train track $(\tau_{i+1},\mu_{i+1})$.
If the switch $s_i$ disappears in $\tau_{i+1}$, then we set $s_{i+1}=\emptyset$ and $\omega_{i+1}=\emptyset$, and the procedure terminates.
Otherwise, we let $s_{i+1}$ be the unique switch of $\tau_{i+1}$ such that $\carrymap[\tau_{i+1}][\tau_i](s_{i+1})=s_i$, and $\omega_{i+1}$ be the orientation of $s_{i+1}$ such that $\carrymap[\tau_{i+1}][\tau_i](\omega_{i+1})=\omega_i$.
We then proceed to the $(i+1)$\=/th step.

This procedure will terminate after a finite number of steps, as \cref{thm:exponential decay} below implies. The output is a sequence of nested \fullymeasured{} train tracks
\[
(\tau_0,\mu_0)\carries(\tau_1,\mu_1)\carries\ldots\carries(\tau_n,\mu_n),
\]
where $(\tau_0,\mu_0)=(\tau,\mu)$, and either $\tau_n$ is the curve represented by $\mu$, or $\tau_n$ has one less switch than $\tau$.
Note that, since the measured train tracks are generated one by one, we can stop the procedure prematurely whenever a condition of our choice is satisfied.
The first $m$ elements of an AHT sequence can be computed with running time
\[
\bigO(\poly(\xi)\cdot m\cdot\norm{\mu}\cdot\log\norm{\mu}).
\]

There are a number of interesting objects associated to each measured train track $(\tau_i,\mu_i)$ in an AHT sequence.
We have already defined the switch $s_i$ and the co\=/orientation $\omega_i$ in the description above; we now introduce the following additional notation for $0\le i<n$:
\begin{itemize}
\item we write $\carrymap_{j,i}$ for $\carrymap[\tau_j][\tau_i]$ for every $i\le j\le n$;
\item $\eta_i$ is the direction of $s_i$ such that the end $\listofends{\tau_i}{\eta_i}[\omega_i]_1$ is splitting in $\tau_i$; in particular, it satisfies
\[
    \mu_i(\listofends{\tau_i}{\eta_i}[\omega_i]_1)\le\mu_i(\listofends{\tau_i}{\bar{\eta_i}}[\omega_i]_1);
\]
\item $\bar{b}_i$ is the branch of $\tau_i$ of which $\listofends{\tau_i}{\bar{\eta_i}}[\omega_i]_1$ is an end;
\item if the branch $\bar{b}_i$ disappears in $\tau_{i+1}$, then we set $\bar{b}_i'=\emptyset$; otherwise, we let $\bar{b}_i'$ be the unique branch of $\tau_{i+1}$ such that $\carrymap_{i+1,i}(\one{\bar{b}_i'})=\one{\bar{b}_i}$.
\end{itemize}

\subsection{Elementary properties of AHT sequences}

Let
\[
(\tau_0,\mu_0)\carries(\tau_1,\mu_1)\carries\ldots\carries(\tau_n,\mu_n)
\]
be the AHT sequence generated by a \fullymeasured{} train track $(\tau,\mu)$, a switch $s\in\switches(\tau)$, and a co\=/orientation $\omega$ of $s$.

\begin{figure}
\centering
    \tikzsetnextfilename{aht-sequence-typical-carrying}%
    \input{figures-source/aht-sequence-typical-carrying.tex}%

\caption{A typical carrying $\tau_i\carried\tau_0$ between train tracks in an AHT sequence.}
\label{fig:aht sequence typical carrying}
\end{figure}
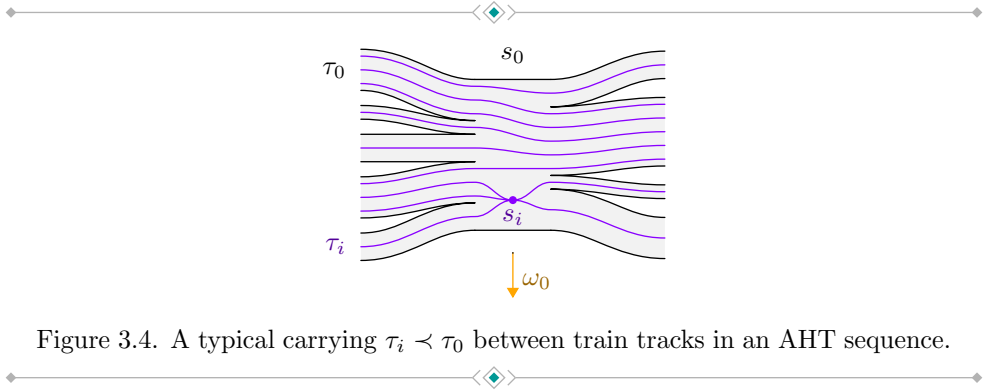

\step{Carrying properties.}
Let us fix an integer $0\le i\le n$.
\Cref{fig:aht sequence typical carrying} shows what the typical local picture of the carrying $\tau_i\carried\tau_0$ looks like near $s$.

We now list a few properties of this carrying of train tracks that follow immediately from the procedure for generating the AHT sequence.
\begin{enumarabic}
\item\label[carrying property]{itm:carrying properties:untouched branch} If a branch $b_0\in\branches(\tau_0)$ is untouched in $\tau_i$, then there is a sequence of branches $b_1\in\branches(\tau_1)$, \dots, $b_i\in\branches(\tau_i)$ such that $\carrymap_{j,j-1}(\one{b_j})=\one{b_{j-1}}$ for every $1\le j\le i$.
Moreover, the branches $b_1,\ldots,b_i$ are unique, and $b_0,\ldots,b_{i-1}$ are not splitting in $\tau_0,\ldots,\tau_{i-1}$ respectively.
\item\label[carrying property]{itm:carrying properties:branches subset} If $B$ is any subset of $\branches(\tau_i)$, then\todo{Used only for $i=1$?}
\[
    \card{\{b\in\branches(\tau_0):\text{$\carrymap_{i,0}(\one{b'})\ge\one{b}$ for some $b'\in B$}\}}\ge\card{B}.
\]
\item\label[carrying property]{itm:carrying properties:height inequality} Suppose that $\tau_0$ is one\=/switch. If $e'$ is an end of $\tau_i$ and $e=\carrymap_{i,0}(e')$, then
\begin{align*}
    \height{\omega_0}(e)\le\height{\omega_i}(e')&&\text{and}&&\height{\omega_i}(e')+\mu_i(e')\le\height{\omega_0}(e)+\mu_0(e).
\end{align*}
\item\label[carrying property]{itm:carrying properties:persistency condition} Suppose that $\tau_0$ is one\=/switch. An end $e$ of $\tau_0$ is persistent in $\tau_i$ if and only if
\[
    \height{\omega_0}(e)<\mu_i(s_i).
\]
In particular, for every direction $\eta$ of $s_0$, the set of ends in $\listofends{\tau_0}{\eta}$ that are persistent in $\tau_i$ is a final segment of $\listofends{\tau_0}{\eta}[\omega_0]$.
\item\label[carrying property]{itm:carrying properties:untouched height} Suppose that $\tau_0$ is one\=/switch. Let $b$ be a branch of $\tau_0$ that is untouched in $\tau_i$. Suppose that the ends $b^\bullet$ and $b_\bullet$ of $b$ satisfy
\[
    \height{\omega_0}(b^\bullet)>\height{\omega_0}(b_\bullet).
\]
Let $b'$ be the (unique by \cref{itm:carrying properties:untouched branch}) branch of $\tau_i$ such that $\carrymap_{i,0}(\one{b'})=\one{b}$, and $b'_\bullet$ be the end of $b'$ such that $\carrymap_{i,0}(b'_\bullet)=b_\bullet$. Then
\[
    \height{\omega_i}({b'}^\bullet)>\height{\omega_i}(b'_\bullet),
\]
and ${b'}^\bullet$ is the unique end of $\tau_i$ whose image under $\carrymap_{i,0}$ is $b^\bullet$.
\end{enumarabic}

\Cref{itm:carrying properties:persistency condition} can be proved as follows. If $e$ is persistent in $\tau_i$, then by \cref{itm:carrying properties:height inequality} we have
\[
    \height{\omega_0}(e)\le\height{\omega_i}(e')<\mu_i(s_i),
\]
where $e'$ is any end of $\tau_i$ such that $\carrymap_{i,0}(e')=e$.
If, conversely, the end $e$ satisfies $\height{\omega_0}(e)<\mu_i(s_i)$, then there is a unique end $e'$ of $\tau_i$ such that $\direction(e)=\direction(\carrymap_{i,0}(e'))$ and
\[
    \height{\omega_i}(e')\le\height{\omega_0}(e)<\height{\omega_i}(e')+\mu_i(e').
\]
Combining these inequalities with \cref{itm:carrying properties:height inequality} yields
\[
    \height{\omega_0}(\carrymap_{i,0}(e'))\le\height{\omega_0}(e)<\height{\omega_0}(\carrymap_{i,0}(e'))+\mu_0(\carrymap_{i,0}(e')),
\]
which necessarily implies that $\carrymap_{i,0}(e')=e$.

The other \namecrefs{itm:carrying properties:height inequality} are immediate when $i=1$, and can be proved by induction on $i$ in the general case.

\step{Exponential decay.}
We now address the issue of how long an AHT sequence can be.
It is quite clear that it has to be finite, since the quantity
\[
    \sum_{b\in\branches(\tau_i)}\mu_i(b)
\]
is a positive integer and strictly decreasing in $i$.
However, the following proposition gives a much better bound than this first naive estimate.

\newconstant{exponential decay n}{N_{\mathrm{dec}}}{Ndec}{6\xi+1}
\begin{proposition}[Exponential decay]\label{thm:exponential decay}
Either $n\le\constant{exponential decay n}$, or there is an integer $1\le k\le\constant{exponential decay n}$ such that $\norm{\mu_k}\le\norm{\mu_0}-1$, where
\[
    \declareconstant{exponential decay n}=\constantvalue{exponential decay n}.
\]
\end{proposition}
\begin{proof}
Let us assume that $n>\constant{exponential decay n}$; in particular, $s_i\neq\emptyset$ for every $0\le i\le\constant{exponential decay n}$.
\begin{substeps}
\item Suppose that $\bar{b}_i'=\emptyset$ for some $0\le i<\constant{exponential decay n}$. Then
\[
    \norm{\mu_{i+1}}=\norm{\mu_i}-\log(\mu_i(\bar{b}_i)+1)\le\norm{\mu_i}-1,
\]
and we can take $k=i+1$.
Consequently, from now on, we will assume that $\bar{b}_i'\neq\emptyset$ for every $0\le i<\constant{exponential decay n}$.
\item Note that, if $2\mu_{i+1}(\bar{b}_i')<\mu_i(\bar{b}_i)$ for some $0\le i<\constant{exponential decay n}$, then $\norm{\mu_{i+1}}\le\norm{\mu_i}-1$, and we can take $k=i+1$.
Therefore, from now on, we will assume that $2\mu_{i+1}(\bar{b}_i')\ge\mu_i(\bar{b}_i)$ for every $0\le i<\constant{exponential decay n}$.
\item Suppose now that $(\tau_{i+1},\mu_{i+1})$ is obtained from $(\tau_i,\mu_i)$ by a \twistmove{} move for some $0\le i<\constant{exponential decay n}$.
Then \cref{eqn:twist move formula} implies that $\mu_{i+1}(\bar{b}_i')$ is the remainder of $\mu_i(\bar{b}_i)$ modulo a positive integer smaller than $\mu_i(\bar{b}_i)$ -- namely, the sum of the weights of the splitting ends of $\tau_i$.
In particular, we have that $2\mu_{i+1}(\bar{b}_i')<\mu_i(\bar{b}_i)$, which contradicts our previous assumption.
Therefore, from now on, we will assume that $(\tau_{i+1},\mu_{i+1})$ is obtained  from $(\tau_i,\mu_i)$ by a \splitmove{} move for every $0\le i<\constant{exponential decay n}$.
\item Let us now fix an integer $0\le i<\constant{exponential decay n}-1$.
Observe that, by definition of \splitmove{} move, the branch $\bar{b}_i'$ is splitting in $\tau_{i+1}$.
By considering the \splitmove{} move from $\tau_{i+1}$ to $\tau_{i+2}$, we see that
\[
    \mu_{i+2}(\bar{b}_{i+1}')\le\mu_{i+1}(\bar{b}_{i+1})-\mu_{i+1}(\bar{b}_i').
\]
Combining this inequality with $2\mu_{i+1}(\bar{b}_i')\ge\mu_i(\bar{b}_i)$ and $2\mu_{i+2}(\bar{b}_{i+1}')\ge\mu_{i+1}(\bar{b}_{i+1})$ yields
\begin{equation}\label{eqn:exponential decay:increasing weight}
    \mu_i(\bar{b}_i)\le\mu_{i+1}(\bar{b}_{i+1})\qquad\text{for every $0\le i<\constant{exponential decay n}-1$}.
\end{equation}

The intuition for why \cref{eqn:exponential decay:increasing weight} gives a contradiction is that, since $\tau_0$ has at most $6\xi<\constant{exponential decay n}$ branches, two of the first $\constant{exponential decay n}$ \splitmove{} moves will operate on the same branch.
More precisely, we claim the following: there are two indices $0\le i<j<\constant{exponential decay n}$ such that $\carrymap_{j,i}(\one{\bar{b}_j})\ge\one{\bar{b}_i}$ as pre\=/measures on $\tau_i$.
In fact, if this were not the case, an easy induction using \cref{itm:carrying properties:branches subset} would show that the cardinality of the set
\[
\{b\in\branches(\tau_i):\text{$\carrymap_{j,i}(\one{\bar{b}_j})\ge\one{b}$ for some $i\le j<\constant{exponential decay n}$}\}
\]
is at least $\constant{exponential decay n}-i$ for $0\le i<\constant{exponential decay n}$, which is of course impossible for $i=0$.

Let then $i<j$ be two indices as above.
Since $\carrymap_{j,i}(\mu_j)=\mu_i$, we immediately get that $\mu_j(\bar{b}_j)\le\mu_i(\bar{b}_i)$, which we can refine to a strict inequality via the following argument.
If $\carrymap_{j,i+1}(\one{\bar{b}_j})\ge\one{\bar{b}_i'}$, then $\mu_j(\bar{b}_j)\le\mu_{i+1}(\bar{b}_i')<\mu_i(\bar{b}_i)$; otherwise, $\mu_i(\bar{b}_i)\ge\mu_{i+1}(\bar{b}_i')+\mu_j(\bar{b}_j)$.
Either way, the strict inequality $\mu_j(\bar{b}_j)<\mu_i(\bar{b}_i)$ contradicts \cref{eqn:exponential decay:increasing weight}.\qedhere
\end{substeps}
\end{proof}
As a corollary, we observe that
\begin{align}\label{eqn:aht sequence length}
    i<\constant{exponential decay n}\cdot(\norm{\mu_0}-\norm{\mu_i}+1)&&\text{for every $0\le i\le n$.}
\end{align}

\subsection{Deep nesting algorithm}

Our aim is now to describe an algorithm that, given as input a measured train track $(\tau,\mu)$, produces a measured train track $(\tau',\mu')\carried(\tau,\mu)$ such that $\tau'$ is deeply nested in $\tau$ and uniformly close to it in $\C(S)$.
This algorithm will play a crucial role in our proof of \cref{thm:distance algorithm for carried curves}.
We start with a preliminary lemma.

\begin{lemma}\label{thm:end of non-untouched branch is non-persistent}
Let $(\tau,\mu)$ be a one\=/switch \fullymeasured{} train track that is not a simple closed curve, and let
\[
(\tau_0,\mu_0)\carries\ldots\carries(\tau_n,\mu_n)
\]
be the AHT sequence generated by $(\tau,\mu)$, the unique switch $s$ of $\tau$, and a co\=/orientation $\omega$ of $s$.
Let $b$ be a branch of $\tau$, and let $b^\bullet$ be the end of $b$ such that
\[
\height{\omega}(b^\bullet)>\height{\omega}(b_\bullet).
\]
Suppose that $b$ is untouched in $\tau_{i-1}$ but not in $\tau_i$ for some $1\le i\le n$. Then the end $b^\bullet$ is not persistent in $\tau_i$.
\end{lemma}
\begin{proof}
Let $b'$ be the branch of $\tau_{i-1}$ such that $\carrymap_{i-1,0}(\one{b'})=\one{b}$ (which is unique by \cref{itm:carrying properties:untouched branch}), and let $b'_\bullet$ be the end of $b'$ such that $\carrymap_{i-1,0}(b'_\bullet)=b_\bullet$.
\Cref{itm:carrying properties:untouched height} implies that
\[
\height{\omega_{i-1}}({b'}^\bullet)>\height{\omega_{i-1}}(b'_\bullet).
\]
There are two possible reasons for $b'$ to not be untouched in $\tau_i$: either $b'=\bar{b}_{i-1}$ and it disappears in $\tau_i$, or the end ${b'}^\bullet$ is splitting in $\tau_{i-1}$.
In both cases, one checks that the end ${b'}^\bullet$ is not persistent in $\tau_i$; since ${b'}^\bullet$ is the only end of $\tau_{i-1}$ that $\carrymap_{i-1,0}$ sends to $b^\bullet$, this implies that $b^\bullet$ is not persistent in $\tau_i$.
\end{proof}

We now present an intermediate algorithm that, loosely speaking, takes the input train track ``one step closer'' to being deeply nested.

\newconstant{deep nesting one step d}{D'_{\mathrm{d.n}}}{Ddn'}{\constant{proximity d}+10}
\begin{proposition}\label{thm:deep nesting one step}
Let $(\tau,\mu)$ be a one\=/switch \fullymeasured{} train track, whose unique switch we denote by $s$.
Let $e$ and $e'$ be two consecutive ends of $\tau$.
There is an algorithm to compute a one\=/switch \fullymeasured{} train track $(\tau',\mu')\carried(\tau,\mu)$, a co\=/orientation $\omega$ of $s$, and a choice $e''\in\{e,e'\}$ such that
\[
    \diamtt{\tau}{\tau'}\le\constant{deep nesting one step d},
\]
and $e''$ and all the ends of $\tau$ which are $\omega$\=/smaller than $e''$ are not persistent in $\tau'$.
 The constant $\constant{deep nesting one step d}$ is given by
 \[
    \declareconstant{deep nesting one step d}=\constantvalue{deep nesting one step d}.
 \]
 The running time of the algorithm is
 \[
    \bigO(\poly(\xi)\cdot(\norm{\mu}-\norm{\mu'}+1)\cdot\norm{\mu}\cdot\log\norm{\mu}).
 \]
\end{proposition}
\begin{proof}
Note that $\tau$ cannot be a simple closed curve. Let $b$ be the branch of which $e$ is an end, and let us write $b^\bullet$ for $e$.
\begin{substeps}
\item Let us first suppose that $b$ does not bound a punctured monogon in $S$.
In this case, the desired co\=/orientation $\omega$ of $s$ is the one such that
\[
    \height{\omega}(b^\bullet)>\height{\omega}(b_\bullet).
\]
Let
\[
    (\tau_0,\mu_0)\carries\ldots\carries(\tau_n,\mu_n)
\]
be the AHT sequence generated by the measured train track $(\tau,\mu)$, the switch $s$, and the co\=/orientation $\omega$.
We set $(\tau',\mu')=(\tau_i,\mu_i)$, where
\[
i=\min\{1\le j\le n:\text{$b$ is not untouched in $\tau_j$}\}.
\]
This is well\=/defined, because $b$ is not untouched in $\tau_n$ -- in fact, the train track $\tau_n$ is a simple closed curve.
By \cref{thm:end of non-untouched branch is non-persistent}, the end $e$ is not persistent in $\tau_i$.
We are now left to show that $\tau_i$ is close to $\tau_0$ in $\C(S)$.
To see this, we refer the reader to \cref{fig:deep nesting one step:no monogon}, which describes how to homotope $b$ to a smooth curve $c$ that intersects every fundamental curve of $\tau_0$ and of $\tau_{i-1}$ at most twice.
This gives the bound
\begin{align*}
\diamtt{\tau_0}{\tau_i}&\le\diam[\C(S)](\fundamental(\tau_0),c)+\diam[\C(S)](c,\fundamental(\tau_{i-1}))+{}\\
&\phantom{{}\le{}}\!\diamtt{\tau_{i-1}}{\tau_i}\\
&\le 2+2+\constant{proximity d}\le\constant{deep nesting one step d},
\end{align*}
where the bound on $\diamtt{\tau_{i-1}}{\tau_i}$ comes from \cref{thm:proximity}.

\begin{figure}
\centering
\begin{subcaptionblock}{.49\linewidth}
\centering
    \tikzsetnextfilename{deep-nesting-one-step-no-monogon-a}%
    \input{figures-source/deep-nesting-one-step-no-monogon-a.tex}%

\caption{}
\label{fig:deep nesting one step:no monogon:a}
\end{subcaptionblock}
\begin{subcaptionblock}{.49\linewidth}
\centering
    \tikzsetnextfilename{deep-nesting-one-step-no-monogon-b}%
    \input{figures-source/deep-nesting-one-step-no-monogon-b.tex}%

\caption{}
\label{fig:deep nesting one step:no monogon:b}
\end{subcaptionblock}
\caption{Since the branch $b$ of $\tau_0$ is untouched in $\tau_{i-1}$, there is a smooth curve $c$ homotopic to $b$ that intersects every fundamental curve of $\tau_0$ and of $\tau_{i-1}$ at most twice. When $b$ is returning, as in \subref{fig:deep nesting one step:no monogon:a}, we can take $c=b$. Otherwise, we take $c$ to be a copy of $b$ that is smoothed out near $s_0$, as shown in \subref{fig:deep nesting one step:no monogon:b}.}
\label{fig:deep nesting one step:no monogon}
\end{figure}

Clearly, the same argument works when $e'$ is an end of a branch of $\tau$ that does not bound a punctured monogon in $S$.
\item Let $\omega'$ be the co\=/orientation of $s$ such that $e\endl{\omega'}e'$, and let $\bar{e}$ be an end of $\tau$ such that $\direction(\bar{e})\neq\direction(e)$ and 
\begin{equation}\label{eqn:deep nesting one step:e bar definition}
    \height{\omega'}(\bar{e})\le\height{\omega'}(e)\le\height{\omega'}(\bar{e})+\mu(\bar{e}).
\end{equation}
We now assume that $\bar{e}$ is an end of a branch $\bar{b}\in\branches(\tau)$ which does not bound a punctured monogon in $S$.

By the previous case, we can find a co\=/orientation $\omega$ of $s$ and an integer $i$ such that, if $(\tau_i,\mu_i)$ is the $i$\=/th element of the AHT sequence generated by $(\tau,\mu)$, $s$, and $\omega$, then
\[
    \diamtt{\tau_0}{\tau_i}\le\constant{deep nesting one step d}
\]
and $\bar{e}$ is not persistent in $\tau_i$.
If $\omega'=\omega$, then \cref{itm:carrying properties:persistency condition,eqn:deep nesting one step:e bar definition} imply that
\[
\mu_i(s_i)\le\height{\omega}(\bar{e})\le\height{\omega}(e)
\]
and hence $e$ is not persistent in $\tau_i$ either.
If $\omega'=\bar{\omega}$, we use \cref{itm:carrying properties:persistency condition,eqn:deep nesting one step:e bar definition} again to deduce that
\begin{equation}
\begin{aligned}\label{eqn:deep nesting one step:flipped height}
    \mu_i(s_i)&\le\height{\omega}(\bar{e})\\
    &=\mu_0(s_0)-\height{\omega'}(\bar{e})-\mu_0(\bar{e})\\
    &\le\mu_0(s_0)-\height{\omega'}(e)\\
    &=\height{\omega}(e'),
\end{aligned}
\end{equation}
hence $e'$ is not persistent in $\tau_i$.

We can then set $(\tau',\mu')=(\tau_i,\mu_i)$, and $e''$ to be $e$ or $e'$ depending on whether $\omega'=\omega$ or not.
\item The only case we are left to consider is when $e$, $e'$, and $\bar{e}$ are all ends of branches of $\tau$ which bound punctured monogons in $S$.
Let us write $\bar{b}^\bullet$ for $\bar{e}$.
Up to swapping $e$ and $e'$, we can assume that $\bar{b}^\bullet\endg{\omega'}\bar{b}_\bullet$.
Set $\omega=\bar{\omega'}$, and let
\[
    (\tau_0,\mu_0)\carries\ldots\carries(\tau_n,\mu_n)
\]
be the AHT sequence generated by $(\tau,\mu)$, $s$, and $\omega$.
We set $(\tau',\mu')=(\tau_i,\mu_i)$, where
\[
i=\min\{1\le j\le n:\text{$b$ or $\bar{b}$ is not untouched in $\tau_j$}\}.
\]
The curve represented by the measure $\one{b}+\one{\bar{b}}$ is a fundamental curve of both $\tau_0$ and $\tau_{i-1}$, hence it intersects every fundamental curve of $\tau_0$ and of $\tau_{i-1}$ at most $4$ times.
The usual estimates then give the bound
\[
    \diamtt{\tau_0}{\tau_i}\le 5+5+\constant{proximity d}=\constant{deep nesting one step d}.
\]
Finally, we argue that one of $e$ and $e'$ is not persistent in $\tau_i$.
If $b$ is not untouched in $\tau_i$, then \cref{thm:end of non-untouched branch is non-persistent} implies that the $\omega$\=/smaller end of $b$ -- which is either $e$ or $e'$, since $b$ bounds a punctured monogon in $S$ -- is not persistent in $\tau_i$.
If, instead, the branch $\bar{b}$ is not untouched in $\tau_i$, then we use \cref{thm:end of non-untouched branch is non-persistent} again to deduce that $\bar{e}$ is not persistent in $\tau_i$.
In this case, the same computation carried out in \cref{eqn:deep nesting one step:flipped height} shows that $e'$ is not persistent in $\tau_i$.
\end{substeps}

Two final remarks are in order.
Firstly, we have shown how to find a \fullymeasured{} train track $(\tau',\mu')\carried(\tau,\mu)$, a co\=/orientation $\omega$ of $s$, and a choice $e''\in\{e,e'\}$ such that the end $e''$ is not persistent in $\tau'$.
Since $(\tau',\mu')$ comes from an AHT sequence, however, \cref{itm:carrying properties:persistency condition} implies that all the ends of $\tau$ which are $\omega$\=/smaller than $e''$ are also not persistent in $\tau'$.
Secondly, as far as the running time of the algorithm is concerned, note that we only need to perform the first $i$ steps of the procedure generating the AHT sequence, where $(\tau',\mu')=(\tau_i,\mu_i)$.
Each step can be performed in time $\bigO(\poly(\xi)\cdot\norm{\mu}\cdot\log\norm{\mu})$, hence \cref{eqn:aht sequence length} gives the desired bound for the running time.
\end{proof}

The final deep nesting algorithm is essentially a corollary of \cref{thm:deep nesting one step}.

\newconstant{deep nesting d}{D_{\mathrm{d.n}}}{Ddn}{12\xi\cdot\constant{deep nesting one step d}}
\begin{proposition}[Deep nesting algorithm]\label{thm:deep nesting}
Let $(\tau,\mu)$ be a one\=/switch \fullymeasured{} train track.
There is an algorithm to compute a one\=/switch \fullymeasured{} train track $(\tau',\mu')\carried(\tau,\mu)$ such that $\tau'$ is deeply nested in $\tau$ and
\[
    \diamtt{\tau}{\tau'}\le\constant{deep nesting d},
\]
where
\[
    \declareconstant{deep nesting d}=\constantvalue{deep nesting d}.
\]
The running time of the algorithm is
 \[
    \bigO(\poly(\xi)\cdot(\norm{\mu}-\norm{\mu'}+1)\cdot\norm{\mu}\cdot\log\norm{\mu})
 \]
\end{proposition}
\begin{proof}
The algorithm works as follows.
Start by setting $(\tau_0,\mu_0)=(\tau,\mu)$.
On the $i$\=/th step of the algorithm, for $i\ge 0$, we compute the set
\[
E_i=\{e\in\ends(\tau_0):\text{$e$ is persistent in $\tau_i$}\}.
\]
If $\card{E_i}=2$ then the algorithm terminates and returns $(\tau',\mu')=(\tau_i,\mu_i)$.
Otherwise, we can find two consecutive ends $e$ and $e'$ of $\tau_i$ such that $\carrymap[\tau_i][\tau_0](e)\neq\carrymap[\tau_i][\tau_0](e')$.
\Cref{thm:deep nesting one step} applied to the measured train track $(\tau_i,\mu_i)$ and the ends $e$ and $e'$ yields a one\=/switch \fullymeasured{} train track $(\tau_{i+1},\mu_{i+1})\carried(\tau_i,\mu_i)$ and a co\=/orientation $\omega$ of $\tau_i$ such that -- without loss of generality -- $e\endl{\omega}e'$, and $e$ and all the ends of $\tau_i$ which are $\omega$\=/smaller than $e$ are not persistent in $\tau_{i+1}$.
By \cref{itm:carrying properties:persistency condition}, it follows that $\carrymap[\tau_i][\tau_0](e)$ is not persistent in $\tau_0$.
In particular, the cardinality of $E_{i+1}$ is strictly smaller than that of $E_i$.
We then proceed to the $(i+1)$\=/th step of the algorithm.

Suppose that the algorithm terminates after $k$ steps -- or, in other words, that $\card{E_k}=2$.
Clearly, the integer $k$ is bounded above by $\card{\ends(\tau)}-2<12\xi$.
Since
\begin{align*}
    \diamtt{\tau_i}{\tau_{i+1}}\le\constant{deep nesting one step d}&&\text{for every $0\le i<k$},
\end{align*}
we have the desired estimate
\[
    \diamtt{\tau}{\tau'}\le k\cdot\constant{deep nesting one step d}<\constant{deep nesting d}.
\]

The bound on the running time follows from that of \cref{thm:deep nesting one step}.
\end{proof}

\subsection{Distances for carried curves}

Our next goal is to algorithmically estimate the distance in $\C(S)$ between a train track and a curve carried by it.
As an intermediate result, we show that if a curve is carried by a train track $\tau$, then it is also carried by a one\=/switch train track $\tau$ which is uniformly close to $\tau$ in $\C(S)$.

\newconstant{one switch d}{D_{\mathrm{1.s}}}{D1s}{48\xi^2\cdot\constant{proximity d}}
\begin{proposition}[One\=/switch algorithm]\label{thm:one switch algorithm}
Let $(\tau,\mu)$ be a \fullymeasured{} train track. There is an algorithm to compute a one\=/switch \fullymeasured{} train track $(\tau',\mu')\carried(\tau,\mu)$ such that
\[
    \diamtt{\tau}{\tau'}\le\constant{one switch d},
\]
where
\[
    \declareconstant{one switch d}=\constantvalue{one switch d}.
\]
The running time of the algorithm is
 \[
    \bigO(\poly(\xi)\cdot(\norm{\mu}-\norm{\mu'}+1)\cdot\norm{\mu}\cdot\log\norm{\mu}).
 \]
\end{proposition}
\begin{proof}
The algorithm works as follows.
Start by setting $(\tau_0,\mu_0)=(\tau,\mu)$.
On the $i$\=/step of the algorithm, for $i\ge 0$, choose a switch $s\in\switches(\tau_i)$ arbitrarily.
If $s$ is the unique switch of $\tau_i$, then the algorithm terminates with $(\tau',\mu')=(\tau_i,\mu_i)$.
Otherwise, we distinguish two cases.
\begin{substeps}
\item If $s$ has a returning branch $b$, then we arbitrarily choose another switch $s'$ of $\tau_i$ and a co\=/orientation $\omega$ of $s'$.
We let $(\tau_{i+1},\mu_{i+1})$ be the last measured train track in the AHT sequence generated by $(\tau_i,\mu_i)$, $s'$, and $\omega$.
We note that $\tau_{i+1}$ has one less switch than $\tau_i$, and moreover the branch $b$ is untouched in $\tau_{i+1}$.
This implies that $b$ is a curve whose intersection number with every fundamental curve of $\tau_i$ and of $\tau_{i+1}$ is at most $2$.
As a consequence, \cref{thm:distance bound from intersections} gives the bound
\begin{equation}\label{eqn:one switch algorithm:distance bound 1}
\diamtt{\tau_i}{\tau_{i+1}}\le\diam[\C(S)](\fundamental(\tau_i),b)+\diam[\C(S)](b,\fundamental(\tau_{i+1}))\le 4.
\end{equation}
\item If $s$ has no returning branches, then we choose a co\=/orientation $\omega$ of $s$, and we let $(\tau_{i+1},\mu_{i+1})$ be the last measured train track in the AHT sequence generated by $(\tau_i,\mu_i)$, $s$, and $\omega$.
Since $s$ has no returning branches, it is easy to see that this AHT sequence only involves \splitmove{} moves, and its length is bounded above by the number of ends of $\tau_i$ which are incident to $s$.
Therefore, by \cref{thm:proximity} we have the bound
\begin{equation}\label{eqn:one switch algorithm:distance bound 2}
    \diamtt{\tau_i}{\tau_{i+1}}\le\card{\ends(\tau_i)}\cdot \constant{proximity d}\le 12\xi\cdot\constant{proximity d}.
\end{equation}
\end{substeps}
We then proceed to the $(i+1)$\=/th step of the algorithm.

Suppose that the algorithm terminates after $k$ steps -- or, in other words, that $\tau_k$ is one\=/switch.
Clearly, the integer $k$ is bounded above by $\card{\switches(\tau)}-1<4\xi$.
From \cref{eqn:one switch algorithm:distance bound 1,eqn:one switch algorithm:distance bound 2} we get the desired inequality
\[
    \diamtt{\tau}{\tau'}\le k\cdot 12\xi\cdot\constant{proximity d}<\constant{one switch d}.
\]
The bound on the running time follows from \cref{eqn:aht sequence length}.
\end{proof}

The tools developed in this section, combined with \cref{thm:nested sequence distance}, find their purpose in the proof of the following theorem.

\begin{theorem}[Distance algorithm for carried curves]\label{thm:distance algorithm for carried curves}
Let $(\tau,\mu)$ be a \fullymeasured{} train track representing a curve $a\in\C_0(S)$. 
There is an algorithm to compute an integer $d\ge 0$ such that
\[
d-\constant{carried curves l+}\le\dist[\C(S)](\fundamental(\tau),a)\le\diam[\C(S)](\fundamental(\tau),a)\le\constant{carried curves l*}\cdot d+\constant{carried curves l+},
\]
where
\begin{align*}
\declareconstant{carried curves l+}&=\constantvalue{carried curves l+},\tag*{}\\
\declareconstant{carried curves l*}&=\constantvalue{carried curves l*},\tag*{}\\
\declareconstant{carried curves d}&=\constantvalue{carried curves d}.\tag*{}
\end{align*}
The running time of the algorithm is
\[
\bigO(\poly(\xi)\cdot\norm{\mu}^2\cdot\log\norm{\mu}).
\]
\end{theorem}
\begin{proof}
We split the algorithm in three steps.
\step{Multi\=/switch to one\=/switch.}
The first step is to apply \cref{thm:one switch algorithm} to obtain a one\=/switch \fullymeasured{} train track $(\tau_0,\mu_0)\carried(\tau,\mu)$ such that
\begin{equation}\label{eqn:distance algorithm for carried curves:tau to tau_0}
\diamtt{\tau}{\tau_0}\le\constant{one switch d}.
\end{equation}
If $\tau_0$ is not large, then $\diam[\C(S)](\fundamental(\tau),a)<\constant{one switch d}+2\le\constant{carried curves l+}$, hence the algorithm can return $d=0$.
Therefore, from now on, we assume that the train track $\tau_0$ is large.
\step{A sequence of nested train tracks.}
We now compute a sequence of measured train tracks
\begin{equation}\label{eqn:distance algorithm for carried curves:sequence}
    (\tau_0,\mu_0)\carries(\tau_0',\mu_0')\carries(\tau_1,\mu_1)\carries(\tau_1',\mu_1')\carries\ldots\carries(\tau_{n-1}',\mu_{n-1}')\carries(\tau_n,\mu_n)\carries(\bar{\tau},\bar{\mu}),
\end{equation}
according to the following procedure.

On the $i$\=/th step of the procedure, for $i\ge 0$, we apply the algorithm in \cref{thm:deep nesting} to $(\tau_i,\mu_i)$ to find a one\=/switch \fullymeasured{} train track $(\tau_i',\mu_i')\carried(\tau_i,\mu_i)$ such that $\tau_i'$ is deeply nested in $\tau_i$ and
\[
    \diamtt{\tau_i}{\tau_i'}\le\constant{deep nesting d}.
\]
If $\tau_i'$ is not large, we set $n=i$ and $(\bar{\tau},\bar{\mu})=(\tau_i',\mu_i')$, and we terminate the procedure.
Otherwise, we consider an AHT sequence generated by $(\tau_i',\mu_i')$, and let $(\tau_{i+1},\mu_{i+1})\carried(\tau_i',\mu_i')$ be the first measured train track in the sequence such that either $\tau_{i+1}$ is a simple closed curve, or $\norm{\mu_{i+1}}\le\norm{\mu_i'}-1$.
By \cref{thm:proximity} and \cref{thm:exponential decay}, we have the bound
\[
    \diamtt{\tau_i'}{\tau_{i+1}}\le\constant{proximity d}\cdot\constant{exponential decay n}.
\]
If $\tau_{i+1}$ is not large, we set $n=i$ and $(\bar{\tau},\bar{\mu})=(\tau_{i+1},\mu_{i+1})$, and we terminate the procedure.
Otherwise, we proceed to the $(i+1)$\=/th step of the procedure.

When the procedure terminates, we are left with a sequence \cref{eqn:distance algorithm for carried curves:sequence} of nested measured train tracks, of which we remark a few properties.
\begin{itemize}
\item Since $\norm{\mu_{i+1}}\le\norm{\mu_i}-1$ for every $0\le i<n$, the length $n$ is bounded above by $\norm{\mu}$; in particular, this shows that we can compute the sequence \cref{eqn:distance algorithm for carried curves:sequence} in time
\[
    \bigO(\poly(\xi)\cdot\norm{\mu}^2\cdot\log\norm{\mu}).
\]
\item For $0\le i<n$, we have the uniform bound
\begin{align*}
    \diamtt{\tau_i}{\tau_{i+1}}&\le\diamtt{\tau_i}{\tau_i'}+\diamtt{\tau_i'}{\tau_{i+1}}\\
    &\le \constant{deep nesting d}+\constant{proximity d}\cdot\constant{exponential decay n}=\constant{carried curves d}.
\end{align*}
\item Similarly, we have that $\diamtt{\tau_n}{\bar{\tau}}\le\constant{carried curves d}$.
\item The train tracks $\tau_i$ and $\tau_i'$ are large for $0\le i<n$, and so is $\tau_n$.
\item The train track $\bar{\tau}$ is not large, hence $\diam[\C(S)](\fundamental(\bar{\tau}),a)\le 2$. As a consequence, we see that
\begin{equation}\label{eqn:distance algorithm for carried curves:tau_n to a}
    \diam[\C(S)](\fundamental(\tau_n),a)\le \constant{carried curves d}+2.
\end{equation}
\end{itemize}
\step{A quasi\=/geodesic subsequence.}
Finally, we show how to compute a subsequence of \cref{eqn:distance algorithm for carried curves:sequence} whose length is a coarse estimate for the distance in $\C(S)$ between $\tau$ and $a$.
Start by setting $i_0=0$.
On the $j$\=/th step of the procedure, for $j\ge 0$, we consider the set
\[
I_j=\left\{
i_j<k\le n:\,
\begin{matrix*}[l]
\text{for every fundamental measure $\nu$ on $\tau_k$, the}\\
\text{measure $\carrymap[\tau_k][\tau_{i_j}](\nu)$ fills a large subtrack of $\tau_{i_j}$}
\end{matrix*}
\right\}.
\]
If $i_j=n$, then we set $m=j$ and terminate the procedure.
If $i_j<n$ but $I_j$ is empty, we set $i_{j+1}=n$ and $m=j+1$, and terminate the procedure.
Otherwise, we set $i_{j+1}=\min I_j$, and proceed to the $(j+1)$\=/th step.

In terms of computational complexity, note that, for a fixed integer $0\le k\le n$, there are only $\bigO(\xi^2)$ fundamental measures on $\tau_k$, and for each $\nu\in\fundamental(\tau_k)$ we can compute the sets
\[
\{b\in\branches(\tau_h):\carrymap[\tau_k][\tau_h](\nu)\ge\one{b}\}\qquad\text{for $0\le h\le k$}
\]
in time $\bigO(\poly(\xi)\cdot k)=\bigO(\poly(\xi)\cdot\norm{\mu})$.
This information is enough to construct the subsequence $\{i_j:0\le j\le m\}$ which, given the sequence \cref{eqn:distance algorithm for carried curves:sequence}, can be done in time
\[
\bigO(\poly(\xi)\cdot\norm{\mu}^2).
\]

Our goal is now to prove that $d=\lfloor m/2\rfloor$ is a valid return value.
\begin{substeps}
\item For the upper bound, fix an integer $1\le j\le m$, and note that there is a fundamental curve $\nu$ of $\tau_{i_j-1}$ which is carried by a non\=/large subtrack $\sigma$ of $\tau_{i_{j-1}}$.
As a consequence, there is an essential curve $c\in\C_0(S)$ contained in the complement of $\sigma$, and we can choose $c$ so that it intersects every branch of $\tau_{i_{j-1}}$ at most twice. In particular, the curve $c$ intersects every fundamental curve of $\tau_{i_{j-1}}$ at most $4$ times, and 
\Cref{thm:distance bound from intersections} then gives the bound
\[
\diam[\C(S)](\fundamental(\tau_{i_{j-1}}),\nu)\le\diam[\C(S)](\fundamental(\tau_{i_{j-1}}),c)+\dist[\C(S)](c,\nu)\le 5+1=6,
\]
which in turn yields
\begin{align*}
\diamtt{\tau_{i_{j-1}}}{\tau_{i_j}}&\le\diam[\C(S)](\fundamental(\tau_{i_{j-1}}),\nu)+\diam[\C(S)](\nu,\fundamental(\tau_{i_j}))\\
&\le 6+\constant{carried curves d}.
\end{align*}
We conclude that
\begin{equation}\label{eqn:distance algorithm for carried curves:tau_0 to tau_n upper}
    \diamtt{\tau_0}{\tau_n}\le(\constant{carried curves d}+6)\cdot m
\end{equation}
\item For the lower bound, we start by proving that $\PN(\tau_{i_j}')\subs\intPN(\tau_{i_{j-1}})$ for every $1\le j<m$.
Let $\sigma'$ be a large recurrent subtrack of $\tau_{i_j}'$, and let $\rho'$ be a recurrent diagonal extension of $\sigma'$.
The carrying $\tau_{i_j}'\carried\tau_{i_j}$ induces a carrying $\sigma'\carried\sigma$ for a large recurrent subtrack $\sigma$ of $\tau_{i_j}$; by choosing the smallest possible such $\sigma$, we can assume that $\sigma'$ fills $\sigma$.
Since $\tau_{i_j}'$ is deeply nested in $\tau_{i_j}$, it is clear that $\sigma'$ is deeply nested in $\sigma$, the two persistent ends of $\sigma$ in $\sigma'$ being exactly the two persistent ends of $\tau_{i_j}$ in $\tau_{i_j}'$.
\Cref{thm:deep nesting lemma} then implies that $\rho'$ is carried by and deeply nested in a recurrent diagonal extension $\rho$ of $\sigma$, in such a way that the two persistent ends of $\rho$ in $\rho'$ are ends of branches of $\sigma$.
Moreover, \cref{thm:deep nesting lemma} applied to the carrying $\rho'\carried\rho$ shows that there is a fundamental measure $\nu$ on $\rho$ such that $\carrymap[\rho'][\rho](\mu')\ge\nu$ for every non\=/zero $\mu'\in\P(\rho')$, and in fact the measure $\nu$ is supported on $\sigma$.
We refer the reader to \cref{fig:distance algorithm for carried curves:deep nesting} for a graphical depiction of the carryings $\sigma'\carried\sigma$ and $\rho'\carried\rho$.

\begin{figure}
\centering
\begin{subcaptionblock}{.49\linewidth}
\centering
    \tikzsetnextfilename{distance-algorithm-for-carried-curves-deep-nesting-a}%
    \input{figures-source/distance-algorithm-for-carried-curves-deep-nesting-a.tex}%

\caption{}
\label{fig:distance algorithm for carried curves:deep nesting:a}
\end{subcaptionblock}
\begin{subcaptionblock}{.49\linewidth}
\centering
    \tikzsetnextfilename{distance-algorithm-for-carried-curves-deep-nesting-b}%
    \input{figures-source/distance-algorithm-for-carried-curves-deep-nesting-b.tex}%

\caption{}
\label{fig:distance algorithm for carried curves:deep nesting:b}
\end{subcaptionblock}
\caption{\subref{fig:distance algorithm for carried curves:deep nesting:a} The subtrack $\sigma'$ is carried by and deeply nested in $\sigma$; the two ends of $\sigma$ which are persistent in $\sigma'$ are highlighted. \subref{fig:distance algorithm for carried curves:deep nesting:b} The diagonal extension $\rho'$ of $\sigma'$ is carried by and deeply nested in a diagonal extension $\rho$ of $\sigma$. The two ends of $\rho$ which are persistent in $\rho'$ are ends of branches of $\sigma$.}
\label{fig:distance algorithm for carried curves:deep nesting}
\end{figure}

The train track $\sigma$ fills a large subtrack $\sigma_*$ of $\tau_{i_{j-1}}$, and $\rho$ is carried by a large recurrent diagonal extension $\rho_*$ of $\sigma_*$.
If $\mu'$ is any non\=/zero measure on $\rho'$, the inequality $\carrymap[\rho'][\rho](\mu')\ge\nu$ implies that $\carrymap[\rho'][\rho_*](\mu')\ge\carrymap[\rho][\rho_*](\nu)$.
Since $\carrymap[\tau_{i_j}][\tau_{i_{j-1}}](\nu)$ fills a large subtrack of $\tau_{i_{j-1}}$, we see that $\carrymap[\rho'][\rho_*](\mu')$ fills a diagonal extension of a large subtrack of $\tau_{i_{j-1}}$.
In other words, we have proved that $\PN(\tau_{i_j}')\subs\intPN(\tau_{i_{j-1}})$, and in particular that $\PN(\tau_{i_{j+1}})\subs\intPN(\tau_{i_{j-1}})$.
By means of \cref{thm:nested sequence distance,eqn:one-switch train track fundamental diameter}, we finally get the lower bound
\begin{equation}\label{eqn:distance algorithm for carried curves:tau_0 to tau_n lower}
\disttt{\tau_0}{\tau_n}\ge\left\lfloor\frac{m}{2}\right\rfloor-\diam[\C(S)](\fundamental(\tau_0))\ge\left\lfloor\frac{m}{2}\right\rfloor-5.
\end{equation}
\end{substeps}

We now combine \cref{eqn:distance algorithm for carried curves:tau to tau_0,eqn:distance algorithm for carried curves:tau_0 to tau_n upper,eqn:distance algorithm for carried curves:tau_0 to tau_n lower,eqn:distance algorithm for carried curves:tau_n to a} to get
\begin{align*}
\diam[\C(S)](\fundamental(\tau),a)&\le(\constant{carried curves d}+6)\cdot m+\constant{one switch d}+\constant{carried curves d}+2,\\
\dist[\C(S)](\fundamental(\tau),a)&\ge\left\lfloor\frac{m}{2}\right\rfloor-\constant{one switch d}-\constant{carried curves d}-7.
\end{align*}
Therefore, $d=\lfloor m/2\rfloor$ is a valid return value.
\end{proof}

\subsection{Distances in the curve graph}

The main application of \cref{thm:distance algorithm for carried curves} is to compute coarse distances between curves in the curve graph.
In particular, we will show in \cref{thm:distance algorithm for short curves} how to estimate the distance between a very ``short'' curve and any other curve.
A theorem of \textcite{lackenby-pants-graph} will then allow us to generalise this result to coarsely compute distances between any two curves.
We start by describing how to find a train track carrying an arbitrary normal curve on a surface $S$ with a triangulation $\TTT$.
This task is essentially trivial when $S$ is punctured, but requires some effort in the closed case.

\begin{proposition}\label{thm:find carrying train track}
Let $a\in\C_0(S)$ be a curve on $S$.
There is an algorithm to find a train track $\tau$ and a measure $\mu\in\intP(\tau)$ representing $a$ such that every edge of $\TTT$ intersects $\tau$ transversely and at most once.
The running time of the algorithm is polynomial in $\xi$ and $\norm{a}$.
\end{proposition}
\begin{proof}
In the context of this proof, by \emph{pre\=/track} we mean a smooth 1\=/complex embedded in $S$ which satisfies all the properties of a train track except for those concerning complementary regions.

When $S$ is punctured, the argument is standard.
We let $\tau_0$ be the pre\=/track obtained by repeating the pattern shown in \cref{fig:find carrying train track:pattern} for each triangle of $\TTT$.
The pre\=/track $\tau_0$ has $\card{\TTT}$ switches, all of which are $4$\=/valent, and three branches per triangle of $\TTT$.
We can define a measure $\mu_0$ on $\tau_0$ by counting the number of normal arcs that run parallel to each branch, as shown in \cref{fig:find carrying train track:measure}.
By removing branches of measure zero from $\tau_0$, we get a train track $\tau$ and an induced measure $\mu\in\intP(\tau)$ representing $a$.

\begin{figure}
\centering
\newsavebox{\largestimage}
\savebox{\largestimage}{%
    \tikzsetnextfilename{find-carrying-train-track-pattern}%
    \input{figures-source/find-carrying-train-track-pattern.tex}%
}
\begin{subcaptionblock}{.35\linewidth}
\centering
\usebox{\largestimage}
\caption{}
\label{fig:find carrying train track:pattern}
\end{subcaptionblock}
\begin{subcaptionblock}{.55\linewidth}
\centering
\raisebox{\dimexpr.5\ht\largestimage-.5\height}{%
    \tikzsetnextfilename{find-carrying-train-track-measure}%
    \input{figures-source/find-carrying-train-track-measure.tex}%
}
\caption{}
\label{fig:find carrying train track:measure}
\end{subcaptionblock}
\caption{\subref{fig:find carrying train track:pattern} The pre\=/track $\tau_0$ has three branches in each triangle of $\TTT$, arranged to form a cusped triangle.
\subref{fig:find carrying train track:measure} The measure $\mu_0$ assigns to each branch of $\tau_0$ weight equal to the number of normal arcs of $a$ which run parallel to it.}
\end{figure}

In general, let us call the pre\=/track $\tau$ obtained from the procedure above the \emph{pre\=/track associated} to the normal curve $a$.
In the closed case, there is no guarantee that the pre\=/track associated to $a$ is an actual train track, because it may have monogon or bigon complementary regions.
Therefore, we need to preemptively modify the normal curve $a$, by possibly isotopying it across the vertex $v$ of $\TTT$.

Let us number the types of normal arcs in $\TTT$ by $1,2,\ldots,n$  in the order in which they appear counter\=/clockwise around $v$; here $n=2\card{\TTT}$ is three times the number of triangles of $\TTT$.
We will call a normal arc of $a$ \emph{innermost} if it is the closest to $v$ amongst all normal arcs of the same type, and \emph{outermost} if it is the farthest away.
For a normal curve $b$ and an integer $1\le i\le n$, let us call $x(b)_i$ the number of normal arcs of $b$ of type $i$.
Clearly, the integers $x(b)_1,\ldots,x(b)_n$ are enough to recover the vector $\TTT(b)$, and vice versa.

We represent the triangles of $\TTT$ as they appear around $v$, as shown in \cref{fig:find carrying train track:monogon and bigon}.
For $1\le i\le n$, let $\Delta_i$ be the triangle of $\TTT$ containing the normal arc type $i$.
We remark that, in this representation, every triangle appears three times, but a triangle can never be adjacent to itself; in other words, $\Delta_i\neq\Delta_{i+1}$ for $1\le i<n$, and $\Delta_n\neq\Delta_1$.
This follows from the fact that, in a one\=/vertex triangulation of a surface, a triangle cannot be glued to itself along two of its edges.

\begin{figure}
\centering
\begin{subcaptionblock}{.45\linewidth}
\centering
    \tikzsetnextfilename{find-carrying-train-track-monogon-and-bigon-a}%
    \input{figures-source/find-carrying-train-track-monogon-and-bigon-a.tex}%

\caption{}
\label{fig:find carrying train track:monogon and bigon:a}
\end{subcaptionblock}
\begin{subcaptionblock}{.45\linewidth}
\centering
    \tikzsetnextfilename{find-carrying-train-track-monogon-and-bigon-b}%
    \input{figures-source/find-carrying-train-track-monogon-and-bigon-b.tex}%

\caption{}
\label{fig:find carrying train track:monogon and bigon:b}
\end{subcaptionblock}
\caption{We draw the triangles of $\TTT$ in the order they appear counter\=/clockwise around the vertex $v$ of $\TTT$. These pictures show what the pre\=/track associated to $a$ looks like when it has a monogon \subref{fig:find carrying train track:monogon and bigon:a} or bigon \subref{fig:find carrying train track:monogon and bigon:b} complementary region.}
\label{fig:find carrying train track:monogon and bigon}
\end{figure}

We claim that the pre\=/track $\tau_0$ associated to a normal curve $b$ has a monogon (respectively bigon) complementary region if and only if $x(b)_i=0$ for exactly one (respectively two) integers $1\le i\le n$.
The reverse implication is clear, as shown in \cref{fig:find carrying train track:monogon and bigon:a} (respectively \cref{fig:find carrying train track:monogon and bigon:b}).
For the forward implication, suppose that a component $Z$ of $S\setminus\tau_0$ is a monogon (respectively bigon).
Since no complementary region entirely contained in a triangle of $\TTT$ is a monogon (respectively bigon), the vertex $v$ must lie in $Z$; since no edge can be entirely contained in a monogon (respectively bigon), every edge of $\TTT$ intersects $\tau_0$ exactly once.
It follows that, if $x(b)_i=0$ for some integer $1\le i\le n$, then $x(b)_j$ and $x(b)_k$ are both positive, where $j$ and $k$ refer to the two types of normal arcs other than $i$ which appear in $\Delta_i$.
Therefore, the region $Z$ is a cusped disc, with as many cusps as there are integers $1\le i\le n$ such that $x(b)_i=0$.

We now describe how to modify the normal curve $a$ to ensure that $\tau_0$ is an actual train track.
\begin{substeps}
\item In order to get rid of monogons, we begin by applying \cite[Theorem 6.3]{lackenby-minimal-position} to put $a$ in minimal position with an arbitrary edge $e$ of $\TTT$; for ease of notation, we keep calling the resulting normal curve $a$.
Since $a$ is essential, we can assume that $x(a)_1=0$.
Suppose, for the sake of contradiction, that the pre\=/track associated to $a$ has a monogon complementary region; in other words, assume that $x(a)_i\ge 1$ for every $2\le i\le n$.
This situation is depicted in \cref{fig:find carrying train track:monogon:a}.

Let $2\le h,k\le n$ be the integers such that the normal arc types $1$, $h$, and $k$ appear counter\=/clockwise in $\Delta_1$.
We observe that, since $\Delta_2\neq\Delta_1$, we have that $x(a)_k=x(a)_{k-1}+x(a)_2\ge 2$.
Similarly, the integer $x(a)_h$ is also at least $2$.
Let $\alpha$ be the subarc of $a$ obtained by taking the union of the innermost normal arcs of types $2,\ldots,n$ and the outermost normal arcs of types $k$ and $h$.
Since $x(a)_h\ge 2$ and $x(a)_k\ge 2$, we see that these normal arcs are all different, and hence the arc $\alpha$ is embedded in $S$.
We can then isotope $\alpha$ across $v$, fixing $a\setminus\alpha$, to reduce the number of intersections with $e$, as depicted in \cref{fig:find carrying train track:monogon:b}.
It follows that $a$ was not in minimal position with $e$.
\begin{figure}
\centering
\begin{subcaptionblock}{.45\linewidth}
\centering
    \tikzsetnextfilename{find-carrying-train-track-monogon-a}%
    \input{figures-source/find-carrying-train-track-monogon-a.tex}%

\caption{}
\label{fig:find carrying train track:monogon:a}
\end{subcaptionblock}
\begin{subcaptionblock}{.45\linewidth}
\centering
    \tikzsetnextfilename{find-carrying-train-track-monogon-b}%
    \input{figures-source/find-carrying-train-track-monogon-b.tex}%

\caption{}
\label{fig:find carrying train track:monogon:b}
\end{subcaptionblock}
\caption{If the pre\=/track associated to $a$ has a monogon complementary region \subref{fig:find carrying train track:monogon:a}, then $a$ is not in minimal position with $e$, since it can be isotoped to reduce the number of intersections with it \subref{fig:find carrying train track:monogon:b}.}
\end{figure}
\item Suppose now that the pre\=/track associated to $a$ has a bigon complementary region.
This situation is depicted in \cref{fig:find carrying train track:bigon:a}.
There is a unique integer $2\le j\le n$ such that $x(a)_j=0$, and we can assume that $2\le j\le n/2+1$.
Let $1\le l,m\le n$ be the integers such that the normal arc types $j$, $l$, and $m$ appear counter\=/clockwise in $\Delta_j$.

We remark that, if $j<h$, then we have that $\{h+1,n\}\cap\{1,j\}=\emptyset$, and hence $x(a)_h=x(a)_{h+1}+x(a)_n\ge 2$.
Similarly, if $j<m$, then $x(a)_m=x(a)_{m-1}+x(a)_{j+1}\ge 2$.

Let $\alpha$ be the subarc of $a$ obtained by taking the union of the innermost normal arcs of types $j+1,\ldots,n$ and the outermost normal arcs of types $h$ and $m$.
By the remark above, all these normal arcs are different, and the arc $\alpha$ is embedded in $S$.
We can then isotope $\alpha$ across $v$, fixing $a\setminus\alpha$, to obtain a new normal curve $a'$, as shown in \cref{fig:find carrying train track:bigon:b}.
The numbers $x(a')_1,\ldots,x(a')_n$ can easily be computed in terms of $x(a)_1,\ldots,x(a)_n$ as follows.
For each $1\le i\le n$, define
\[
\NewDocumentCommand{\cardicap}{m}{\card{\{i\}\cap\{#1\}}}
    r(i)=\cardicap{2,\ldots,j-1}+\cardicap{k,l}-\cardicap{j+1,\ldots,n}-\cardicap{h,m}.
\]
Then $x(a')_i=x(a)_i+r(i)$.
In particular, note that $x(a')_1=x(a')_j=0$.
This means that either $x(a')_i=0$ for at least three integers $1\le i\le n$, or we can repeat this arc isotopy procedure with $a'$ instead of $a$.
By induction, if we let
\[
    t=\min\left(\{x(a)_i:1\le i\le n,r(i)=-1\}\cup\left\{\left\lfloor\frac{x(a)_i}{2}\right\rfloor:1\le i\le n,r(i)=-2\right\}\right),
\]
then there is a normal curve $a''$ with $x(a'')_i=x(a)_i+t\cdot r(i)$, obtained by iterating the arc isotopy procedure $t$ times.
This curve satisfies $x(a'')_1=x(a'')_j=0$, and $x(a'')_i=0$ for a third integer $2\le i\le n$.
In particular, the normal curve $a''$ is isotopic to $a$, and its associated pre\=/track has no monogon or bigon complementary regions.
\end{substeps}

\begin{figure}
\centering
\begin{subcaptionblock}{.45\linewidth}
\centering
    \tikzsetnextfilename{find-carrying-train-track-bigon-a}%
    \input{figures-source/find-carrying-train-track-bigon-a.tex}%

\caption{}
\label{fig:find carrying train track:bigon:a}
\end{subcaptionblock}
\begin{subcaptionblock}{.45\linewidth}
\centering
    \tikzsetnextfilename{find-carrying-train-track-bigon-b}%
    \input{figures-source/find-carrying-train-track-bigon-b.tex}%

\caption{}
\label{fig:find carrying train track:bigon:b}
\end{subcaptionblock}
\caption{If the pre\=/track associated to $a$ has a bigon complementary region \subref{fig:find carrying train track:bigon:a}, we can isotope a subarc $\alpha$ of $a$ across the vertex of the triangulation \subref{fig:find carrying train track:monogon and bigon:b}.}
\end{figure}

Two concluding remarks in the closed case are in order.
Firstly, note that the curve $a''$ -- or, more precisely, the integers $x(a'')_1,\ldots,x(a'')_n$ -- can be computed in one single step, without actually performing the arc isotopy procedure $t=\bigO(2^{\norm{a}})$ times.
This guarantees that the running time of the algorithms remains polynomial in $\norm{a}$.
Secondly, for the bigon removal argument to work, it is crucial that the genus of $S$ is at least $2$.
In fact, this condition guarantees that $n>6$ and hence that $r(i)\in\{-1,-2\}$ for some $1\le i\le n$, thus making $t$ well\=/defined.
When $S$ is a torus, it may happen that $a'=a$ as normal curves, thus rendering the arc isotopy procedure ineffective.
\end{proof}

As anticipated above, \cref{thm:find carrying train track} immediately gives an algorithm to estimate the distance in $\C(S)$ between two curves $a$ and $b$, one of which ($a$ in the statement below) is short.

\begin{proposition}[Distance algorithm for short curves]\label{thm:distance algorithm for short curves}
Let $a,b\in\C_0(S)$ be curves on $S$.
There is an algorithm to compute an integer $d\ge 0$ such that
\[
d-\frac{6\norm{a}_1+6}{\log\xi}-\constant{carried curves l+}-2\le\dist[\C(S)](a,b)\le\constant{carried curves l*}\cdot d+\frac{6\norm{a}_1+6}{\log\xi}+\constant{carried curves l+}+2.
\]
The running time of the algorithm is polynomial in $\xi$ and $\norm{b}$.
\end{proposition}
\begin{proof}
We use \cref{thm:find carrying train track} to compute a \fullymeasured{} train track $(\tau,\mu)$ representing $b$, such that $\tau$ intersect every edge of $\TTT$ at most once.
Note that, as a consequence, every fundamental curve of $\tau$ can be realised so that it intersects $a$ at most
\[
\sum_{e\in\TTT}2\TTT(a)_e
\]
times.
\Cref{thm:distance bound from intersections} immediately implies that
\[
\diam[\C(S)](a,\fundamental(\tau))\le\frac{6\norm{a}_1+6}{\log\xi}+2.
\]
It is then clear that the integer $d$ returned by the algorithm in \cref{thm:distance algorithm for carried curves} applied to $(\tau,\mu)$ satisfies the desired inequalities.
\end{proof}

\begin{remark}\label{rmk:distance algorithm polynomial degree}
The statement of \cref{thm:distance algorithm for short curves} claims that the running time of the algorithm is polynomial in $\norm{b}$, without any estimate on the degree of the polynomial.
The reason for this omission is that the proof relies on the algorithm of \cref{thm:find carrying train track}, which in turn uses \cite[Theorem 6.3]{lackenby-minimal-position} in the case where $S$ is closed; \citeauthor{lackenby-minimal-position} does not give any upper bound to the degree of the polynomial in the running time of the latter.

However, we can be more precise if we assume that $S$ is punctured.
In this case, the algorithm of \cref{thm:find carrying train track} runs in linear time in $\norm{b}$.
As a consequence, the running time of the algorithm in \cref{thm:distance algorithm for short curves} is the same as that of \cref{thm:distance algorithm for carried curves}, that is
\[
\bigO(\poly(\xi)\cdot\norm{b}^2\cdot\log\norm{b}).\qedhere
\]
\end{remark}

In order to get rid of the dependence on $\norm{a}_1$ of the estimate in \cref{thm:distance algorithm for short curves}, we use an algorithm of \citeauthor{lackenby-pants-graph} to find a triangulation of $S$ with respect to which $a$ is very short.
We remark that a similar result was already proved by \textcite{bell-simplifying-triangulations}, but without the polynomial dependence on the complexity of $S$.

\begin{theorem}[{\cite[Theorems 7.4 and 7.5]{lackenby-pants-graph}}]
\label{thm:simplifying triangulation}
Let $a\in\C_0(S)$ be a curve.
There is an algorithm to produce a path in the flip and twist graph of $S$ between $\TTT$ and a triangulation $\TTT'$ of $S$ such that $a$ intersects every edge of $\TTT'$ at most twice.
The length of the path and the running time of the algorithm are polynomial in $\xi$ and $\norm[\TTT]{a}$.
\end{theorem}

Note that Theorems 7.4 and 7.5 of \cite{lackenby-pants-graph} are not stated in terms of flip and twist graph.
However, the sequences of triangulations of $S$ that \citeauthor{lackenby-pants-graph} construct can be interpreted as paths in $\GGG(S)$, where two consecutive triangulations are joined by an arc of length $\bigO(\log(\xi)+\norm[\TTT,1]{a})$.
Since the lengths of their sequences and the running times of their algorithms are polynomial in $\xi$ and $\norm[\TTT]{a}$, this justifies the statement of \cref{thm:simplifying triangulation}.
We also remark that, for the closed case, Theorem 7.4 of \cite{lackenby-pants-graph} produces a triangulation $\TTT'$ of which $a$ is an edge.
However, every edge of $\TTT'$ can be isotoped to a normal curve intersecting every (other) edge at most twice, as we require in \cref{thm:simplifying triangulation}.

By combining \cref{thm:distance algorithm for short curves,thm:simplifying triangulation}, we finally have a general algorithm to estimate distances in the curve graph.

\begin{corollary}[Distance algorithm]\label{thm:distance algorithm}
Let $a,b\in\C_0(S)$ be curves.
There is an algorithm to compute an integer $d\ge 0$ such that
\[
d-\constant{curves l+}\le\dist[\C(S)](a,b)\le\constant{carried curves l*}\cdot d+\constant{curves l+},
\]
where
\[
\declareconstant{curves l+}=\constantvalue{curves l+}.
\]
The running time of the algorithm is polynomial in $\xi$, $\norm{a}$, and $\norm{b}$.
\end{corollary}
\begin{proof}
We use \cref{thm:simplifying triangulation} to find a path in $\GGG(S)$ between $\TTT$ and a triangulation $\TTT'$, every edge of which intersects $a$ at most twice.
The path has polynomial length in $\xi$ and $\norm[\TTT]{a}$, and we use it to compute the vectors $\TTT'(a)$ and $\TTT'(b)$.
We observe that the sum of the entries of $\TTT'(a)$ is at most $2\card{\TTT'}\le-6\chi+6$; as a consequence, some easy algebra shows that
\[
    \frac{6\norm[\TTT',1]{a}+6}{\log\xi}\le 34.
\]
Moreover, note that the complexity $\norm[\TTT']{b}$ is polynomial in $\norm[\TTT]{b}$ and in the length of the path in $\GGG(S)$ from $\TTT$ to $\TTT'$; in particular, we have that
\[
    \norm[\TTT']{b}=\bigO(\poly(\xi)\cdot\poly(\norm[\TTT]{a})\cdot\poly(\norm[\TTT]{b})).
\]
Finally, we apply \cref{thm:distance algorithm for short curves} to the curves $a$ and $b$ with respect to the triangulation $\TTT'$; the estimate $d$ that we get satisfies the required bounds for this algorithm.
\end{proof}

\begin{remark}
For the sake of simplicity, we have stated \cref{thm:distance algorithm for carried curves} as a coarse distance algorithm.
However, it is clear from the proof that the same algorithm can be used to produce an explicit quasi\=/geodesic between a fundamental curve $\nu$ of $\tau$ and $a$.
In fact, using the notation from the proof, we can arbitrarily choose a fundamental curve $\nu_j$ of $\tau_{i_j}$ for each $0\le j\le m$.
It is then easy to see that $\nu,\nu_0,\ldots,\nu_m,a$ is a $(\constant{carried curves l*},\constant{carried curves l+})$\=/quasi\=/geodesic from $\nu$ to $a$, in the sense that
\[
    |j-k|-\constant{carried curves l+}\le\dist[\C(S)](\nu_j,\nu_k)\le\constant{carried curves l*}\cdot|j-k|+\constant{carried curves l+}\qquad\text{for $-1\le j,k\le m+1$,}
\]
where $\nu_{-1}=\nu$ and $\nu_{m+1}=a$.

Note that, since all the carrying maps in the sequence \cref{eqn:distance algorithm for carried curves:sequence} are computed explicitly, we can think of the curves $\nu_{-1},\ldots,\nu_{m+1}$ as being carried by $\tau$.
It is then straightforward to go through the proofs of \cref{thm:distance algorithm for short curves,thm:distance algorithm} and obtain an algorithm to compute a $(\constant{carried curves l*},\constant{curves l+})$\=/quasi\=/geodesic between two given curves $a,b\in\C_0(S)$.
\end{remark}

%% file: figures-source/split-move.tex
\bgroup
\def\tikzttsetup#1{
\tikzset{
tt/new=#1,
tt/#1/.cd,
set sep=.5cm,
set width=.5cm,
add switch={0,0},
add switch={1.8,.8},
set switch thickness={2}{.75cm},
add branch standard={1}{2},
add bunch={1}{left}{5,spread=1.2cm,sep=.2cm},
add bunch={1}{below 1}{2,spread=1.2cm,sep=.4cm},
add bunch={2}{above 1}{1,spread=1.2cm,sep=.8cm},
add bunch={2}{right}{3,spread=.8cm,sep=.4cm},
set branch weights={2=1,3=1,4=1,5=2,6=3,1=4,7=2,8=2,9=2,10=3,11=1,12=2},
finalise,
}
\fill[grey background] [tt/#1/use branch background all but={2,3,4,1},tt/#1/use switch background all];
\fill[\colourmakebackground{3}] [tt/#1/use branch background={2,3,4}];
\fill[\colourmakebackground{2}] [tt/#1/use branch background={1}];
\draw[black,line] [tt/#1/use branch contour all but={2,3,4,1},tt/#1/use switch contour all];
\draw[\colourmakedark{3},line] [tt/#1/use branch contour={2,3,4}];
\draw[\colourmakedark{2},line] [tt/#1/use branch contour={1}];
}
\begin{tikzpicture}
\tikzttsetup{a}
\node[\colourmakedark{3},left,xshift=-2pt,yshift=2pt] at (spath cs:\ttbranch{a}{2} 1) {$b_1$};
\node[\colourmakedark{3},left,xshift=-2pt] at (spath cs:\ttbranch{a}{4} 1) {$b_h$};
\node[\colourmakedark{2},below right,yshift=-.35cm,xshift=2pt,label node] at (spath cs:\ttbranch{a}{1} 1) {$\bar{b}$};
\draw[pin line] (label) -- (spath cs:\ttbranch{a}{1} .8);
\node[black,above=.5cm] at (a-1) {$s$};
\node[black,left,xshift=-2pt] at (spath cs:\ttbranch{a}{6} 1) {$\tau$};
\coordinate (p) at ($(a-1)-(0,.75)$);
\draw[palette 4,line,-{Triangle[]}] ++(p) -- +(-.6,0) node[\colourmakedark{4},below] {$\eta$};
\draw[palette 4,line,-{Triangle[]}] ++(p) -- +(0,-.6) node[\colourmakedark{4},above right] {$\omega$};
\fill[black] (p) circle (.3pt);

\begin{scope}[shift={(7,0)}]
\tikzttsetup{b}
\tikzset{
tt/b/carried={
name=c,
switch={1}{-.1cm}{4,5}{4,5,6},
switch={2}{0}{1,...,5}{1,2,3},
1=4,2=1,3=1,4=1,5=1,6=1,7=1,8=1,9=1,10=1,11=1,12=1,
},
}
\draw[palette 3,line] [tt/c/use branch={2,3,4}];
\draw[palette 2,line] [tt/c/use branch={1}];
\draw[palette 1,line] [tt/c/use branch all but={1,2,3,4}];
\pic[palette 1] at (c-1) {dot};
\pic[palette 1] at (c-2) {dot};
\node[\colourmakedark{3},left,xshift=-2pt,yshift=2pt] at (spath cs:\ttbranch{b}{2} 1) {$b'_1$};
\node[\colourmakedark{3},left,xshift=-2pt] at (spath cs:\ttbranch{b}{4} 1) {$b'_h$};
\node[\colourmakedark{2},below right,yshift=-.35cm,xshift=2pt,label node] at (spath cs:\ttbranch{b}{1} 1) {$\bar{b}'$};
\draw[pin line] (label) -- (spath cs:\ttbranch{c}{1} .6);
\node[\colourmakedark{1},below,yshift=1pt] at (c-1) {\contour{grey background}{$s'$}};
\node[\colourmakedark{1},left,xshift=-2pt] at (spath cs:\ttbranch{b}{6} 1) {$\tau'$};
\coordinate (p) at ($(b-1)-(0,.75)$);
\draw[palette 4,line,-{Triangle[]}] ++(p) -- +(0,-.6) node[\colourmakedark{4},above right] {$\omega'$};
\fill[black] (p) circle (.3pt);
\end{scope}
\draw[line,-{Triangle[]}] (3.5,0) -- (4.5,0) node[anchor=base,yshift=6pt,midway] {\splitmove{}};
\end{tikzpicture}
\egroup

%% file: figures-source/split-move-branch-disappear.tex
\bgroup
\def\tikzttsetup#1{
\tikzset{
tt/new=#1,
tt/#1/.cd,
set sep=.5cm,
set width=.5cm,
add switch={0,0},
add switch={1.8,.8},
set switch thickness={2}{.75cm},
add branch standard={1}{2},
add bunch={1}{left}{5,spread=1.2cm,sep=.2cm},
add bunch={1}{below 1}{2,spread=1.2cm,sep=.3cm},
add bunch={2}{above 1}{1,spread=1.2cm,sep=.7cm},
add bunch={2}{right}{3,spread=.8cm,sep=.3cm},
set branch weights={2=1,3=1,4=1,5=2,6=3,1=3,7=3,8=2,9=3,10=2,11=2,12=2},
finalise,
}
\fill[grey background] [tt/#1/use branch background all but={2,3,4,1},tt/#1/use switch background all];
\fill[\colourmakebackground{3}] [tt/#1/use branch background={2,3,4}];
\fill[\colourmakebackground{2}] [tt/#1/use branch background={1}];
\draw[black,line] [tt/#1/use branch contour all but={2,3,4,1},tt/#1/use switch contour all];
\draw[\colourmakedark{3},line] [tt/#1/use branch contour={2,3,4}];
\draw[\colourmakedark{2},line] [tt/#1/use branch contour={1}];
\foreach \b/\dir[count=\i] in {2/right contour,4/left contour,1/right contour,1/left contour} {
    \path[spath/save=tmp] [tt/#1/use branch \dir/.expanded=\b];
    \coordinate (\i) at (spath cs:tmp 0) {};
}
\draw[line,dotted,dash expand off,black,opacity=.5] (2) -- (3);
}
\begin{tikzpicture}
\tikzttsetup{a}
\node[\colourmakedark{3},left,xshift=-2pt,yshift=2pt] at (spath cs:\ttbranch{a}{2} 1) {$b_1$};
\node[\colourmakedark{3},left,xshift=-2pt] at (spath cs:\ttbranch{a}{4} 1) {$b_h$};
\node[\colourmakedark{2},below right,yshift=-.35cm,xshift=2pt,label node] at (spath cs:\ttbranch{a}{1} 1) {$\bar{b}$};
\draw[pin line] (label) -- (spath cs:\ttbranch{a}{1} .8);
\node[black,above=.5cm] at (a-1) {$s$};
\node[black,left,xshift=-2pt] at (spath cs:\ttbranch{a}{6} 1) {$\tau$};
\coordinate (p) at ($(a-1)-(0,.75)$);
\draw[palette 4,line,-{Triangle[]}] ++(p) -- +(-.6,0) node[\colourmakedark{4},below] {$\eta$};
\draw[palette 4,line,-{Triangle[]}] ++(p) -- +(0,-.6) node[\colourmakedark{4},above right] {$\omega$};
\fill[black] (p) circle (.3pt);
\draw[line,palette 4,decorate,decoration={brace,amplitude=1pt,raise=1pt}] (1) -- (2);
\draw[line,palette 4,decorate,decoration={brace,amplitude=1pt,raise=1pt}] (3) -- (4);
\draw[line,palette 4,double distance=1pt,nfold] ++([xshift=-2pt]$(1)!.5!(3)$) -- +(4pt,0);

\begin{scope}[shift={(7,0)}]
\tikzttsetup{b}
\tikzset{
tt/b/carried={
name=c,
switch={1}{-.1cm}{4,5}{4,5},
switch={2}{0}{1,...,4}{1,2,3},
1=3,2=1,3=1,4=1,5=1,6=1,7=1,8=1,9=1,10=1,11=1,12=1,
},
}
\draw[palette 3,line] [tt/c/use branch={1,2,3}];
\draw[palette 1,line] [tt/c/use branch all but={1,2,3}];
\pic[palette 1] at (c-1) {dot};
\pic[palette 1] at (c-2) {dot};
\node[\colourmakedark{3},left,xshift=-2pt,yshift=2pt] at (spath cs:\ttbranch{b}{2} 1) {$b'_1$};
\node[\colourmakedark{3},left,xshift=-2pt] at (spath cs:\ttbranch{b}{4} 1) {$b'_h$};
\node[\colourmakedark{1},below,yshift=1pt] at (c-1) {\contour{grey background}{$s'$}};
\node[\colourmakedark{1},left,xshift=-2pt] at (spath cs:\ttbranch{b}{6} 1) {$\tau'$};
\coordinate (p) at ($(b-1)-(0,.75)$);
\draw[palette 4,line,-{Triangle[]}] ++(p) -- +(0,-.6) node[\colourmakedark{4},above right] {$\omega'$};
\fill[black] (p) circle (.3pt);
\end{scope}
\draw[line,-{Triangle[]}] (3.5,0) -- (4.5,0) node[anchor=base,yshift=6pt,midway] {\splitmove{}};
\end{tikzpicture}
\egroup

%% file: figures-source/twist-move.tex
\def\tikzttsetup#1{
\tikzset{
tt/new=#1,
tt/#1/.cd,
set sep=.5cm,
set width=1cm,
add switch={0,0},
set switch thickness={1}{1.5cm},
add bunch={1}{left}{4,spread=1.2cm,sep=.2cm},
add bunch={1}{right}{3,spread=1.2cm,sep=.3cm},
set branch weights={1=1,2=1,3=7,4=2,5=7,6=3,7=1},
finalise,
debug
}
\fill[grey background] [tt/#1/use branch background all but={1,2,3,5},tt/#1/use switch background all];
\fill[\colourmakebackground{3}] [tt/#1/use branch background={1,2}];
\fill[\colourmakebackground{2}] [tt/#1/use branch background={3,5}];
\draw[black,line] [tt/#1/use branch contour all but={1,2,3,5},tt/#1/use switch contour all];
\draw[\colourmakedark{3},line] [tt/#1/use branch contour={1,2}];
\draw[\colourmakedark{2},line] [tt/#1/use branch contour={3,5}];
}
\begin{tikzpicture}

\tikzttsetup{a}
\node[\colourmakedark{3},left,xshift=-2pt,yshift=4pt] at (spath cs:\ttbranch{a}{1} 1) {$b_1$};
\node[\colourmakedark{3},left,xshift=-2pt] at (spath cs:\ttbranch{a}{2} 1) {$b_h$};
\node[\colourmakedark{2},left,xshift=-2pt] at (spath cs:\ttbranch{a}{3} 1) {$\bar{b}$};
\node[\colourmakedark{2},right,xshift=2pt] at (spath cs:\ttbranch{a}{5} 1) {$\bar{b}$};
\node[black,above=.8cm] at (a-1) {$s$};
\node[black,left,xshift=-2pt] at (spath cs:\ttbranch{a}{4} 1) {$\tau$};
\coordinate (p) at ($(a-1)-(0,1.05)$);
\draw[palette 4,line,-{Triangle[]}] ++(p) -- +(-.6,0) node[\colourmakedark{4},below] {$\eta$};
\draw[palette 4,line,-{Triangle[]}] ++(p) -- +(0,-.6) node[\colourmakedark{4},above right] {$\omega$};
\fill[black] (p) circle (.3pt);
\foreach \b/\dir[count=\i] in {1/right contour,2/left contour} {
    \path[spath/save=tmp] [tt/a/use branch \dir/.expanded=\b];
    \coordinate (\i) at (spath cs:tmp 0) {};
}
\tikzset{span/.style={line,|<->|,arrows={[scale=.4]},shorten >=.5pt,shorten <=.5pt}}
\draw[palette 4,span] (1) -- (2);
\foreach \i[evaluate=\i as \j using int(\i+1)] in {2,3,4} {
    \coordinate (\j) at ($(\i)+(2)-(1)$);
    \draw[palette 1,span] (\i) -- (\j);
}
\draw[line,\colourmakedark{1},decorate,decoration={brace,amplitude=1pt,raise=3pt}] (2) -- (5) node[midway,right=4pt] {$k$};

\begin{scope}[shift={(6.5,0)}]
\tikzttsetup{b}
\tikzset{
tt/b/carried={
name=c,
switch={1}{-.35cm}{7,...,10}{7,8,9},
1=1,2=1,3=7,4=1,5=7,6=1,7=1,
},
}
\draw[palette 3,line] [tt/c/use branch={1,3,5,7}];
\draw[palette 3!70!black,line] [tt/c/use branch={2,4,6,8}];
\draw[palette 2,line] [tt/c/use branch={11,9}];
\draw[palette 1,line] [tt/c/use branch={10,12,13}];
\pic[palette 1] at (c-1) {dot};
\node[\colourmakedark{3},left,xshift=-2pt,yshift=4pt] at (spath cs:\ttbranch{b}{1} 1) {$b'_1$};
\node[\colourmakedark{3}!70!black,left,xshift=-2pt] at (spath cs:\ttbranch{b}{2} 1) {$b'_h$};
\node[\colourmakedark{2},left,xshift=-4pt,label node] at (spath cs:\ttbranch{b}{3} 1) {$\bar{b}'$};
\draw[pin line] (label) -- (spath cs:\ttbranch{c}{9} .9);
\node[\colourmakedark{2},right,xshift=4pt,label node] at (spath cs:\ttbranch{b}{5} 1) {$\bar{b}'$};
\draw[pin line] (label) -- (spath cs:\ttbranch{c}{11} .9);
\node[\colourmakedark{1},below,yshift=1pt] at (c-1) {\contour{grey background}{$s'$}};
\node[\colourmakedark{1},left,xshift=-2pt] at (spath cs:\ttbranch{b}{4} 1) {$\tau'$};
\coordinate (p) at ($(b-1)-(0,1.05)$);
\draw[palette 4,line,-{Triangle[]}] ++(p) -- +(0,-.6) node[\colourmakedark{4},above right] {$\omega'$};
\fill[black] (p) circle (.3pt);
\end{scope}
\draw[line,-{Triangle[]}] ++(2.7,0) -- +(1,0) node[anchor=base,yshift=6pt,midway] {\twistmove{}};
\end{tikzpicture}

%% file: figures-source/aht-sequence-typical-carrying.tex
\begin{tikzpicture}
\tikzset{
tt/new=a,
tt/a/.cd,
set sep=.2cm,
set spread=1.5cm,
set width=1cm,
add switch={0,0},
set switch thickness={1}{2cm},
add bunch={1}{left}{5},
add bunch={1}{right}{4,sep=.25cm},
set branch weights={1=3,2=1,3=2,4=3,5=2,6=2,7=5,8=1,9=3},
finalise,
carried={
    name=b,
    switch={1}{-.6cm}{7,8,9}{7,8},
    1=3,2=1,3=1,4=3,5=1,6=1,7=5,8=1,9=1,
}
}
\fill[grey background] [tt/a/use branch background all,tt/a/use switch background all];
\draw[black,line] [tt/a/use branch contour all,tt/a/use switch contour all];
\draw[line,palette 1] [tt/b/use branch all];
\pic[palette 1] at (b-1) {dot};
\node[black,left,xshift=-2pt] at (spath cs:\ttbranch{a}{1} 1) {$\tau_0$};
\node[\colourmakedark{1},left,xshift=-2pt] at (spath cs:\ttbranch{a}{5} 1) {$\tau_i$};
\node[black,above=1.1cm] at (a-1) {$s_0$};
\node[\colourmakedark{1},below,yshift=0pt] at (b-1) {\contour{grey background}{$s_i$}};
\coordinate (p) at ($(a-1)-(0,1.3)$);
\draw[palette 4,line,-{Triangle[]}] ++(p) -- +(0,-.6) node[\colourmakedark{4},above right] {$\omega_0$};
\fill[black] (p) circle (.3pt);
\end{tikzpicture}

%% file: figures-source/deep-nesting-one-step-no-monogon-a.tex
\begin{tikzpicture}
\tikzset{
tt/new=a,
tt/a/.cd,
set sep=.2cm,
set spread=1.5cm,
set width=1cm,
add switch={0,0},
set switch thickness={1}{2cm},
add bunch={1}{left}{5},
add bunch={1}{right}{5},
set branch weights={1=5,2=4,3=4,4=3,5=3,6=3,7=5,8=3,9=4,10=4},
finalise,
carried={
    name=b,
    switch={1}{-.6cm}{9,...,13}{9,...,13},
    1=3,7=3,
    2=3,9=3,
    3=3,10=3,
    4=2,6=2,
    5=2,8=2,
},
}
\fill[grey background] [tt/a/use branch background all but={3,10},tt/a/use switch background all];
\fill[\colourmakebackground{2}] [tt/a/use branch background={3,10}];
\draw[black,line] [tt/a/use branch contour all,tt/a/use switch contour all];
\draw[line,palette 1] [tt/b/use branch all but={9,16}];
\draw[line,palette 2] [tt/b/use branch={9,16}];
\pic[palette 1] at (b-1) {dot};
\node[black,left,xshift=-2pt] at (spath cs:\ttbranch{a}{1} 1) {$\tau_0$};
\node[\colourmakedark{1},left,xshift=-2pt] at (spath cs:\ttbranch{a}{5} 1) {$\tau_{i-1}$};
\node[\colourmakedark{2},left,xshift=-2pt] at (spath cs:\ttbranch{a}{3} 1) {$b$};
\node[black,right,xshift=2pt] at (spath cs:\ttbranch{a}{10} 1) {$\mathcolor{\colourmakedark{2}}{b}=\mathcolor{\colourmakedark{3}}{c}$};
\node[black,above=1.1cm] at (a-1) {$s_0$};
\node[\colourmakedark{1},below,yshift=0pt] at (b-1) {\contour{grey background}{$s_{i-1}$}};
\coordinate (p) at ($(a-1)-(0,1.3)$);
\draw[palette 4,line,-{Triangle[]}] ++(p) -- +(0,-.6) node[\colourmakedark{4},above right] {$\omega_0$};
\fill[black] (p) circle (.3pt);
\end{tikzpicture}

%% file: figures-source/deep-nesting-one-step-no-monogon-b.tex
\begin{tikzpicture}
\tikzset{
tt/new=a,
tt/a/.cd,
set spread=1.5cm,
set width=1cm,
add switch={0,0},
set switch thickness={1}{2cm},
add bunch={1}{left}{4,sep=.25cm},
add bunch={1}{right}{6,sep=.15cm},
set branch weights={1=7,2=2,3=7,4=4,5=2,6=4,7=4,8=3,9=4,10=3},
finalise,
carried={
    name=b,
    switch={1}{-.4cm}{5,...,8}{5,7,9,10,11,12},
    1=3,3=3,
    2=1,5=1,
    4=1,6=1,
    7=4,9=4,
    8=1,10=1,
},
}
\fill[grey background] [tt/a/use branch background all but={7,9},tt/a/use switch background all];
\fill[\colourmakebackground{2}] [tt/a/use branch background={7,9}];
\draw[black,line] [tt/a/use branch contour all,tt/a/use switch contour all];
\draw[line,palette 1] [tt/b/use branch all but={10,14,9,15}];
\draw[line,palette 2] [tt/b/use branch={10,14}];
\draw[line,palette 3] [tt/b/use branch={9,15}];
\draw[line over=grey background,palette 3,spath/save=tmp]  (spath cs:\ttbranch{b}{9} 0) to[out=180,in=180,in looseness=2] (spath cs:\ttbranch{b}{15} 0);
\pic[palette 1] at (b-1) {dot};
\node[\colourmakedark{3},above=.2cm,xshift=-4pt,label node] at (b-1) {\contour{grey background}{$c$}};
\draw[pin line] (label) -- (spath cs:tmp .5);
\node[black,left,xshift=-2pt] at (spath cs:\ttbranch{a}{1} 1) {$\tau_0$};
\node[\colourmakedark{1},left,xshift=-2pt] at (spath cs:\ttbranch{a}{4} 1) {$\tau_{i-1}$};
\node[\colourmakedark{2},right,xshift=6pt,label node] at (spath cs:\ttbranch{a}{8} 1) {$b$};
\draw[pin line] (label) -- (spath cs:\ttbranch{b}{10} .9);
\draw[pin line] (label) -- (spath cs:\ttbranch{b}{14} .9);
\node[black,above=1.1cm] at (a-1) {$s_0$};
\node[\colourmakedark{1},below,yshift=-2pt] at (b-1) {\contour{grey background}{$s_{i-1}$}};
\coordinate (p) at ($(a-1)-(0,1.3)$);
\draw[palette 4,line,-{Triangle[]}] ++(p) -- +(0,-.6) node[\colourmakedark{4},above right] {$\omega_0$};
\fill[black] (p) circle (.3pt);
\end{tikzpicture}

%% file: figures-source/distance-algorithm-for-carried-curves-deep-nesting-a.tex
\begin{tikzpicture}
\tikzset{
tt/new=a,
tt/a/.cd,
set spread=1.5cm,
set sep=.3cm,
set width=1cm,
add switch={0,0},
set switch thickness={1}{1.8cm},
add bunch={1}{left}{3},
add bunch={1}{right}{3},
set branch weights={1=3,2=8,3=9,4=9,5=3,6=8},
finalise,
carried={
    name=b,
    switch={1}{.4cm}{2,3,4}{2,3,4},
    1=1,5=1,
    2=3,6=3,
    3=6,4=6,
},
}
\fill[grey background] [tt/a/use branch background all but={2,4},tt/a/use switch background all];
\fill[\colourmakebackground{2}] [tt/a/use branch background={2,4}];
\draw[black,line] [tt/a/use branch contour all,tt/a/use switch contour all];
\draw[line,palette 1] [tt/b/use branch all];
\pic[palette 1] at (b-1) {dot};
\node[black,left,xshift=-2pt] at (spath cs:\ttbranch{a}{3} 1) {$\sigma$};
\node[\colourmakedark{1},left,xshift=-2pt] at (spath cs:\ttbranch{a}{1} 1) {$\sigma'$};
\end{tikzpicture}

%% file: figures-source/distance-algorithm-for-carried-curves-deep-nesting-b.tex
\begin{tikzpicture}
\tikzset{
tt/new=a,
tt/a/.cd,
set spread=1.5cm,
set sep=.2cm,
set width=1cm,
add switch={0,0},
set switch thickness={1}{1.8cm},
add bunch={1}{left}{4},
add bunch={1}{right}{4},
set branch weights={1=3,2=8,3=2,4=9,5=9,6=2,7=3,8=8},
finalise,
carried={
    name=b,
    switch={1}{.4cm}{2,3,4,5}{2,3,4,5}, 
    1=1,7=1,
    2=4,8=4,
    4=8,5=8,
    3=1,6=1,
},
}
\fill[grey background] [tt/a/use branch background all but={3,6},tt/a/use switch background all];
\fill[\colourmakebackground{4}] [tt/a/use branch background={3,6}];
\draw[black,line] [tt/a/use branch contour all,tt/a/use switch contour all];
\draw[line,palette 1] [tt/b/use branch all but={3,6,9,12,17}];
\draw[line,palette 3] [tt/b/use branch={3,6,9,12,17}];
\pic[palette 1] at (b-1) {dot};
\node[black,left,xshift=-2pt] at (spath cs:\ttbranch{a}{4} 1) {$\sigma$};
\node[\colourmakedark{1},left,xshift=-2pt] at (spath cs:\ttbranch{a}{1} 1) {$\sigma'$};
\node[\colourmakedark{4},right,xshift=2pt] at (spath cs:\ttbranch{a}{6} 1) {$\rho\setminus \sigma$};
\node[\colourmakedark{3},right,xshift=6pt,label node,yshift=.2cm] at (spath cs:\ttbranch{a}{8} 1) {$\rho'\setminus\sigma'$};
\draw[pin line] (label) -- (spath cs:\ttbranch{b}{12} .1);
\end{tikzpicture}

%% file: figures-source/find-carrying-train-track-pattern.tex
\begin{tikzpicture}
\edef\myopacity{1}

\tikzset{
triangle/.code n args={3}{
    \begin{scope}[name prefix=triangle]
        \coordinate (a) at (#1);
        \coordinate (b) at (#2);
        \coordinate (c) at (#3);
        \tikzmath{
            coordinate \p;
            \p=(a)-(b);
            \a1=atan2(\py,\px)-90;
            \p=(c)-(b);
            \a2=atan2(\py,\px)+90;
            \p=(a)-(c);
            \a3=atan2(\py,\px)+90;
        }
        \draw[opacity=\myopacity,palette 1,line,every to/.style={looseness=.6}] ($(a)!.5!(b)$) to[out=\a1,in=\a2] ($(b)!.5!(c)$) to[out=\a2,in=\a3] ($(c)!.5!(a)$) to[out=\a3,in=\a1] ($(a)!.5!(b)$);
    \end{scope}
    \edef\myopacity{.3}
}}

\coordinate (1) at (90:1.2);
\coordinate (2) at (-150:1.2);
\coordinate (3) at (-30:1.2);
\coordinate (4) at (1.5,1);
\coordinate (5) at (170:1.5);
\coordinate (6) at (-80:1.5);
\coordinate (7) at (125:1.6);

\tikzset{
triangle={1}{2}{3},
triangle={1}{3}{4},
triangle={1}{5}{2},
triangle={2}{6}{3},
triangle={1}{7}{5},
every path/.style={black,line,line join=round,opacity=.3}
}
\draw[opacity=1] (1) -- (2) -- (3) -- (1);
\draw (1) -- (7) -- (5) -- (1) (5) -- (2) -- (6) -- (3) -- (4) -- (1);

\end{tikzpicture}

%% file: figures-source/find-carrying-train-track-measure.tex
\begin{tikzpicture}
\foreach \i in {1,2,3} {\coordinate (\i) at (-30+120*\i:1.2);}

\foreach \i/\j/\k in {1/2/7,2/3/4,3/1/5} {
    \foreach \l in {1,...,\k} {
        \tikzmath{\t=\l/(\k+1)*.5+.25;}
        \coordinate (\i-\l) at ($(\i)!\t!(\j)$);
    }
}
{
    \tikzset{every path/.style={line,draw=palette 2}}
    \foreach \i in {1,...,4} {
        \tikzmath{\j=6-\i;}
        \draw (1-\i) to[out=-30,in=-150] (3-\j);
    }
    \foreach \i in {1,...,3} {
        \tikzmath{\j=8-\i;}
        \draw (1-\j) to[out=-30,in=90] (2-\i);
    }
    \draw (2-4) to[out=90,in=-150] (3-1);
}

\draw[black,line] (1) -- (2) -- (3) -- (1);

{
    \tikzset{brace/.style={draw=\colourmakedark{2},line,decorate,decoration={brace,mirror,raise=2pt,amplitude=1pt},every node/.style={midway,\colourmakedark{2}}}}
    \draw[brace] ([shift=(60:1pt)]1-1) -- ([shift=(60:-1pt)]1-4) node[midway,above left=1pt] {$x_1$};
    \draw[brace] ([shift=(60:1pt)]1-5) -- ([shift=(60:-1pt)]1-7) node[midway,above left=1pt] {$x_2$};
    \draw[brace] ([shift=(120:-2.5pt)]3-1) -- ([shift=(120:2.5pt)]3-1) node[midway,above right=1pt] {$x_3$};
}

\draw[xshift=1.75cm,yshift=.3cm,{Triangle[]}-{Triangle[]}] (-.4,0) --(.4,0);

\begin{scope}[xshift=3.5cm]
\foreach \i in {1,2,3} {\coordinate (\i) at (-30+120*\i:1.2);}

\draw[line,palette 1,line join=round] ($(1)!.5!(2)$) to[out=-30,in=90] ($(2)!.5!(3)$) to[out=90,in=-150] ($(3)!.5!(1)$) to[out=-150,in=-30] cycle;
\foreach \i in {1,2,3} {
    \node[\colourmakedark{1}] at ({-30+120*\i}:.5) {$x_\i$};
}

\draw[line] (1) -- (2) -- (3) -- (1);
\end{scope}

\end{tikzpicture}

%% file: figures-source/find-carrying-train-track-monogon-and-bigon-a.tex
\begin{tikzpicture}[every node/.style={scale=1}]
\tikzset{ovt/new={
    n=8,
    radius=1.8cm,
    arcs={4,2,1,2,3,1,2,1,1,3,0,2,3,2,2,0,5,2,1,1,2,3,0,3}
}}
\draw[black,line] [ovt/use radial edge all];
\draw[palette 3,line] [ovt/use arc all];
\draw[palette 1,line over={white},spath/save=tt] [ovt/use train track];
\draw[black,line] [ovt/use outer edge all];
\node[below left=2pt,\colourmakedark{3},label node] at ($(ovt-corner-4)!.5!(ovt-corner-5)$) {$a$};
\draw[pin line] (label) -- (spath cs:ovt-arc-26 .7);
\node[below right=2pt,\colourmakedark{1},label node] at ($(ovt-corner-6)!.5!(ovt-corner-7)$) {$\tau$};
\draw[pin line] (label) -- (spath cs:tt 0.83);
\pic[palette 2] at (ovt-centre) {dot};
\node[\colourmakedark{2},above right] at (ovt-centre) {\contour{white}{$v$}};
\end{tikzpicture}

%% file: figures-source/find-carrying-train-track-monogon-and-bigon-b.tex
\begin{tikzpicture}
\tikzset{ovt/new={
    n=8,
    arcs={0,3,3,3,0,2,1,4,5,3,3,0,1,2,1,0,2,2,2,0,3,3,2,0}
}}
\draw[black,line] [ovt/use radial edge all];
\draw[palette 3,line] [ovt/use arc all];
\draw[palette 1,line over={white},spath/save=tt] [ovt/use train track];
\draw[black,line] [ovt/use outer edge all];
\node[below left=2pt,\colourmakedark{3},label node] at ($(ovt-corner-4)!.5!(ovt-corner-5)$) {$a$};
\draw[pin line] (label) -- (spath cs:ovt-arc-29 .7);
\node[below right=2pt,\colourmakedark{1},label node] at ($(ovt-corner-6)!.5!(ovt-corner-7)$) {$\tau$};
\draw[pin line] (label) -- (spath cs:tt 0.86);
\pic[palette 2] at (ovt-centre) {dot};
\node[\colourmakedark{2},above right] at (ovt-centre) {\contour{white}{$v$}};
\end{tikzpicture}

%% file: figures-source/find-carrying-train-track-monogon-a.tex
\begin{tikzpicture}
\tikzset{ovt/new={
    n=8,
    arcs={2,1,0,2,0,0,2,0,1,3,0,1,2,2,2,0,4,4,3,1,2,3,2,0}
}}
\draw[black,line] [ovt/use radial edge all but={3,6}];
\draw[palette 4,line] [ovt/use radial edge={3,6}];
\draw[palette 3,line] [ovt/use arc all but={1,4,6,9,13,22,26,27,33}];
\draw[palette 2,line] [ovt/use arc={1,4,6,9,13,22,26,27,33}];
\draw[black,line] [ovt/use outer edge all];
\tikzset{
    ovt/.cd,
    label in corner={3}{3}{\colourmakedark{4}}{$e$},
    every label in corner/.style={\colourmakedark{1}},
    label in corner={5}{1}{}{$n$},
    label in corner={6}{1}{}{$1$},
    label in corner={6}{2}{}{$h$},
    label in corner={6}{3}{}{$k$},
}
\node[below left=2pt,\colourmakedark{3},label node] at ($(ovt-corner-4)!.5!(ovt-corner-5)$) {$a$};
\draw[pin line] (label) -- (spath cs:ovt-arc-15 .7);
\node[below=4pt,\colourmakedark{2},label node] at ($(ovt-corner-5)!.5!(ovt-corner-6)$) {$\alpha$};
\draw[pin line] (label) -- (spath cs:ovt-arc-22 0.8);
\pic[black] at (ovt-centre) {dot};
\end{tikzpicture}

%% file: figures-source/find-carrying-train-track-monogon-b.tex
\begin{tikzpicture}
\tikzset{ovt/new={
    n=8,
    arcs={2,1,0,2,0,0,2,0,1,3,0,1,2,2,2,0,4,4,3,1,2,3,2,0}
}}
\draw[black,line] [ovt/use radial edge all but={3,6}];
\draw[palette 4,line] [ovt/use radial edge={3,6}];
\draw[palette 3,line] [ovt/use arc all but={1,4,6,9,13,22,26,27,33}];
\draw[palette 2,line,spath/save=tmp] ++(spath cs:ovt-arc-22 1) -- +(0,.65) to[out=90,in=-90] +(-.08,.9) to[out=90,in=180] (0,-.75);
\draw[palette 2,line] [spath/use={tmp,transform={xscale=-1}}];
\draw[black,line] [ovt/use outer edge all];
\tikzset{
    ovt/.cd,
    label in corner={3}{3}{\colourmakedark{4}}{$e$},
    every label in corner/.style={\colourmakedark{1}},
    label in corner={5}{1}{}{$n$},
    label in corner={6}{1}{}{$1$},
    label in corner={6}{2}{}{$h$},
    label in corner={6}{3}{}{$k$},
}
\node[below left=2pt,\colourmakedark{3},label node] at ($(ovt-corner-4)!.5!(ovt-corner-5)$) {$a$};
\draw[pin line] (label) -- (spath cs:ovt-arc-15 .7);
\pic[black] at (ovt-centre) {dot};
\node[below=4pt,\colourmakedark{2},label node] at ($(ovt-corner-5)!.5!(ovt-corner-6)$) {\phantom{$\alpha$}};
\end{tikzpicture}

%% file: figures-source/find-carrying-train-track-bigon-a.tex
\begin{tikzpicture}
\tikzset{ovt/new={
    n=8,
    arcs={0,2,2,2,0,2,4,0,0,4,0,0,3,1,1,0,4,3,3,0,0,2,1,0}
}}
\draw[black,line] [ovt/use radial edge all];
\draw[palette 3,line] [ovt/use arc all but={4,5,9,13,17,25,2,32,29,28}];
\draw[palette 2,line] [ovt/use arc={4,5,9,13,17,25}];
\draw[black,line] [ovt/use outer edge all];
\tikzset{
    ovt/.cd,
    every label in corner/.style={\colourmakedark{1}},
    label in corner={1}{1}{}{$j$},
    label in corner={1}{2}{}{$l$},
    label in corner={1}{3}{}{$m$},
    label in corner={5}{1}{}{$n$},
    label in corner={6}{1}{}{$1$},
    label in corner={6}{2}{}{$h$},
    label in corner={6}{3}{}{$k$},
}
\node[below left=2pt,\colourmakedark{3},label node] at ($(ovt-corner-4)!.5!(ovt-corner-5)$) {$a$};
\draw[pin line] (label) -- (spath cs:ovt-arc-21 .7);
\pic[black] at (ovt-centre) {dot};
\node[below=4pt,\colourmakedark{2},label node] at ($(ovt-corner-5)!.5!(ovt-corner-6)$) {$\alpha$};
\draw[pin line] (label) -- (spath cs:ovt-arc-25 0.8);
\end{tikzpicture}

%% file: figures-source/find-carrying-train-track-bigon-b.tex
\begin{tikzpicture}
\tikzset{ovt/new={
    n=8,
    arcs={0,2,2,2,0,2,4,0,0,4,0,0,3,1,1,0,4,3,3,0,0,2,1,0}
}}
\draw[black,line] [ovt/use radial edge all];
\draw[palette 3,line] [ovt/use arc all but={4,5,9,13,17,25,2,32,29,28}];
\draw[palette 2,line] [ovt/use arc={2,32,29,28}];
\draw[black,line] [ovt/use outer edge all];
\tikzset{
    ovt/.cd,
    every label in corner/.style={\colourmakedark{1}},
    label in corner={1}{1}{}{$j$},
    label in corner={1}{2}{}{$l$},
    label in corner={1}{3}{}{$m$},
    label in corner={5}{1}{}{$n$},
    label in corner={6}{1}{}{$1$},
    label in corner={6}{2}{}{$h$},
    label in corner={6}{3}{}{$k$},
}
\node[below left=2pt,\colourmakedark{3},label node] at ($(ovt-corner-4)!.5!(ovt-corner-5)$) {$a$};
\draw[pin line] (label) -- (spath cs:ovt-arc-21 .7);
\pic[black] at (ovt-centre) {dot};
\node[below=4pt,label node] at ($(ovt-corner-5)!.5!(ovt-corner-6)$) {\phantom{$\alpha$}};
\end{tikzpicture}

%% file: nielsen-thurston-classification.tex
\section{\texorpdfstring{Nielsen\=/Thurston}{Nielsen-Thurston} classification}
\label{sec:classification}

A remarkable theorem of Thurston (see \cite{fathi-laudenbach-poenaru}) states that every mapping class in $\MCG(S)$ is either periodic, reducible (and non\=/periodic), or pseudo\=/Anosov; the category to which a mapping class belongs is called its Nielsen\=/Thurston type.
We will say that a homeomorphism of $S$ is periodic, reducible, or pseudo\=/Anosov if its isotopy class in $\MCG(S)$ is.
In this section, we endeavour to describe an efficient algorithm to decide the Nielsen\=/Thurston type of a surface homeomorphism, by studying its action on the curve graph of $S$ and, in particular, its stable translation length (see \cref{sec:classification:translation length}).
Naturally, the distance algorithm devised in \cref{sec:distance algorithm} will play a crucial role here.

The argument to turn a distance algorithm for the curve graph into a surface homeomorphism classification algorithm is not novel.
In fact, this section closely follows \cite[\S4.2]{bell-webb-algorithms}.
The three main differences are that we write down the proofs in greater detail, we make an effort to keep the constants and the running time bounds polynomial in $\xi$, and we adapt the arguments to make use of a coarse distance algorithm rather than an exact one.
We also remark that, unlike \citeauthor{bell-webb-algorithms}'s, our algorithm does not produce a fixed multicurve in the reducible case.

\subsection{Periodic homeomorphisms}

Of the three Nielsen\=/Thurston types, periodic homeomorphisms are the easiest to identify, thanks to the following bound on their order.

\newconstant{periodic order bound}{C_{\mathrm{per}}}{Cper}{-6\chi}
\begin{proposition}\label{thm:periodic homeomorphism order bound}
If a homeomorphism $\map{f}{S}{S}$ is periodic, then $f^k$ is isotopic to the identity for some $1\le k\le \constant{periodic order bound}$, where
\[
    \declareconstant{periodic order bound}=\constantvalue{periodic order bound}.
\]
\end{proposition}
\begin{proof}
If $S$ is closed and has genus $g$, then it was proved by \textcite{wiman} that the order of $f$ in $\MCG(S)$ is bounded above by $4g+2$, which is smaller than $-6\chi$ for $g\ge 2$.
Otherwise, let $f$ be a homeomorphism of the punctured surface $S$ with finite order $k$. By the Nielsen realisation theorem \cite{nielsen-realisation}, we can assume that $f$ is in fact an isometry of some hyperbolic metric on $S$, such that $f^k$ is the identity. If we let $H$ be the cyclic group of isometries of $S$ generated by $f$, then the quotient $S'=\quotient{S}{H}$ is a topological surface, and the projection $\map{\pi}{S}{S'}$ is a branched covering. The Riemann\=/Hurwitz formula implies that
\[
    -\frac{\chi(S)}{k}=-\chi(S')+\sum_{x\in S'}\left(1-\frac{1}{e_x}\right),
\]
where $e_x$ is the cardinality of the stabilisers in $H$ of the preimages of $x$ under $\pi$. Note that, since $S$ is punctured, then so is $S'$, and hence $\chi(S')\le 1$. A simple case analysis shows that, if the right\=/hand side of the equality is positive, then it is not smaller than $1/6$. The desired inequality $k\le-6\chi(S)$ follows immediately.
\end{proof}

Thanks to the bound of \cref{thm:periodic homeomorphism order bound}, we can use the Alexander method described in \cite[\S2.3]{farb-margalit} to detect periodic homeomorphisms.

\begin{proposition}\label{thm:periodic recognition}
Let $\map{f}{S}{S}$ be a homeomorphism. There is an algorithm to decide if $f$ is periodic. The running time of the algorithm is polynomial in $\xi$ and $\norm{f}$.
\end{proposition}
\begin{proof}
We show how to construct a small set $\AAA$ of short curves on $S$ such that if a homeomorphism $\umap{S}{S}$ fixes each curve in $\AAA$ up to isotopy then it has finite order.
By \cref{thm:periodic homeomorphism order bound}, this is enough to recognise periodic homeomorphisms.
\begin{substeps}
\item Suppose first that $S$ is punctured, and recall that the triangulation $\TTT$ is ideal.
In this case, we take $\AAA$ to be 
\[
    \AAA=\{a:\text{$e$ is an edge of $\TTT$, $a$ is boundary component of a neighbourhood of $e$ in $S$}\},
\]
where we think of punctures of $S$ as marked points on a closed surface.
The curves in $\AAA$ are easy to construct algorithmically, and have complexity $\bigO(\xi)$; moreover, $\card{\AAA}=\bigO(\xi)$.
Finally, we note that a homeomorphism $\umap{S}{S}$ fixes all the curves in $\AAA$ up to isotopy if and only if it fixes all the edges of $\TTT$ up to isotopy.
By the Alexander method \cite[Proposition 2.8]{farb-margalit}, a homeomorphism $\umap{S}{S}$ fixing the isotopy classes of the edges of an ideal triangulation of $S$ is isotopic to the identity.
\item The closed case is more involved.
Let $e_1,\ldots,e_m$ be the edges of the one\=/vertex triangulation $\TTT$.
We start by considering the set of curves
\[
\AAA_1=\left\{a:\,
\begin{matrix*}[l]
\text{$1\le i\le m$, $a$ is a boundary component of a}\\
\text{neighbourhood of $e_1\cup\ldots\cup e_i$ in $S$ and is essential}
\end{matrix*}
\right\}.
\]
We can assume that the curves in $\AAA_1$ are pairwise non\=/isotopic, since \cref{thm:deciding isotopy of curves} provides a procedure to remove duplicates.

The elements of $\AAA_1$ can be constructed by tracing the boundary of a neighbourhood of $e_1\cup\ldots\cup e_i$ and then normalising the resulting curves.
Since these curves are short -- they intersect the triangulation $\bigO(\xi)$ times -- the normalising procedure can be performed directly, without resorting to the machinery of \cite[Theorem 1.3]{lackenby-minimal-position}.

By \cite[Lemma 3.8]{lackenby-pants-graph}, the curves in $\AAA_1$ form a pants decomposition of the surface $S$.
There is a multicurve $A_2$ on $S$ which intersects every pair of pants in three arcs as shown in \cref{fig:periodic recognition:A_2}.
This multicurve can be constructed algorithmically by tracing three arcs in each pair of pants, then joining arcs in adjacent pairs of pants at their endpoints, and normalising the resulting curves.
Again, since the curves in $\AAA_1$ are short, these operations can be performed directly in polynomial time in $\xi$.
Moreover, the multicurve $A_2$ will have $\bigO(\xi)$ components, each with complexity $\bigO(\xi\cdot\log\xi)$.

\begin{figure}
\centering
    \tikzsetnextfilename{periodic-recognition-A_2}%
    \input{figures-source/periodic-recognition-A_2.tex}%

\caption{The multicurve $A_2$ intersects each pair of pants $P$ in the pants decomposition induced by $\AAA_1$ in three arcs, one for each pair of boundary components of $P$.}
\label{fig:periodic recognition:A_2}
\end{figure}
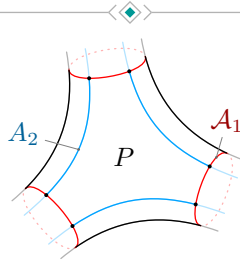

Finally, we let $\AAA_2$ be the set of components of $A_2$, and set $\AAA=\AAA_1\cup\AAA_2$.
We invoke the Alexander method \cite[Proposition 2.8]{farb-margalit} again to conclude that a homeomorphism $\umap{S}{S}$ fixing the isotopy class of each curve in $\AAA$ is periodic. \qedhere
\end{substeps}
\end{proof}

\subsection{Subsurface projection}\label{sec:classification:subsurface projection}
A key ingredient in \citeauthor{bell-webb-algorithms}'s algorithm -- and, therefore, in ours as well -- is the bounded geodesic image theorem, as proved by \textcite{masur-minsky-2} and later refined by \textcite{webb-effective-geometry}. We quickly introduce the relevant notations and definitions, following the first reference \cite[Sections 2.3 and 2.4]{masur-minsky-2}.

A subspace $Y\subs S$ is an \emph{essential surface} if it is a component of $S\setminus\Sigma$ for some essential multicurve $\Sigma\subs S$; note that we do not assume that $\xi(Y)\ge 2$. Suppose that $Y$ is not an annulus nor a sphere with $3$ punctures.
The \emph{subsurface projection map} is a function
\[
    \map{\pi_Y}{\C(S)}{\PPP\C_0(Y)},
\]
where $\PPP\C_0(Y)$ denotes the set of finite subsets of $\C_0(Y)$, satisfying the following properties for every $a\in\C_0(S)$:
\begin{enumarabic}
\item $\pi_Y(a)=\emptyset$ if and only if $a$ can be isotoped outside $Y$;
\item $f(\pi_Y(a))=\pi_{f(Y)}(f(a))\in\PPP\C_0(f(Y))$ for every homeomorphism $\map{f}{S}{S}$.
\end{enumarabic}
We refer the reader to \cite[Section 2.3]{masur-minsky-2} for the actual definition.

If $Y$ is an annulus, we define $\C_0(Y)$ to be the set of arcs in the closure $\closure{Y}$ of $Y$ connecting the two boundary components of $\closure{Y}$, modulo isotopies fixing $\boundary Y$ pointwise. The curve graph $\C(Y)$ is obtained by connecting two classes of arcs in $\C_0(S)$ by an edge of length $1$ if they admit representatives with disjoint interiors. Let $c$ be the core curve of $Y$. The subsurface projection map
\[
    \map{\pi_Y}{\C_0(S)}{\PPP\C_0(Y)}
\]
can be defined also in this case, and satisfies the following property: for every curve $a\in\C_0(S)$, the projection $\pi_Y(a)$ is empty if and only if $i(a,c)=0$. Moreover, the following lemma shows how the subsurface projection map behaves nicely under homeomorphisms that are isotopic to Dehn twists in a sufficiently large neighbourhood of $Y$.

\begin{lemma}\label{thm:dehn twist annular projection}
Let $Y\subs S$ be an essential annulus with core curve $c$, and $Z\subs S$ be an essential subsurface such that $Y\subs Z$ and $c$ is an essential curve in $Z$.
Let $\map{f}{S}{S}$ be a homeomorphism preserving $Z$ such that $f|_Z$, considered as a homeomorphism of $Z$, is isotopic to $\twist_c^k$ for some integer $k\neq 0$. Then
\[
    \dist[\C(Y)](\pi_Y(a),\pi_Y(f(a)))\ge|k|-1
\]
for every curve $a\in\C_0(S)$ such that $\pi_Y(a)\neq\emptyset$.
\end{lemma}
\begin{proof}
We recall from \cite[Section 2.4]{masur-minsky-2} the definition of $\pi_Y$ in the annular case, which is not needed anywhere in this article except for the proof of this lemma. Let us fix a hyperbolic metric for $S$. Let $\map{q}{A}{S}$ be the cover associated to the group $\pi_1Y<\pi_1 S$. The surface $A$ is topologically an open annulus, and can be compactified to a closed annulus $\bar{A}$ using the hyperbolic metric inherited from $S$. Then the subsurface projection map is the function
\[
    \map{\pi_Y}{\C_0(S)}{\PPP\C_0(A)}
\]
sending a curve $a\in\C_0(S)$ to the set of lifts of $a$ to $A$ whose closures in $\bar{A}$ connect the two components of $\boundary\bar{A}$. This is a finite set of diameter at most $1$ in $\C(A)$, and is empty precisely when $a$ can be isotoped away from $c$. We can think of $\pi_Y$ as a map to $\PPP\C_0(Y)$ if we fix an orientation\=/preserving homeomorphism $\umap{\bar{A}}{\closure{Y}}$.

Let $\bar{c}\subs A$ be the unique lift of $c$ which is a core curve of $\bar{A}$, and let $\bar{Y}$ and $\bar{Z}$ be the components of $q^{-1}(Y)$ and $q^{-1}(Z)$ respectively containing $\bar{c}$. Finally, let $\map{\bar{f}}{\bar{A}}{\bar{A}}$ be the lift of $f$ that preserves $\bar{c}$. Note that, since $c$ is essential in $Z$, the closure of $\bar{Z}$ in $\bar{A}$ intersects both boundary components of $\bar{A}$.

Up to isotopy, we can assume that $f$ restricts to the identity on $Z\setminus Y$. It follows that $\bar{f}$ is the identity on the two components of $\bar{Z}\setminus q^{-1}(Y)$ adjacent to $\bar{Y}$. In particular, the homeomorphism $\bar{f}$ will be the identity on 
\[
    W=\boundary\bar{Y}\cup b_1\cup b_2,
\]
where $b_1,b_2\subs\closure{\bar{Z}\setminus q^{-1}(Y)}$ are arcs connecting the two boundary components of $\bar{Y}$ to $\boundary\bar{A}$.
See \cref{fig:annular projection} for a graphical representation of the objects introduced so far.

\begin{figure}
\centering
\begin{subcaptionblock}{.45\linewidth}
\centering
    \tikzsetnextfilename{annular-subsurface-projection-a}%
    \input{figures-source/annular-subsurface-projection-a.tex}%

\caption{}
\label{fig:annular projection:a}
\end{subcaptionblock}
\begin{subcaptionblock}{.45\linewidth}
\centering
    \tikzsetnextfilename{annular-subsurface-projection-b}%
    \input{figures-source/annular-subsurface-projection-b.tex}%

\caption{}
\label{fig:annular projection:b}
\end{subcaptionblock}\\
\begin{subcaptionblock}{.45\linewidth}
\centering
    \tikzsetnextfilename{annular-subsurface-projection-c}%
    \input{figures-source/annular-subsurface-projection-c.tex}%

\caption{}
\label{fig:annular projection:c}
\end{subcaptionblock}
\begin{subcaptionblock}{.45\linewidth}
\centering
    \tikzsetnextfilename{annular-subsurface-projection-d}%
    \input{figures-source/annular-subsurface-projection-d.tex}%

\caption{}
\label{fig:annular projection:d}
\end{subcaptionblock}
\caption{\subref{fig:annular projection:a} Several lifts of $c$ to $\bar{A}$, including the core curve $\bar{c}$. \subref{fig:annular projection:b} The annulus $\bar{Y}\subs\bar{A}$. \subref{fig:annular projection:c} The lift $\bar{Z}$ of $Z$ containing $\bar{c}$. \subref{fig:annular projection:d} A subspace $W\subs\bar{A}$ on which $\bar{f}$ is the identity.}
\label{fig:annular projection}
\end{figure}

Let $\bar{a}$ be any arc in $\C_0(A)$.
Note that $A\setminus W$ is the disjoint union of two open discs and the annulus $\bar{Y}$; since $\bar{f}$ is the identity on $W$ and is the $k$\=/th power of the Dehn twist about $\bar{c}$ on $\bar{Y}$, we can isotope $\bar{f}(\bar{a})$, fixing $\boundary\bar{A}$, so that
\begin{itemize}
\item it does not intersect $\bar{a}$ in $A\setminus\closure{\bar{Y}}$,
\item it intersects $\bar{a}$ twice in $\boundary\bar{Y}$, and
\item it intersects $\bar{a}$ exactly $|k|-1$ times in $\bar{Y}$, always with the same sign,
\end{itemize}
as shown in \cref{fig:annular projection isotopy}

\begin{figure}
\centering
\begin{subcaptionblock}{.45\linewidth}
\centering
    \tikzsetnextfilename{annular-subsurface-projection-e}%
    \input{figures-source/annular-subsurface-projection-e.tex}%

\caption{}
\label{fig:annular projection isotopy:a}
\end{subcaptionblock}
\begin{subcaptionblock}{.45\linewidth}
\centering
    \tikzsetnextfilename{annular-subsurface-projection-f}%
    \input{figures-source/annular-subsurface-projection-f.tex}%

\caption{}
\label{fig:annular projection isotopy:b}
\end{subcaptionblock}
\caption{Isotoping $\bar{f}(\bar{a})$ (shown in \subref{fig:annular projection isotopy:a}) so that it intersects $\bar{a}$ at least $|k|-1$ times \subref{fig:annular projection isotopy:b}.}
\label{fig:annular projection isotopy}
\end{figure}

It follows that the absolute value of the algebraic intersection number of $\bar{a}$ and $\bar{f}(\bar{a})$ is bounded below by $|k|-1$ (note that the two intersections on $\boundary\bar{Y}$ might not be transverse). Equation (2.3) in \cite{masur-minsky-2} then implies that
\[
    \dist[\C(A)](\bar{a},\bar{f}(\bar{a}))\ge|k|.
\]

Finally, we deduce the inequality in the statement of the lemma by recalling that $\pi_Y(a)$ has diameter at most $1$ in $\C(A)$, and realising that $\pi_Y(f(a))=\bar{f}(\pi_Y(a))$.
\end{proof}

We can now state the bounded geodesic image theorem, combining Theorems 4.1.8 and 4.2.1 of \textcite{webb-effective-geometry}'s thesis.

\newconstant{bounded geodesic image}{D_{\mathrm{b.g.i}}}{Dbgi}{100}
\begin{theorem}\label{thm:bounded geodesic image}
Let $Y$ be an essential subsurface of $S$ which is not a sphere with $3$ punctures. Let $\gamma\subs\C(S)$ be a geodesic such that $\pi_Y(a)\neq\emptyset$ for every vertex $a\in\gamma\cap\C_0(S)$. Then
\[
    \diam[\C(Y)](\pi_Y(\gamma\cap\C_0(S)))\le\constant{bounded geodesic image},
\]
where
\[
    \declareconstant{bounded geodesic image}=\constantvalue{bounded geodesic image}.
\]
\end{theorem}

The reason why we take $\constant{bounded geodesic image}$ to be $\constantvalue{bounded geodesic image}$ is that we (following \citeauthor{masur-minsky-2}) define the subsurface projection map in the non\=/annular case to be the composition of \citeauthor{webb-effective-geometry}'s map $\map{\kappa_Y}{\C_0(S)}{\PPP\AAA\CCC_0(Y)}$ to finite subsets of the \emph{arc and curve graph} with the $2$\=/Lipschitz function $\map{\psi}{\AAA\CCC_0(Y)}{\PPP\C_0(Y)}$ introduced in \cite[Lemma 2.2]{masur-minsky-2}. Therefore, \citeauthor{webb-effective-geometry}'s bound of $50$ for non\=/annular subsurfaces must be doubled.

\subsection{Stable translation length}
\label{sec:classification:translation length}
\newconstant{min translation length}{\lambda_0}{l0}{\frac{1}{162\chi^2}}

The quantity we will use to distinguish reducible homeomorphisms from pseudo\=/Anosov ones is the stable translation length of their action on the curve graph. For a homeomorphism $\map{f}{S}{S}$, define its \emph{stable translation length} to be
\[
    \lambda(f)=\liminf_{k\to\infty}\frac{\dist[\C(S)](a,f^k(a))}{k},
\]
where $a$ is an arbitrary curve in $\C_0(S)$, the choice of which is clearly immaterial. It is easy to see that, if $f$ is periodic or reducible, then $\lambda(f)=0$. The converse is also true, as proved by \textcite[Proposition 4.6]{masur-minsky-1}. In fact, they show that if $f$ is pseudo\=/Anosov then $\lambda(f)\ge\constant{min translation length}$ for some constant $\constant{min translation length}>0$ depending only on $S$. Later, \textcite{gadre-min-translation-length} compute an explicit lower bound for $\constant{min translation length}$, which we take as an equality for the purpose of this article:
\[
    \declareconstant{min translation length}=\constantvalue{min translation length}.
\]
We remark that the result by \citeauthor{gadre-min-translation-length} is only stated for surfaces with $\xi(S)\ge 2$. However, Theorem 18 in \cite{baik-torus-translation-length} implies that the stable translation lengths of pseudo\=/Anosov homeomorphisms of the torus with one puncture and of the sphere with $4$ punctures are integers. In particular, the lower bound for $\constant{min translation length}$ stated above holds for surfaces $S$ with $\xi(S)=1$ as well.

When $f$ is pseudo\=/Anosov, we have the following easy estimate.

\begin{proposition}\label{thm:pseudo-anosov distance bound}
Let $S$ be a surface with $\xi(S)\ge 1$, and let $\map{f}{S}{S}$ be a pseudo\=/Anosov homeomorphism. Then for every integer $h\ge 0$ and every curve $a\in\C_0(S)$ the following inequality holds:
\[
    \dist[\C(S)](a,f^h(a))\ge \constant{min translation length}\cdot h.
\]
\end{proposition}
\begin{proof}
Note that, for every integer $h_1\ge 0$, we have
\[
    \dist[\C(S)](a,f^{hh_1}(a))\le h_1\dist[\C(S)](a,f^h(a)).
\]
As a consequence, we find that
\[
    \constant{min translation length}\le\lambda(f)\le\liminf_{h_1\to\infty}\frac{\dist[\C(S)](a,f^{hh_1}(a))}{h\cdot h_1}\le\frac{\dist[\C(S)](a,f^h(a))}{h}.\qedhere
\]
\end{proof}

Our goal is now to estimate how quickly the ratio $\dist[\C(S)](a,f^k(a))/k$ converges to $0$ as $k\to\infty$ when $f$ is a non\=/periodic reducible homeomorphism. The following proposition is based on \cite[Proposition 4.4]{bell-webb-algorithms}, with some additional effort on our part to keep the constants polynomial in $\xi$.

\newconstant{reducible bound}{C_{\mathrm{red}}}{Cred}{2(\constant{bounded geodesic image}+2)\cdot\constant{min translation length}^{-1}\cdot\xi}
\begin{proposition}\label{thm:reducible distance bound}
Let $\map{f}{S}{S}$ be a non\=/periodic reducible homeomorphism. Then there exists an integer $1\le h_0\le\constant{reducible bound}$ such that for every integer $h\ge 0$ and every curve $a\in\C_0(S)$ the inequality
\[
    \dist[\C(S)](a,f^{hh_0}(a))\le\dist[\C(S)](a,f^{h_0}(a))+2
\]
holds. The constant $\constant{reducible bound}$ is given by
\[
    \declareconstant{reducible bound}=\constantvalue{reducible bound}.
\]
\end{proposition}
\begin{proof}
Let $\Sigma$ be the \emph{canonical reduction system} of $f$ as defined by \textcite[\S2]{handel-new-nielsen} (there, the homeomorphism is called $\tau$ and the reduction system $\Gamma$). More explicitly, let $\Sigma'$ be the set of all vertices $c\in\C_0(S)$ whose orbit under $f$ is finite, and let $\Sigma$ be the union of all the curves in $\Sigma'$ that can be realised disjointly from all other curves in $\Sigma'$. Lemma 2.2 of \cite{handel-new-nielsen} shows that $\Sigma$ is non\=/empty, and moreover satisfies the following property: for every component $Z$ of $S\setminus\Sigma$ and every integer $k\ge 1$ such that $f^k(Z)=Z$, the homeomorphism $\map{f^k|_Z}{Z}{Z}$ is periodic or pseudo\=/Anosov.

Note that $\Sigma$ has at most $\xi$ components; let $c$ be one of them. There is an integer $1\le h_1\le 2\xi$ such that the homeomorphism $f_1=f^{h_1}$ preserves $c$ and the two (not necessarily distinct) components $Z_1$ and $Z_2$ of $S\setminus\Sigma$ which are adjacent to $c$. Let $f_2=f_1^{h_2}$, where $h_2$ is a positive integer that will be defined shortly; its value will only depend on the quantities introduced so far. Let $a\in\C_0(S)$ be any curve, and fix a geodesic $\gamma\subs\C(S)$ between $a$ and $f_2(a)$. We show that $\dist[\C(S)](\gamma,c)\le 1$ by analysing two cases.
\begin{substeps}
\item Suppose that the restriction of $f_1$ to (say) $Z_1$ is pseudo\=/Anosov; in particular, this entails that $\xi(Z_1)\ge 1$. Let $h_2=\lfloor\constant{bounded geodesic image}\cdot\constant{min translation length}^{-1}\rfloor+1$. Suppose, for a contradiction, that $\dist[\C(S)](\gamma,c)\ge 2$. Then every vertex of $\gamma$ must intersect $Z_1$, and \cref{thm:bounded geodesic image} gives the bound
\[
    \diam[\C(Z_1)](\pi_{Z_1}(\gamma\cap\C_0(S)))\le\constant{bounded geodesic image}.
\]

Fix now a curve $\bar{a}\in\pi_{Z_1}(a)$. We can apply \cref{thm:pseudo-anosov distance bound} to the surface $Z_1$ and the homeomorphism $f_1|_{Z_1}$ and deduce the chain of inequalities
\begin{align*}
\constant{min translation length}(S)\cdot h_2&\le\constant{min translation length}(Z_1)\cdot h_2\\
&\le\dist[\C(Z_1)](\bar{a},f_2(\bar{a}))\\
&\le\diam[\C(Z_1)](\pi_{Z_1}(a),f_2(\pi_{Z_1}(a)))\\
&=\diam[\C(Z_1)](\pi_{Z_1}(a),\pi_{Z_1}(f_2(a)))\\
&\le\diam[\C(Z_1)](\pi_{Z_1}(\gamma\cap\C_0(S)))\le\constant{bounded geodesic image},
\end{align*}
from which we derive the contradiction $h_2\le\constant{bounded geodesic image}\cdot\constant{min translation length}^{-1}$.
\item Suppose now that the restrictions of $f_1$ to $Z_1$ and $Z_2$ are periodic. Let $Z=Z_1\cup Z_2\cup c$, and let $Y\subs Z$ be a small annular neighbourhood of $c$. There exists an integer
\[
    1\le h_3\le \constant{periodic order bound}(Z_1)\cdot\constant{periodic order bound}(Z_2)\le\constant{periodic order bound}(S)^2
\]
such that $f_1^{h_3}|_{Z_1}$ and $f_1^{h_3}|_{Z_2}$ are isotopic to the identity. Note that $f_1^{h_3}|_Z$ is isotopic to $\twist_c^k$ on $Z$ for some $k\in\ZZ$; the integer $k$ is non\=/zero, otherwise $f_1|_Z$ would be periodic and $c$ would not be part of the canonical reduction system $\Sigma$.

Let $h_2=(\constant{bounded geodesic image}+2)\cdot h_3$, and suppose for a contradiction that $\dist[\C(S)](\gamma,c)\ge 2$. Then every vertex of $\gamma$ intersects $Y$, and again by \cref{thm:bounded geodesic image} we get the bound
\[
    \diam[\C(Y)](\pi_Y(\gamma\cap\C_0(S)))\le\constant{bounded geodesic image}.
\]
On the other hand, since $c$ is clearly essential in $Z$, \cref{thm:dehn twist annular projection} applied to $f_2$ implies that
\[
    \dist[\C(Y)](\pi_Y(a),\pi_Y(f_2(a)))\ge|(\constant{bounded geodesic image}+2)\cdot k|-1\ge\constant{bounded geodesic image}+1,
\]
contradicting \cref{thm:bounded geodesic image}.
\end{substeps}

In both cases, we have shown the existence of an integer
\[
    1\le h_2\le(\constant{bounded geodesic image}+2)\cdot\constant{min translation length}^{-1}
\]
such that $\dist[\C(S)](\gamma,c)\le 1$; let $b\in\gamma\cap\C_0(S)$ be at distance at most $1$ from $c$. We set $h_0=h_1\cdot h_2$, noting that the bound $h_0\le\constant{reducible bound}$ is satisfied. If now $h\ge 1$ is an arbitrary integer, we have the following chain of inequalities:
\begin{align*}
\dist[\C(S)](a,f_2^h(a))&\le\dist[\C(S)](a,c)+\dist[\C(S)](c,f_2^h(a))\\
&=\dist[\C(S)](a,c)+\dist[\C(S)](c,f_2(a))\\
&\le\dist[\C(S)](a,b)+\dist[\C(S)](b,f_2(a))+2\dist[\C(S)](b,c)\\
&\le\dist[\C(S)](a,f_2(a))+2.
\end{align*}
The equality on the second line follows from the fact that $c$ is fixed by $f_2$; the last inequality holds because $b$ lies on a geodesic between $a$ and $f_2(a)$, and $\dist[\C(S)](b,c)\le 1$. This concludes the proof.
\end{proof}

\subsection{Algorithmic classification}

We are now ready to present our algorithm for the Nielsen\=/Thurston classification of surface homeomorphisms.

\begin{theorem}[Nielsen\=/Thurston classification algorithm]\label{thm:nielsen-thurston classification algorithm}
Let $\map{f}{S}{S}$ be a homeomorphism.
There is an algorithm to decide the Nielsen\=/Thurston type of $f$.
The running time of the algorithm is polynomial in $\xi$ and $\norm{f}$.
\end{theorem}
\begin{proof}
\Cref{thm:periodic recognition} gives a procedure to decide whether $f$ is periodic.
Assuming it is not, we now need to decide if $f$ is reducible and non\=/periodic or pseudo\=/Anosov.
To this aim, we first find a curve $a\in\C_0(S)$ such that $\TTT(a)_e\le 2$ for every edge $e$ of $\TTT$; for instance, we can take $a$ to be a component of a neighbourhood of an edge of $\TTT$.
In particular, note that $\norm{a}_1\le\log(-6\chi+6)$.
For every integer $1\le h_0\le\constant{reducible bound}$, we perform the following operations.
\begin{itemize}
\item Compute
\[
    h=\left\lfloor\frac{\constant{carried curves l*}\cdot(2\norm{f^{h_0}(a)}_1+\constant{curves l+}+6)+\constant{curves l+}}{\constant{min translation length}\cdot h_0}\right\rfloor+1.
\]
\item Apply \cref{thm:distance algorithm for short curves} to the curves $a$ and $f^{hh_0}(a)$ to obtain an integer $d\ge 0$ such that
\[
d-\constant{curves l+}\le\dist[\C(S)](a,f^{hh_0}(a))\le\constant{carried curves l*}\cdot d+\constant{curves l+}.
\]
\item Check whether the inequality
\begin{equation}\label{eqn:nielsen-thurston classification algorithm:inequality}
d\le2\norm{f^{h_0}(a)}_1+\constant{curves l+}+6
\end{equation}
is satisfied.
\end{itemize}

We claim that $f$ is reducible if and only if \cref{eqn:nielsen-thurston classification algorithm:inequality} is satisfied for at least one value of $1\le h_0\le\constant{reducible bound}$.
Since the procedure described above can be performed in polynomial time in $\xi$ and $\norm{f}$, it is enough to prove this claim.
\begin{substeps}
\item If $f$ is reducible (and non\=/periodic), then \cref{thm:reducible distance bound} guarantees the existence of an integer $1\le h_0\le\constant{reducible bound}$ such that
\begin{align*}
d-\constant{curves l+}&\le\dist[\C(S)](a,f^{hh_0}(a))\\
&\le\dist[\C(S)](a,f^{h_0}(a))+2\\
&\le2\log i(a,f^{h_0}(a))+4\\
&\le2\norm{f^{h_0}(a)}_1+6,
\end{align*}
which obviously implies that \cref{eqn:nielsen-thurston classification algorithm:inequality} holds for this value of $h_0$.
\item If $f$ is pseudo\=/Anosov and $h_0$ is any positive integer, then \cref{thm:pseudo-anosov distance bound} gives
\begin{align*}
    \constant{carried curves l*}\cdot d+\constant{curves l+}&\ge\dist[\C(S)](a,f^{hh_0}(a))\\
    &\ge\constant{min translation length}\cdot h_0\cdot h\\
    &>\constant{carried curves l*}\cdot(2\norm{f^{h_0}(a)}_1+\constant{curves l+}+6)+\constant{curves l+},
\end{align*}
which contradicts \cref{eqn:nielsen-thurston classification algorithm:inequality}. \qedhere
\end{substeps}
\end{proof}

\begin{remark}
In the statement of \cref{thm:nielsen-thurston classification algorithm}, we give no estimate for the degree of the polynomial dependence of the running time on $\norm{f}$.
This omission stems from three main reasons: firstly, since we are not imposing any specific representation of surface homeomorphisms, we have no control over how expensive it is to compute images of curves under $f$; secondly, we make use of \cref{thm:distance algorithm for short curves}, for which -- as pointed out in \cref{rmk:distance algorithm polynomial degree} -- we do not give an explicit polynomial dependence in the closed case; thirdly, the proof of \cref{thm:periodic recognition} for closed surfaces relies on the algorithm of \cite[Theorem 1.2]{lackenby-minimal-position}, for which \citeauthor{lackenby-minimal-position} does not give an explicit polynomial bound on the running time.

However, with slight loss of generality, we can give a more precise upper bound to the running time of the algorithm of \cref{thm:nielsen-thurston classification algorithm}, assuming that $S$ is punctured, and $f$ is given as a sequence of $\norm{f}$ flips (and no Dehn twists) in the flip and twist graph.
In this setting, it is easy to see that, given an essential curve $c$ on $S$, we can compute $f(c)$ in time $\bigO(\poly(\xi)\cdot\norm{f}\cdot\norm{c})$; moreover, the complexity of $f(c)$ satisfies the bound
\[
    \norm{f(c)}_1\le\norm{f}+\norm{c}_1.
\]
It follows that the integer $h$ used in the proof of \cref{thm:nielsen-thurston classification algorithm} satisfies $h=\bigO(\poly(\xi)\cdot\norm{f})$, and hence
\[
    \norm{f^{hh_0}(a)}_1=\bigO(\poly(\xi)\cdot\norm{f}^2).
\]
Using the bound given in \cref{rmk:distance algorithm polynomial degree} for the running time of \cref{thm:distance algorithm for short curves} for punctured surfaces -- and noting that the algorithm of \cref{thm:periodic recognition} runs in linear time in $\norm{f}$ when $S$ is punctured -- we conclude that, in this restricted setting, the Nielsen\=/Thurston classification algorithm has running time
\[
    \bigO(\poly(\xi)\cdot\norm{f}^4\cdot\log\norm{f}).\qedhere
\]
\end{remark}

%% file: figures-source/periodic-recognition-A_2.tex
\begin{tikzpicture}

\foreach \i in {1,2,3,4} {
    \foreach \j/\name in {1/l,-1/r,.5/ll,-.5/rr} {
        \coordinate (\i-\name) at ($({100+120*\i}:1.5)+({190+120*\i}:{\j*.5})$);
    }
}
\foreach \i/\j in {1/2,2/3,3/1} {
    \path[spath/save global=\i] (\i-l) to[bend left=30] (\j-r);
}
\foreach \i/\j in {1/2,2/3,3/1} {
    \coordinate (a) at (spath cs:{\i} .9);
    \coordinate (b) at (spath cs:{\j} .1);
    \draw[line,palette 2,spath/save global=a-\i] (a) to[bend left=95,looseness=.6] (b);
    \draw[line,opacity=.3,palette 2,dotted] (b) to[bend left=85,looseness=.6] (a);
}
\foreach \i/\j in {1/2,2/3,3/1} {
    \draw[black,line,opacity=.3] [spath/use=\i];
    \draw[line,palette 3,opacity=.3,spath/save global=b-\i] (\i-ll) to[bend left=30,looseness=1.1] (\j-rr);
}
\foreach \i/\j in {1/3,2/1,3/2} {
    \tikzset{spath/split at intersections with={\i}{a-\j},spath/get components of={\i}{\components},spath/split at intersections with={\getComponentOf{\components}{2}}{a-\i},spath/get components of={\getComponentOf{\components}{2}}{\components}}
    \draw[line,black] [spath/use=\getComponentOf{\components}{1}];
    \tikzset{spath/split at intersections with={b-\i}{a-\j},spath/get components of={b-\i}{\components},spath/split at intersections with={\getComponentOf{\components}{2}}{a-\i},spath/get components of={\getComponentOf{\components}{2}}{\components}}
    \draw[line,palette 3] [spath/use=\getComponentOf{\components}{1}];
    \foreach \k in {0,1} {
        \fill[black] (spath cs:\getComponentOf{\components}{1} \k) circle(.75pt);
    }
}

\node[black] at (0,0) {$P$};
\node[label node,\colourmakedark{2},above=6pt] at (spath cs:a-1 1) {$\mathcal{A}_1$};
\draw[pin line] (label) -- (spath cs:a-1 .8);
\node[label node,\colourmakedark{3},left=6pt] at (spath cs:3 .5) {$A_2$};
\draw[pin line] (label) -- (spath cs:b-3 .55);

\end{tikzpicture}

%% file: figures-source/annular-subsurface-projection-a.tex
\begin{tikzpicture}
\tikzsetupannulussubsurfaceprojection{}

\annuluscontourbehind

\cbarbehind

\otherliftsofc
\cbar
\annuluscontour

\node[\colourmakedark{2},label node,above=4pt,xshift=-4pt] at (spath cs:top edge .5) {$\bar{c}$};
\draw[pin line] (label) -- (spath cs:c bar .3);

\node[black,label node,below=4pt,xshift=4pt] at (spath cs:bottom edge .5) {$\bar{A}$};

\end{tikzpicture}

%% file: figures-source/annular-subsurface-projection-b.tex
\begin{tikzpicture}

\tikzsetupannulussubsurfaceprojection{}

\Ybarbackground

\annuluscontourbehind

\cbarbehind
\otherliftsofc
\cbar
\Ybarcontour
\annuluscontour

\node[\colourmakedark{3},label node,below=4pt,xshift=4pt] at (spath cs:bottom edge .5) {$\bar{Y}$};
\draw[pin line] (label) -- (.05,-.3);

\end{tikzpicture}

%% file: figures-source/annular-subsurface-projection-c.tex
\begin{tikzpicture}

\tikzsetupannulussubsurfaceprojection{}

\Zbarbackground

\annuluscontourbehind

\cbarbehind

\Zbarcontour
\otherliftsofc
\cbar
\annuluscontour

\node[\colourmakedark{4},label node,below=4pt,xshift=4pt] at (spath cs:bottom edge .5) {$\bar{Z}$};
\draw[pin line] (label) -- (.3,-.3);

\end{tikzpicture}

%% file: figures-source/annular-subsurface-projection-d.tex
\begin{tikzpicture}

\tikzsetupannulussubsurfaceprojection{}

\Zbarbackground

\annuluscontourbehind

\cbarbehind
\Wbehind

\Zbarcontour
\otherliftsofc
\cbar
\W
\annuluscontour

\node[\colourmakedark{1},label node,below=4pt] at (spath cs:bottom edge .5) {$W$};
\draw[pin line] (label) -- (spath cs:Y bar contour left .7);
\draw[pin line] (label) -- (spath cs:Y bar contour right .7);

\end{tikzpicture}

%% file: figures-source/annular-subsurface-projection-e.tex
\begin{tikzpicture}

\tikzsetupannulussubsurfaceprojection{}

\annuluscontourbehind

\abarbehind
\Wbehind
\fabarwigglybehind

\W
\abar
\fabarwiggly
\annuluscontour

\node[\colourmakedark{1},label node,below=4pt] at (spath cs:bottom edge .5) {$W$};
\draw[pin line] (label) -- (spath cs:Y bar contour left .7);
\draw[pin line] (label) -- (spath cs:Y bar contour right .7);
\node[\colourmakedark{2},label node,left=4pt,yshift=-12pt] at (spath cs:a bar 0) {$\bar{a}$};
\draw[pin line] (label) -- (spath cs:a bar .07);
\node[\colourmakedark{3},label node,left=4pt,yshift=-2pt] at (spath cs:f a bar 0) {$\bar{f}(\bar{a})$};
\draw[pin line] (label) -- (f a bar wiggly pin);

\end{tikzpicture}

%% file: figures-source/annular-subsurface-projection-f.tex
\begin{tikzpicture}

\tikzsetupannulussubsurfaceprojection{}

\annuluscontourbehind

\abarbehind
\Wbehind
\fabarbehind
\W
\abar
\fabar
\annuluscontour

\intersectionsofabarandfabar

\end{tikzpicture}